\documentclass[3p,preprint,review,11pt]{elsarticle}
\usepackage{amssymb}
\usepackage{amsmath}
\usepackage{amsfonts}
\usepackage{mathrsfs}
\usepackage{physics}
\usepackage{url}
\usepackage{hyperref}
\usepackage{xcolor}
\usepackage[T1]{fontenc}
\newdefinition{rmk}{Remark}
\journal{Applied Mathematics and Computation}
\date{}
\begin{document}

\begin{frontmatter}



\title{Performance of affine-splitting pseudo-spectral methods for fractional complex Ginzburg-Landau equations}

\author[FCEIA]{Lisandro A. Raviola\corref{cor1}}
\cortext[cor1]{Corresponding author.}
\ead{raviola@fceia.unr.edu.ar}
%\affiliation[FCEIA]{organization={Departamento de F{\'i}sica y Qu{\'i}mica}}
\affiliation[FCEIA]{organization={Departamento de F{\'i}sica y Qu{\'i}mica, Escuela de Formaci{\'o}n B{\'a}sica, Facultad de Ciencias Exactas, Ingenier{\'i}a y Agrimensura (FCEIA), Universidad Nacional de Rosario}, addressline={Av. Pellegrini 250},city={Rosario}, postcode={S2000BTP}, state={Santa Fe}, country={Argentina}}
\author[INMABB]{Mariano F. De Leo}
\ead{mariano.deleo@uns.edu.ar}
\affiliation[INMABB]{organization={Instituto de Matem{\'a}tica de Bah{\'i}a Blanca (INMABB--CONICET), Departamento de Matem{\'a}tica, Universidad Nacional del Sur},%Department and Organization
addressline={Av. Alem 1253}, 
city={Bah{\'i}a Blanca},
postcode={B8000CPB}, 
state={Buenos Aires},
country={Argentina}}

\begin{abstract}
We evaluate the performance of novel numerical methods for solving one-dimensional nonlinear fractional dispersive and dissipative evolution equations. The methods are based on affine combinations of time-splitting integrators and pseudo-spectral discretizations using Hermite and Fourier expansions. We show the effectiveness of the proposed methods by numerically computing the dynamics of soliton solutions of the the standard and fractional variants of the nonlinear Schr{\"o}dinger equation (NLSE) and the complex Ginzburg-Landau equation (CGLE), and by comparing the results with those obtained by standard splitting integrators. An exhaustive numerical investigation shows that the new technique is competitive when compared to traditional composition-splitting schemes for the case of Hamiltonian problems both in terms accuracy and computational cost. Moreover, it is applicable straightforwardly to irreversible models, outperforming high-order symplectic integrators which could become unstable due to their need of negative time steps. Finally, we discuss potential improvements of the numerical methods aimed to increase their efficiency, and  possible applications to the investigation of dissipative solitons that arise in nonlinear optical systems of contemporary interest. Overall, the method offers a promising alternative for solving a wide range of evolutionary partial differential equations.
\end{abstract}


%\item We introduce a novel affine time-splitting pseudospectral numerical method for solving nonlinear fractional dispersive PDEs in unbounded unidimensional domains.
%\item We showcase the method by studying the nonlinear space-fractional Schrödinger equation and the fractional complex Ginzburg-Landau equation
%\item We ...

\end{frontmatter}

\section{Introduction}
Nonlinear partial differential equations play an essential role as mathematical devices to describe and analyze complex phenomena that depend continuously on space and time. In the last decades, many scientific and technological problems have prompted the investigation of relevant nonlinear models, like the nonlinear Schr{\"o}dinger equation (NLSE) for the study of quasi-monochromatic light propagation in nonlinear optical media and the dynamics of water and plasma waves \cite{sulem1999nls, scott2004encyclopedia, ablowitz2011nonlinearwaves, agrawal2013nonlinearfiber}; the closely related Gross-Pitaevskii equation (GPE) for the quantum description of Bose-Einstein condensation \cite{scott2004encyclopedia, Bao2003GPE, bao2005splitting}; and the complex Ginzburg-Landau equation (CGLE) for non-equilibrium phenomena in hydrodynamics, chemical reactions and dissipative optics \cite{aranson2002review, garcia-morales2012CGLE}. While rigorous theoretical results for various instances of these models have been established on solid bases, and many of its experimental consequences verified in the laboratory, numerical methods continue to play a crucial role by bridging the gap between theory and experiments, facilitating the development of analytical results and suggesting experimental pathways through computer simulations. Among the various numerical techniques employed, time-splitting methods \cite{mclahan2002splitting, hairer2006geometric, blanes2008splitting, holden2010splitting, borgna2015splitting, glowinski2017splitting, blanes2017concise} combined with pseudo-spectral space discretizations \cite{fornberg1996pseudospectral, trefethen2000spectral, boyd2000chebyshev, canuto2007spectral, hesthaven2007spectral, kopriva2009implementing, shen2011spectral} have been widely utilized to solve many cases of the aforementioned models due to their versatility and remarkable accuracy. For instance, the split-step Fourier method was successfully applied to the NLSE in \cite{taha1984analytical, weideman1986split, pathria1990pseudo, muslu2005higher}, while the GPE has been solved using second- and fourth-order split-step Fourier and Hermite pseudo-spectral schemes by Bao \emph{et al.} \cite{Bao2003GPE, bao2005splitting}, to mention only some of the pioneering works. The geometric preserving properties of the time-splitting schemes are particularly useful due to the existence of conserved quantities in these models. In this regard, Antoine \emph{et al.} \cite{antoine2013computational} have done an extensive review of the various numerical methods, devoting particular attention to the conservation of dynamical invariants at the discrete level.

On the other hand, complex Ginzburg-Landau equations (CGLE) with cubic and cubic-quintic nonlinearities --which can be considered generalizations of the NLSE that account for linear and nonlinear dissipative effects-- were investigated both theoretically and numerically by Akhmediev \emph{et al.} \cite{akhmediev1996singularities, akhmediev2001pulsating} in the context of nonlinear dissipative optics phenomena. By means of extensive numerical simulations based on low-order split-step Fourier schemes they identified and characterized many kinds of stable localized solutions called \emph{dissipative solitons} \cite{akhmediev2001solitons, akhmediev2005dissipative, akhmediev2008dissipative, ferreira2022dissipative}, anticipating the experimental observations made years later in mode-locked lasers \cite{grelu2012dissipative}. 

More recently, there has been a growing interest in generalizations of the previous models by incorporating non-local effects, in particular fractional dispersion and diffusion in space \cite{laskin2000fractional, weitzner2003fractional,  guo2012standingfNLS, frank2013fractional, klein2014fractional, longhi2015fractional, duo2016spectralfNLS, mao2017fractional, qiu2020solitonfCGLE, liu2023experimental}. The investigation of these fractional generalizations is an ongoing endeavor and many open questions remain to be answered, therefore the construction of robust and dependable numerical schemes to strengthen this research line is of utmost relevance.

While standard time-splitting methods possess many desirable properties, namely preservation of geometric invariants, those of order greater than two unavoidably require negative time steps \cite{goldman1996nonreversible, blanes2005necessity}, a fact that makes them unsuitable in irreversible problems due to possible unstable behavior, for the backward evolution is ill-posed. In a recent paper, De Leo \emph{et al.} \cite{deleo2016affine} introduced a novel numerical technique based on affine combinations of Lie-Trotter integrators which retain the benefits of the splitting strategy while using only positive time steps and thus allowing for the application of high-order schemes to irreversible equations. Their paper is mainly devoted to prove rigorously the convergence order of the proposed schemes, which the authors claim to be well-suited for the solution of irreversible pseudo-differential problems, the fractional CGLE (fCGLE) being a paradigmatic example. Although they suggest that the new schemes' performance could be comparable to standard splitting methods without the associated drawbacks, this claim is not systematically explored. Therefore, the present article aims to thoroughly evaluate the efficiency and accuracy of the high-order affine-splitting schemes introduced in \cite{deleo2016affine} and assess their effectiveness for solving the NLSE and CGLE in both the standard and fractional cases. As a byproduct of our investigation, we also obtain new numerical results on the fCGLE that have not been reported elsewhere.

The structure of the article is the following: in Section \ref{sec:methods} we introduce the problems that will be addressed and the methods used to numerically solve them, while in Section \ref{sec:results} we show and discuss the obtained results for the investigated models. Section \ref{sec:conclusions} summarizes the main conclusions and perspectives for further developments and applications. In \ref{sec:pseudo-spectral} we briefly describe the pseudo-spectral discretization, and in  \ref{sec:propagators} we include some analytical results and additional details for the numerical implementation.

\section{Numerical methods}
\label{sec:methods}
%\subsection{Model problem}
In this paper, we will be concerned with the numerical solution of one-dimensional initial value problems (IVP) of the form
\begin{equation}
\left\{ \quad
\begin{aligned}
    \label{eq:IVP}
    \mathrm{i} \partial_t \ u(x,t) &= A u(x,t) + B(u(x,t)),\qquad x \in  \mathbb{R}, t \geq 0, \\ 
    u(x,0) &= u_0(x),
\end{aligned}\right.
\end{equation}
where $A$ is a linear {\it Fourier multiplier} (a pseudo-differential operator) with symbol $\mathcal{A}(k) \in \mathbb{C}$,
\begin{equation}
    \widehat{A u}(k,t) = \mathcal{A}(k)\  \hat{u}(k,t),
    \label{eq:A_symbol}
\end{equation}
defined by means of the unitary Fourier Transform (FT) in $x$  
\begin{equation}
    \label{eq:unitaryFT}
    \hat{u}(k, t) := \mathcal{F} \big\{ u (x,t)\big\} =  \frac{1}{\sqrt{2 \pi}}\int_\mathbb{R} u(x,t)\  \mathrm{e}^{-\mathrm{i} k x} \  \mathrm{d}x,
\end{equation}
and $B$ is a possibly nonlinear operator acting on $u$. The general conditions that $u_0$, $A$ and $B$ must fulfill in order for the numerical methods to be applicable are detailed in reference \cite{deleo2016affine}, but suffice to say here that the framework is sufficiently general to encompass a vast array of evolution equations including diffusion, reaction-diffusion, linear and nonlinear Schrödinger, Gross-Pitaevskii and complex Ginzburg-Landau models. For example, we get the one-dimensional Gross-Pitaevskii equation with a potential $V(x)$ when $\mathcal{A}(k)=\frac{1}{2}|k|^2$ and $B(u)=V(x) u -\vert  u \vert^2 u$. In this work we consider $u_0 \in L^2(\mathbb{R})$.

By virtue of the possibly non-local character of linear operator $A$, defined in frequency space by equation \eqref{eq:A_symbol}, we resort to pseudo-spectral methods for the space discretization of states and operators. These methods approximate the solution globally by means of an expansion in orthogonal functions that satisfy appropriate boundary conditions. Even though their applicability is generally restricted to problems in simple geometries, they are nonetheless well suited for the problems addressed in this paper, standing out for having excellent convergence properties –especially when the solutions are highly regular– and for requiring less computational resources than finite difference methods for achieving a given numerical precision \cite{fornberg1996pseudospectral, trefethen2000spectral, boyd2000chebyshev, canuto2007spectral, hesthaven2007spectral, kopriva2009implementing, shen2011spectral}. On the other hand, we approximate the evolution in time by means of time-splitting methods. In what follows, we briefly summarize the various  splitting schemes used in this work, particularly those introduced in \cite{deleo2016affine}, which will form the basis of the proposed numerical methods.  The details of the pseudo-spectral discretization are deferred to \ref{sec:pseudo-spectral} for the interested reader. Note that while alternative discretizations and higher dimensional problems can be addressed along with the presented splitting schemes, in this paper we will stay within the context of one-dimensional problems to highlight the basic structure of numerical methods in their most elementary form.  Results concerning the generalizations mentioned above will be presented elsewhere.

\subsection*{Evolution in time: operator-splitting methods}
Operator-splitting methods (also called \emph{time-splitting}, \emph{split-step} or \emph {fractional-step} methods) have a long history in the solution of IVP for ordinary and partial differential equations, with origins in problems of celestial mechanics, molecular dynamics, particle accelerators and fluid mechanics. They are an instance of the old {\em divide-and-conquer} strategy and rely on the {\em splitting} of the IVP \eqref{eq:IVP} into two partial initial value subproblems
\begin{align}
 \mathrm{i} \partial_t u(x,t) &= A u(x,t),  \quad u(x,0) = u_0(x), \label{eq:sub_A}\\ 
 \mathrm{i} \partial_t u(x,t) &= B (u(x,t)), \quad u(x,0) = u_0(x), \label{eq:sub_B}
\end{align}
which are in some sense easier or more convenient to solve than the complete problem, allowing for the application of a specialized solver for each subproblem to exploit and/or preserve some aspect of its particular structure (e.g. qualitative properties). The splitting criteria can stem from a physical (different phenomena), dimensional (different coordinates) and/or mathematical (linear-nonlinear separation) rationale. Once obtained, either in exact or approximate form, the partial solutions are appropriately combined to give an approximate solution of the complete problem. Building upon this basic idea numerous schemes had been constructed, and there exists an ample bibliography, so we give here only a brief account of their main features \cite{mclahan2002splitting, hairer2006geometric, blanes2008splitting, holden2010splitting, borgna2015splitting, glowinski2017splitting, blanes2017concise}. 

If $\phi(t)$ is the flow associated with the IVP \eqref{eq:IVP}, namely if $u(x,t)=\phi(t)(u_0(x))$ is its solution, a first order approximation to $\phi$ by means of operator splitting is given by the {\em Lie--Trotter} propagator
\begin{equation} 
    \Phi_\mathrm{LT}(\Delta t) = \phi_B(\Delta t) \circ \phi_A(\Delta t) = \phi(\Delta t) + \mathcal{O}(\Delta t^2),
\end{equation}
where $\phi_A$ and $\phi_B$ are the flows associated with the partial subproblems \eqref{eq:sub_A} and \eqref{eq:sub_B}, respectively. Another popular method is the {\em Strang-Marchuk} second-order scheme
\begin{align}
      \Phi_\mathrm{SM}(\Delta t)  &= \phi_A(\frac{1}{2}\Delta t)\circ\phi_B(\Delta t)\circ \phi_A(\frac{1}{2}\Delta t) = \phi(\Delta t) + \mathcal{O}(\Delta t^3).
\end{align}
In general, similar composition-based splittings can be cast in the form
\begin{equation}
    \Phi_\mathrm{Symp}(\Delta t) = \phi_B(b_s \Delta t) \circ \phi_A(a_s \Delta t)\circ \cdots \circ \phi_B(b_1 \Delta t)\circ \phi_A(a_1 \Delta t),
    \label{eq:composition_splitting}
\end{equation}
where $s$ is the number of stages, $a_1+\cdots+a_s=b_1+\cdots +b_s=1$, and the coefficients are chosen to achieve the desired order of convergence $q$. Within this framework, a third-order scheme was proposed by Ruth \cite{ruth1983splitting} with coefficients $$b_3=\frac{7}{24}, a_3=\frac{2}{3}, b_2=\frac{3}{4}, a_2=-\frac{2}{3}, b_1=-\frac{1}{24}, a_1=1,$$ 
and Neri \cite{neri1987lie} devised an efficient fourth-order scheme using 
\begin{equation}
a_1=a_4=\frac{1}{2(2-2^{1/3})},a_2=a_3=-\frac{2^{1/3}-1}{2(2-2^{1/3})}, b_1=b_3=2 a_1, b_2=-2^{1/3} b_1, b_4=0.
\label{eq:neri_coeffs}
\end{equation}
Even-order composition methods can be obtained in a systematic fashion by means of the construction technique introduced by Yoshida in \cite{yoshida1990symplectic}, where in particular an eight-stages sixth-order ($s=8,q=6$) method is reported (which we call Yoshida6 in what follows). A valuable practical resource for these and other composition splittings is \cite{auzinger2023splittings}, where coefficients and references of many schemes can be found.

It is well known that schemes of the form \eqref{eq:composition_splitting} require some of the coefficients $a_i,b_i$ to be negative when the order is greater than two \cite{goldman1996nonreversible, blanes2005necessity}, so their stability is not guaranteed in irreversible problems which are ill-posed when integrated backwards in time. In order to overcome this drawback while still retaining the versatility and accuracy of high-order splitting schemes, De Leo \emph{et al.} \cite{deleo2016affine} have introduced a new type of additive schemes based on affine combinations of Lie-Trotter propagators that avoid negative steps, called \emph{affine time-splitting methods}. The construction of affine schemes proceed as follows: given the partial propagators $\phi_A, \phi_B$, define the pair of adjoint {\em Lie--Trotter} propagators as
\begin{align}
    \Phi^+_1(\Delta t) &= \phi_B(\Delta t) \circ \phi_A(\Delta t), \label{eq:prop_plus} \\
    \Phi^-_1(\Delta t) &= \phi_A(\Delta t) \circ \phi_B(\Delta t), \label{eq:prop_minus}
\end{align}
and recursively construct the compositions
\begin{equation}
    \Phi_m^\pm(\Delta t) = \Phi^\pm_1(\Delta t) \circ \Phi^\pm_{m-1}(\Delta t), \qquad m = 2, 3, \dots
\end{equation}
Then, the following  {\em symmetric} affine combinations are introduced, involving only \emph{positive} time steps: \footnote{Reference \cite{deleo2016affine} also introduces {\em asymmetric} combinations. Given that symmetric methods outperform the asymmetric ones, in this paper we use exclusively the former.}
\begin{equation}
    \Phi_\mathrm{SA}^s(\Delta t) = \sum_{j=1}^s \gamma_j \Big(\Phi_j^+\big(\frac{\Delta t}{j}\big) + \Phi_j^-\big(\frac{\Delta t}{j}\big)\Big). \label{eq:sym_affine_combination}
\end{equation}
It is proved rigorously in \cite{deleo2016affine} that, under appropriate hypothesis on the operators and initial data, the flow defined by \eqref{eq:sym_affine_combination} converges with order $q=2n$ to the exact flow $\phi(\Delta t)$, if and only if the set of coefficients $\gamma_j$ satisfy the consistency and order conditions
\begin{align}
    \sum_{j=1}^s \gamma_j &=  \frac{1}{2}, \\
    \sum_{j=1}^s \frac{\gamma_j}{j^{2k}} &= 0, \qquad 1 \leq k \leq n-1.
    \label{eq:sym_order_cond}
\end{align}
As discussed in the reference, the method is highly parallelizable, with the most consuming (last) stage requiring the calculation of $2s$ Lie-Trotter evolutions. Besides, as pointed out in \cite{stillfjord2018adaptive}, high order methods embed trivially (with different coefficients $\gamma_j$) the lower order ones, thus allowing for the cheap computation of local error estimates, which in turn can be used to control the step size for the construction of adaptive schemes.

In the present paper we give compelling numerical evidence showing that high-order affine methods, although not designed to preserve the underlying geometric structure of the models, exhibit a competitive performance when compared with composition splittings in Hamiltonian problems. Moreover, they are generally superior in terms of accuracy, stability and computational cost in problems with irreversible dynamics. We show that this is the case even without additional optimizations (namely parallelization and time-adaptivity), a subject we will discuss in another paper. We concentrate mainly in schemes of order 2, 4, and 6, which are those of most practical value for the intended applications. 

\section{Example models and numerical results}
\label{sec:results}
In this section we compare the properties and performance of numerical schemes based on composition (Strang, Neri and Yoshida6) and symmetric affine splittings of order 2, 4 and 6, using Hermite and Fourier pseudo-spectral discretizations of states and operators (see \ref{sec:pseudo-spectral} for details). With this aim, we apply each scheme to approximate the solution of the IVP given in \eqref{eq:IVP} for several instances of the standard and fractional NLSE and CGLE. To assess quantitatively the relative merit of each method, we define various metrics. In the first place, given a reference solution $u_\mathrm{ref}(x,t)$ of the equation (known either exactly or numerically with high accuracy), we compute with the chosen method the approximate numerical solution $u_\mathrm{num}(x_n, t)$ on the spatial grid $x_n$ associated with the pseudo-spectral discretization and obtain its absolute error with respect to the reference solution in the discrete $L^\infty_x$-norm, i.e. at time $t$ we compute the maximum absolute error on the grid
\begin{equation}
    \mathcal{E}_{\infty}(t)= \max_{\{x_n\}} \big\lvert u_\mathrm{ref}(x_n,t) - u_\mathrm{num}(x_n,t) \big\rvert.
\end{equation}
If the model possesses dynamical invariants, we also assess the quality of the numerical solutions by computing the relative errors of the conserved quantities. This is particularly relevant since we want to compare the performance of the proposed affine schemes with respect to geometric integrators explicitly designed to preserve some of these quantities. To this end, we define the relative error in the computation of the conserved quantity $Q$ at time $t$  as
\begin{equation}
    \epsilon_Q(t) = \Big|\frac{Q(t)}{Q(0)}-1\Big|.
\end{equation}
To measure the efficiency of each scheme with respect to the previous metrics, we determine the {\em computational cost} required to achieve a given error. This cost is estimated by counting the total number of evaluations of the propagators $\phi_A$ and $\phi_B$, which serve as the fundamental computational units for all the splitting schemes. It is worth noticing that symmetric affine schemes evaluate both propagators an equal number of times at each step, but this is not the case for the even-order composition schemes analyzed in this article. For instance, the fourth-order scheme of Neri given in \eqref{eq:neri_coeffs} requires three and four evaluations of each propagator per step, as one of its coefficients is zero. Additionally, the cost of evaluating each propagator can differ significantly (for example, $\phi_A $ may be implemented as an $\mathcal{O}(N^2)$ matrix-vector product or an $\mathcal{O}(N \log_2 N)$ Fast Fourier Transform, while $\phi_B$ may require an $\mathcal{O}(N)$ element-wise vector-vector multiplication). Therefore, the performance of composition schemes can be influenced by the sequence in which the propagators are computed. To conduct a fair comparison between the two families of splitting schemes, it would be necessary to determine the total computing time per step for each approach. However, since this time can be highly dependent on factors such as implementation details, computer architecture, and software environment, we consider more fruitful to estimate the computational cost as the number of evaluations of the most expensive propagator. For composition schemes, this propagator is chosen as the {\em least} invoked one, which establishes a lower bound on their computational cost, bearing in mind that this approach may underestimate the actual cost for these schemes.

All numerical experiments were performed using both Hermite and Fourier pseudo-spectral discretizations. When high accuracy is sought, the Fourier method achieves an absolute error $\mathcal{E}_{\infty}$ lower than the Hermite method by virtue of its lower computational complexity (and thus less round-off error), but in all cases the results are qualitatively similar and lead to identical general conclusions. A brief account of the pseudo-spectral method is given in \ref{sec:pseudo-spectral}. In  \ref{sec:propagators} we give details on the calculation of partial propagators for the various models. 

The algorithms were implemented in Python 3.10.9, utilizing the numerical libraries \texttt{numpy} 1.22.4 \cite{harris2020numpy} and \texttt{scipy} 1.7.3 \cite{virtanen2020scipy}, and the experiments were conducted on a laptop computer with an Intel Core i5-10210U processor and 16 GB of RAM. The source code needed for generating the presented results is available at the repository \url{https://github.com/raviola/pseudosplit-paper} under a free software license \cite{raviola2023pseudosplit}.

\subsection{The nonlinear Schrödinger equation}
\subsubsection{Cubic nonlinear Schrödinger equation (NLSE3)}
The one-dimensional cubic nonlinear Schrödinger equation (NLSE3) is a well-known Hamiltonian model and thus constitutes an ideal first benchmark for the proposed schemes. The NLSE3 governs the evolution of the envelope of slowly varying quasi-monochromatic wave packets in weakly nonlinear media with dispersion and negligible dissipation, and arises in the context of nonlinear optics, hydrodynamics and plasma waves phenomena \cite{sulem1999nls, scott2004encyclopedia, ablowitz2011nonlinearwaves, agrawal2013nonlinearfiber}. It is usually given in the standard form 
\begin{equation}
    \mathrm{i} \partial_t u = \frac{1}{2}(-\partial_{x}^2) u \pm |u|^{2} u,
    \label{eq:NLSE3}
\end{equation}
where the minus sign corresponds to the focusing or attracting case. Its complete integrability, obtained by the Inverse Scattering Transform (IST), implies the existence of an infinite number of conserved quantities, among them \cite{ablowitz2011nonlinearwaves} 
\begin{enumerate} 
    \item the \emph{mass}, \emph{particle number} or \emph{energy} (depending on the physical context)
    \begin{equation}
        M(u) = \int_\mathbb{R} |u|^2 \mathrm{d}x;
        \label{eq:NLSE3_mass}
    \end{equation}

    \item the {\em Hamiltonian} or \emph{energy}
    \begin{equation}
        H(u) = \frac{1}{2}  \int_\mathbb{R} \left( |\partial_x u|^2  \pm |u|^{4} \right) \mathrm{d}x.
        \label{eq:NLSE3_hamiltonian}
    \end{equation}
\end{enumerate}
In the focusing case, the NLSE3 admits localized pulse-like waves or {\em soliton} solutions \cite{sulem1999nls, scott2004encyclopedia, ablowitz2011nonlinearwaves}
\begin{equation}
    u^\mathrm{NLS}_\mathrm{sol}(x,t) = \eta \  \mathrm{sech}\left[\eta (x - ct - x_0)\right] \mathrm{e}^{\mathrm{i}(cx-\omega_\mathrm{sol}t+\phi_0)}
    \label{eq:NLSE3_soliton}
\end{equation}
where $\omega_\mathrm{sol}=(c^2 - \eta^2)/2$, $\eta$ and $c$ are the arbitrary amplitude and speed of the soliton, and $x_0,\phi_0$ are the soliton's center position and phase at $t=0$.

\begin{figure}[t]
    \centering
    \includegraphics[width=0.45\textwidth]{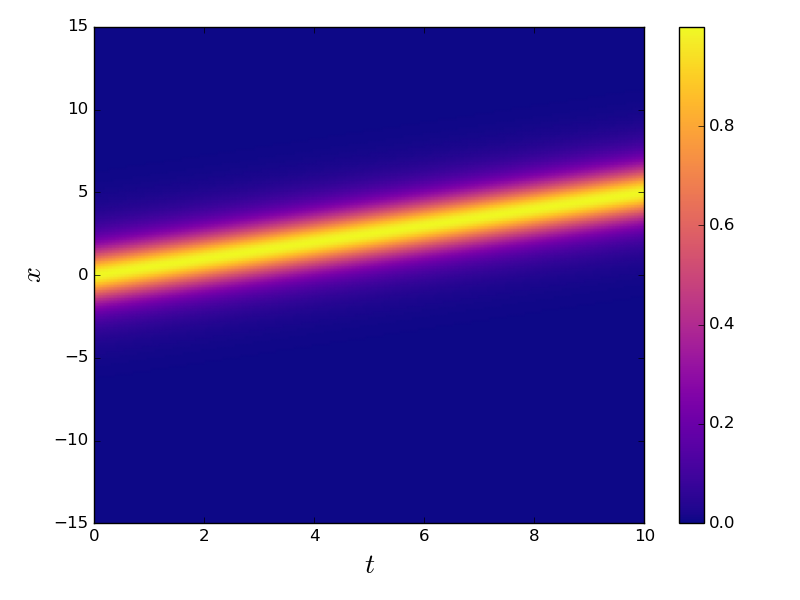} \includegraphics[width=0.54\textwidth]{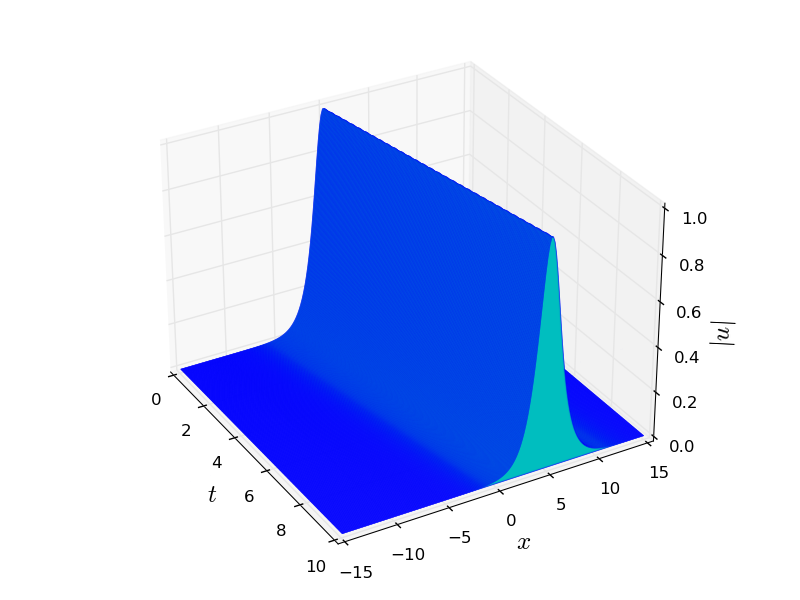}
    \caption{(Left) Colormap for the absolute value of the numerical solution of the NLSE3 in  $t\in[0,10]$ with initial value $u_0(x)=\sech(x)\mathrm{e}^{\mathrm{i}\frac{1}{2}x}$, corresponding to a soliton with parameters  $\eta=1, c=0.5, \phi_0=x_0=0$.  (Right) Soliton evolution up to $t=10$.}
    \label{fig:NLSE3_abs_solution}
\end{figure}

For the numerical example of this subsection, we chose as reference solution $u_\mathrm{ref}(x,t)$ the traveling soliton given by equation \eqref{eq:NLSE3_soliton} with parameters $\eta=1, \phi_0=x_0=0$, moving with speed $c=0.5$. The discretization in space is obtained by using the Fourier pseudo-spectral method with $N=2^{11}$ modes on the interval $I=[-50, 50]\subset \mathbb{R}$ (see \ref{sec:pseudo-spectral} for details on the pseudo-spectral discretization and criteria for parameters selection).  The evolution of the soliton is shown in Figure \ref{fig:NLSE3_abs_solution}, as calculated with the sixth-order affine splitting with $\Delta t=0.025$ for $t\in[0,10]$.  It is observed the expected result, namely, the soliton moves without changing its shape and its center position goes from $x_0=0$ to $x=5$ at $t=10$.

\begin{figure}[t]
    \centering
    \includegraphics[width=0.49\textwidth]{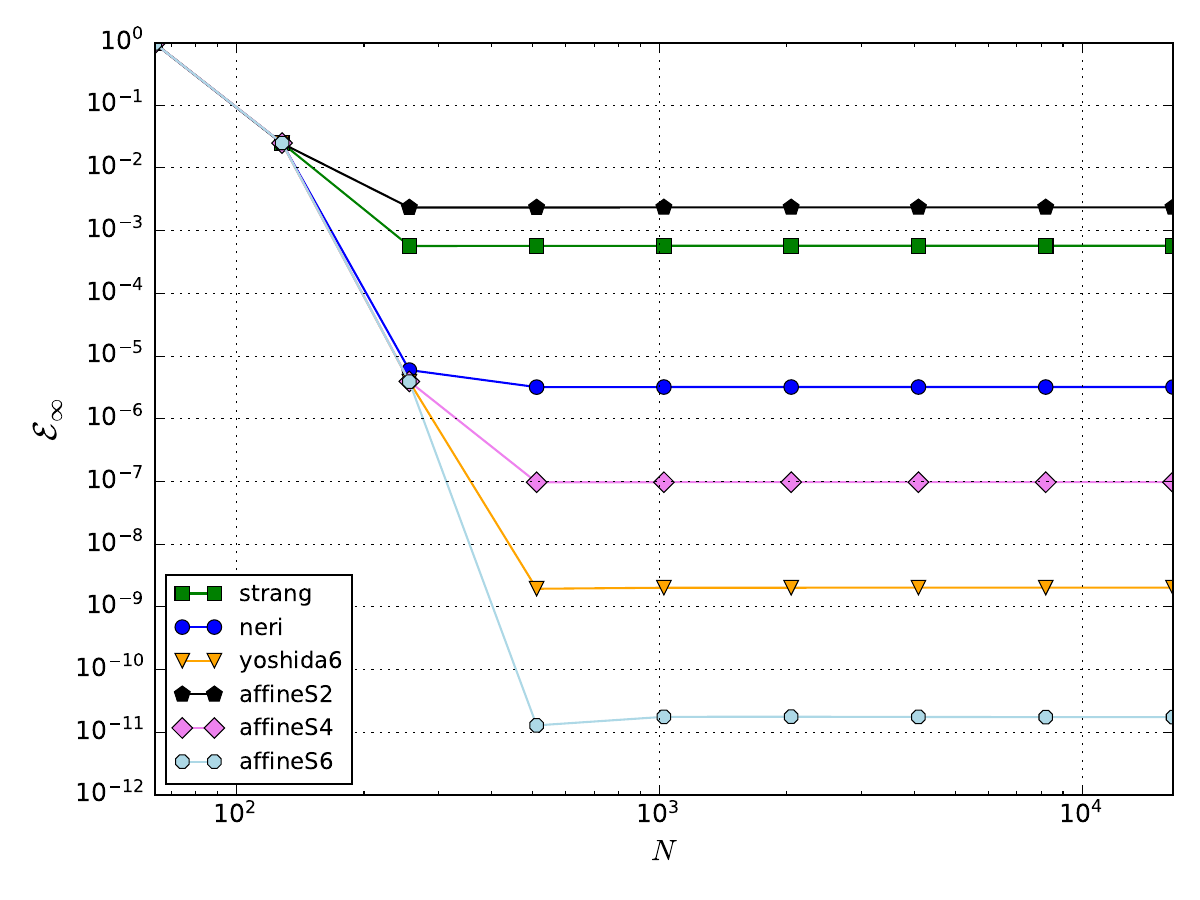}
    \includegraphics[width=0.49\textwidth]{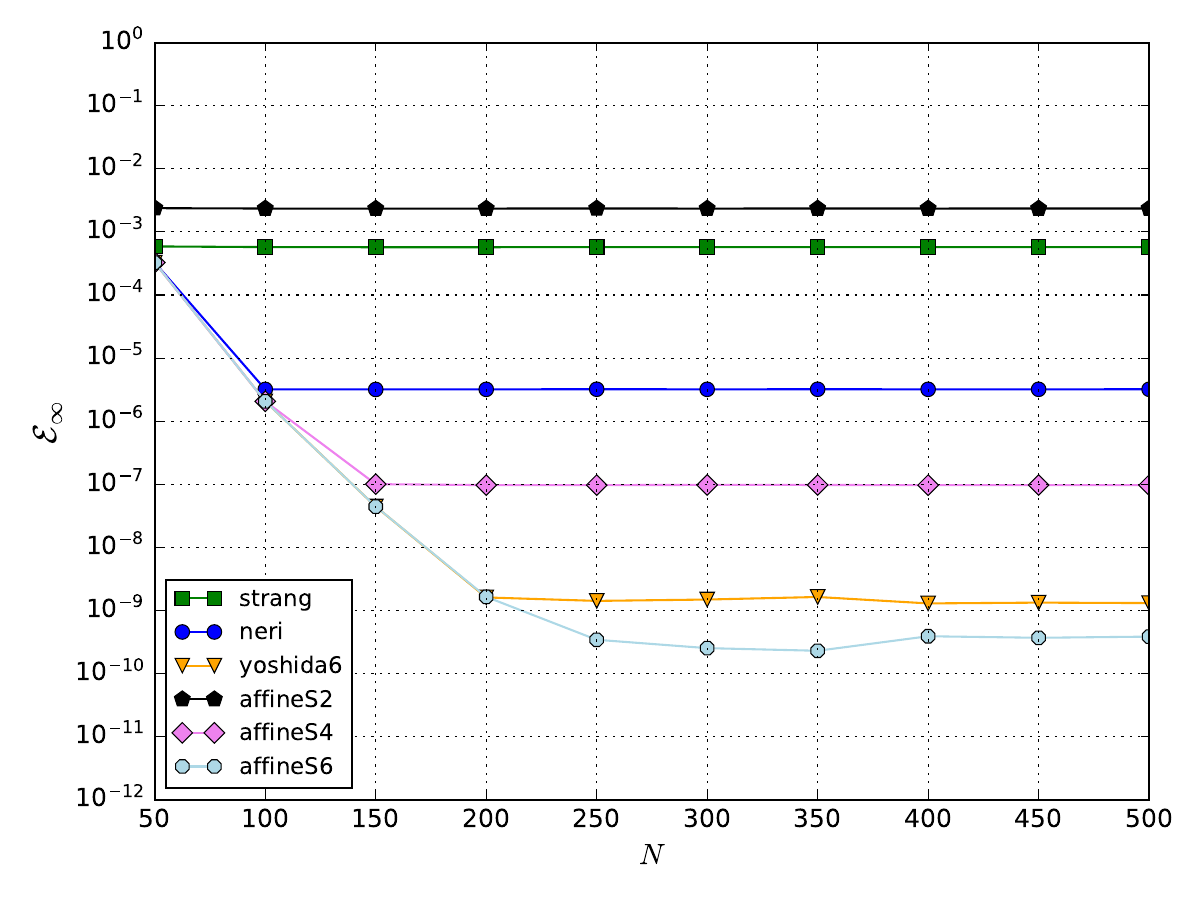}
    \caption{(Left) Absolute error $\mathcal{E}_\infty$ at time $t=10$, as a function of the of the Fourier pseudo-spectral basis dimension $N=2^m, m=6,\dots,14$ (in log scale), for different composition and symmetric affine time-splittings. The interval of the Fourier basis is $I=[-50,50]$. (Right) Absolute error as a function of the dimension of the Hermite basis (in linear scale). The scaling parameter of the basis is $s=1.25$. In both cases, the error is computed against an exact reference soliton \ref{eq:NLSE3_soliton} with parameters $\eta=1.0, c=0.5, \phi_0=x_0=0$, and the time step is $\Delta t=0.025$. } 
    \label{fig:NLSE3_fourier_soliton_err_inf_vs_N}
\end{figure}

The experiments that follow were conducted using both the Fourier and Hermite pseudo-spectral methods with similar results, but some minor differences between both discretizations are worth discussing. To this end, in Figure \ref{fig:NLSE3_fourier_soliton_err_inf_vs_N}  we compare the behavior of the pseudo-spectral approximations. The left panel of this figure shows the max-norm error of the numerical solution $\mathcal{E}_\infty$ at $t=10$  as a function of the Fourier basis dimension $N$, while the right panel gives the same information for the Hermite basis using a scaling factor $s=1.25$ (this last value was selected empirically to get the best results, see \cite{boyd2000chebyshev, tang1993hermite, ma2005scaling, tang2018fractional} for a discussion). Spectral convergence can be clearly appreciated for both methods (as expected owing to the infinite smoothness of the solution) up to the point where the error in time or round-off dominate the total error.  Comparison of left and right panels reveals that, for a given splitting scheme and time step $\Delta t$,  both discretizations yield very similar errors for a sufficiently large number $N$ of basis functions, unless the main source of error is round-off accumulation. This latter case is exemplified by the bottom (light blue) lines corresponding to the  highly accurate sixth-order affine splitting. It can be seen that the Fourier pseudo-spectral method yields a lower error than the Hermite method, and the difference can be credited to its lower computational cost, which is   $\mathcal{O}(N\log_{2}N)$ versus  $\mathcal{O}(N^{2})$. Indeed, the lowest error is attained by the Fourier method with $N=2^{10}\approx10^{3}$ when combined with the sixth-order affine method, while for the Hermite method this happens with $N\approx3\times10^2$. Then the cost of the Fourier method is $\mathcal{O}(10^{4})$ while for the Hermite case it is $\mathcal{O}(10^{5})$. In conclusion, when the error is due mainly to round-off it can be expected a greater accuracy with the Fourier method owing to its lower operations count.

In the rest of this subsection and  for the sake of brevity we will only show the results of the Fourier method, which will be described in certain detail, as a similar procedure will be used for analyzing all the following examples.
\begin{figure}[t]
    \centering
    \includegraphics[width=0.49\textwidth]{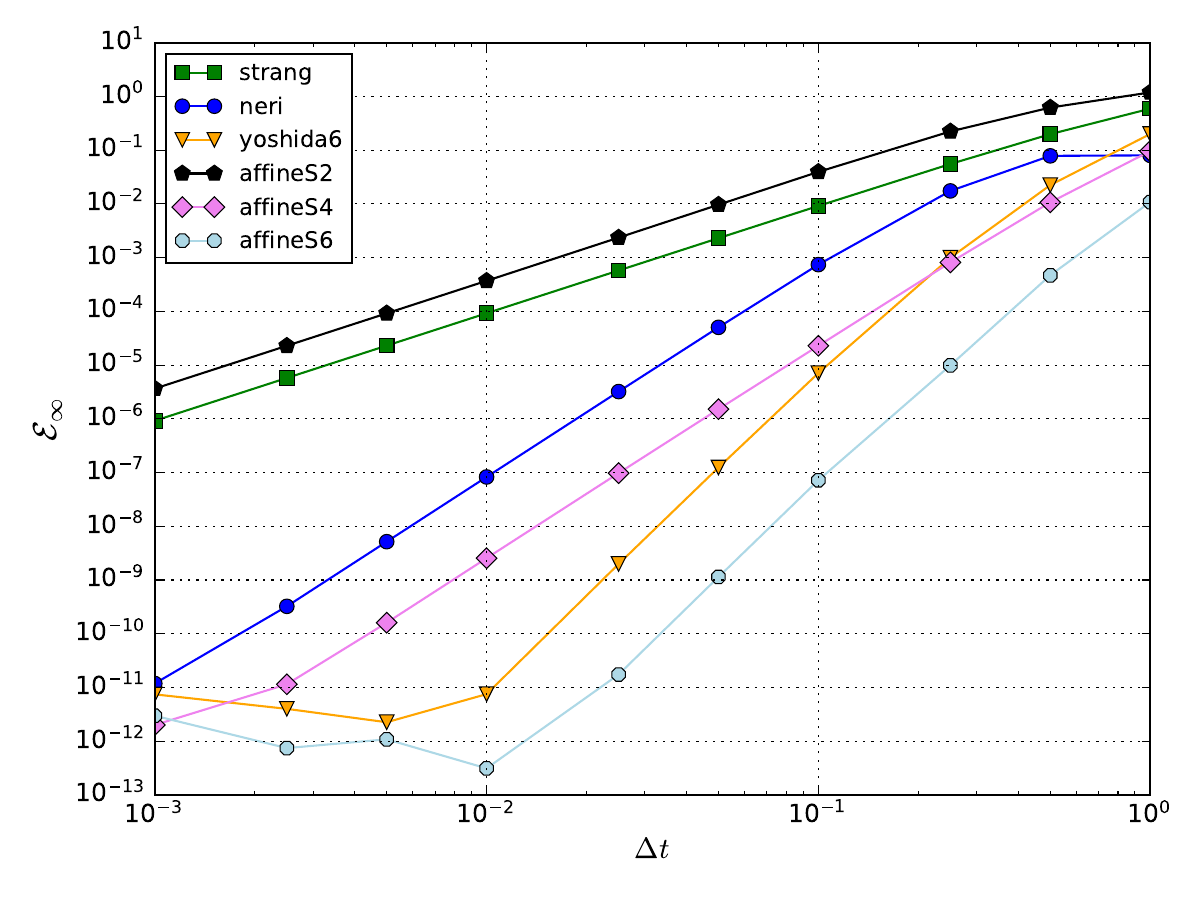}
    \includegraphics[width=0.49\textwidth]{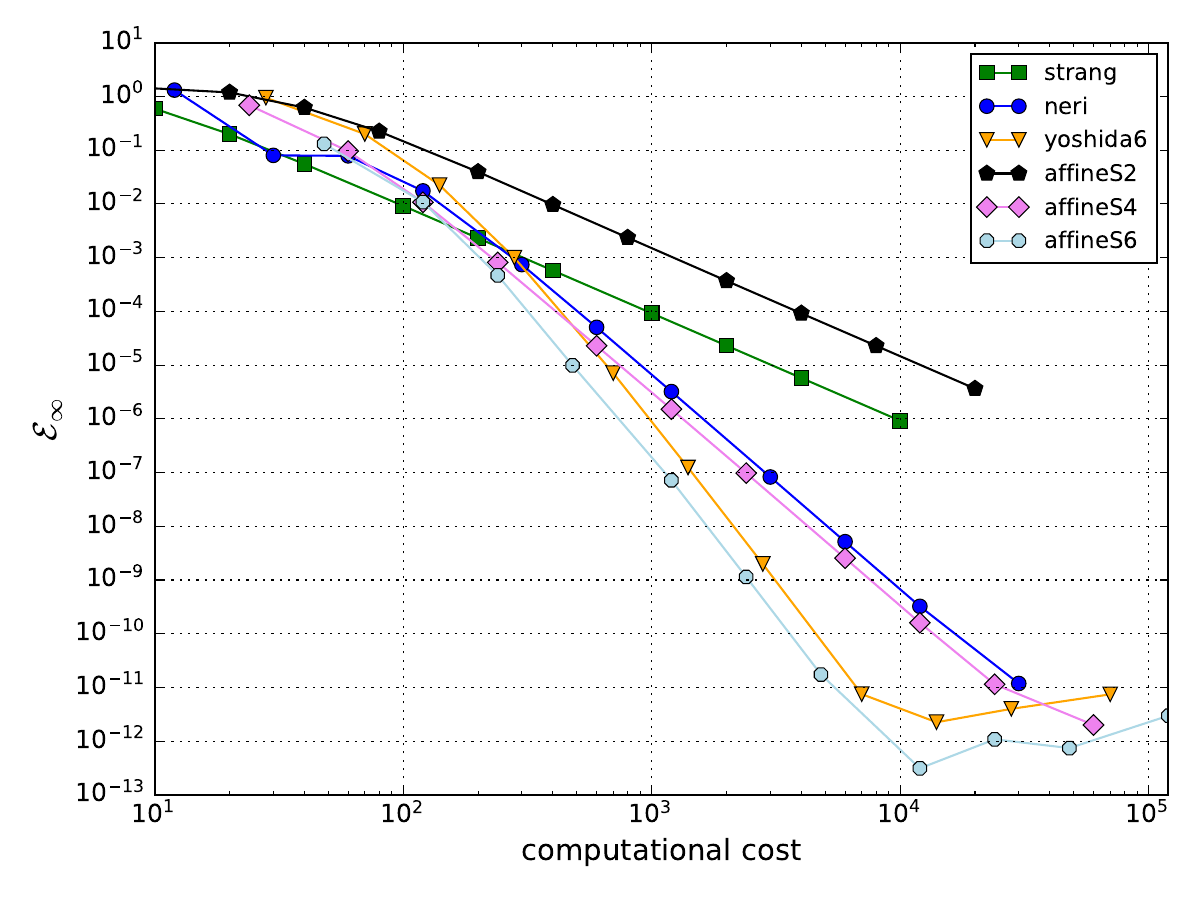}
    \caption{(Left) Absolute error $\mathcal{E}_\infty$ at time $t=10$ for the numerical solution of the focusing NLSE3, as a function of the time step $\Delta t$, for different composition and symmetric affine splittings. (Right) Absolute error $\mathcal{E}_\infty(t=10)$ vs. computational cost (efficiency plot) for the same splittings. The reference solution is the soliton \eqref{eq:NLSE3_soliton} with $c=0.5, x_0=\phi_0=0$ and $\eta=1$. States and operators are discretized by means of a Fourier pseudo-spectral method with $N=2^{11}$ modes over the interval $[-50, 50]\subset \mathbb{R}$.}
    \label{fig:NLSE3_soliton_fourier_err_inf}
\end{figure}

The left panel of Figure \ref{fig:NLSE3_soliton_fourier_err_inf} shows the numerical solution error $\mathcal{E}_\infty$ at $t=10$ as function of the time step $\Delta t$ for the investigated splitting schemes, in log-log scale. The actual convergence order $q$ of each scheme can be inferred from the slope of the corresponding line in the central region of the figure, and is consistent with the theoretical one.
It is clearly seen in this figure that, except for the second order scheme, affine schemes yield lower errors than the corresponding composition schemes of the same order for a given step size $\Delta t$. Notably, the difference in favor of the affine schemes can be of more than one order of magnitude. Moreover, the efficiency plot shown in the right panel of Figure \ref{fig:NLSE3_soliton_fourier_err_inf} reveals that the greater accuracy of the fourth and sixth-order affine schemes is also attained with lower computational cost than for composition schemes of the same order. When the computational cost of the high-order schemes reaches approximately $10^4$, the error then inverts its decreasing tendency due to accumulation of round-off errors, which begin to dominate the total error.

\begin{figure}[t]
    \centering
    \includegraphics[width=0.49\textwidth]{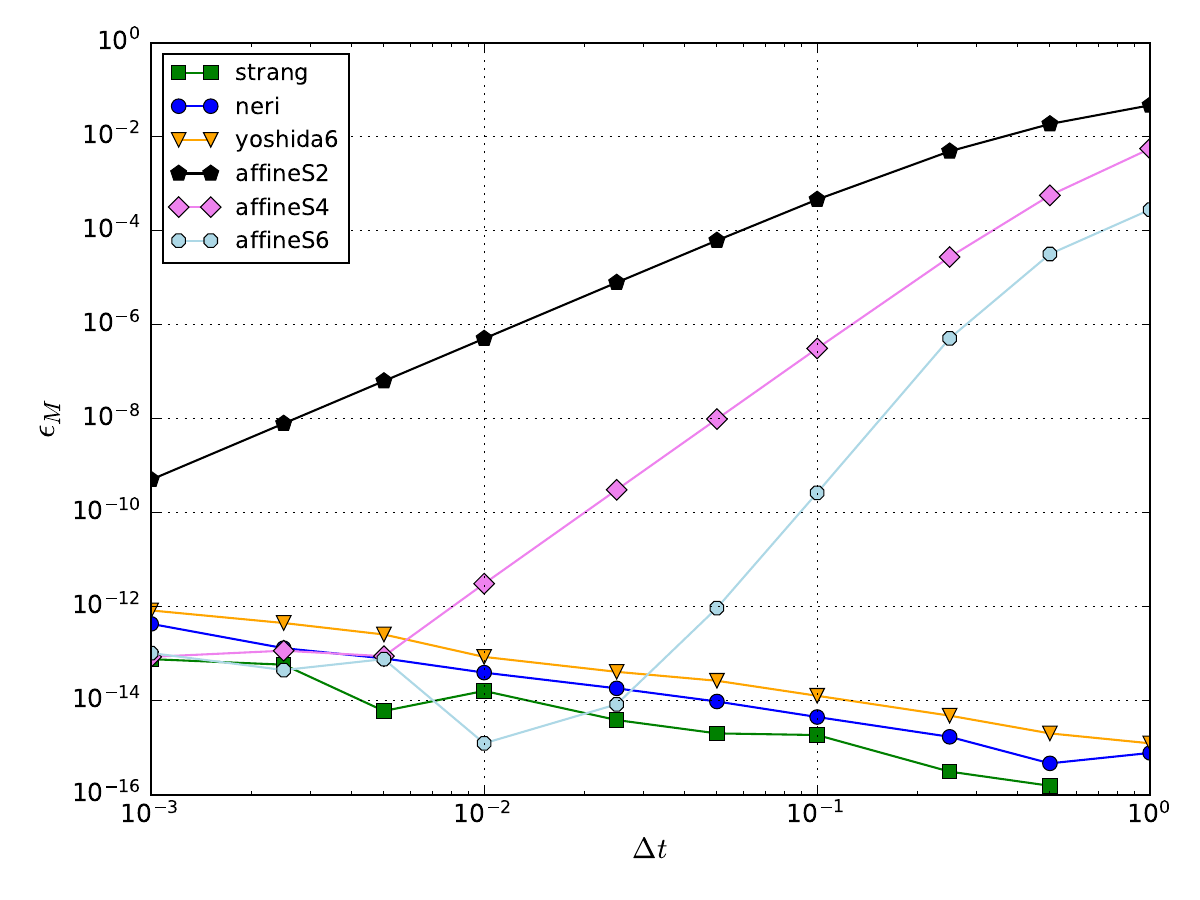}
    \includegraphics[width=0.49\textwidth]{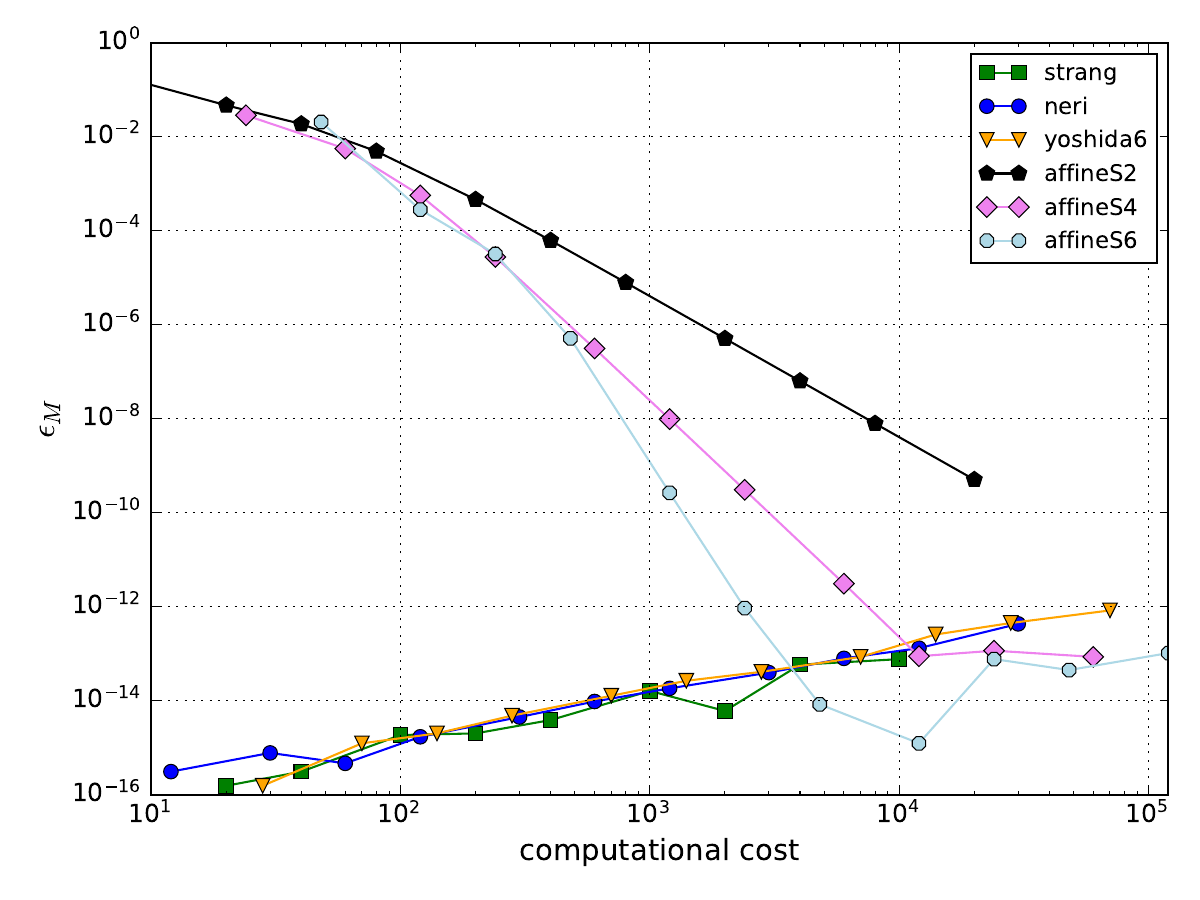}
    \caption{(Left) Mass relative error $\epsilon_M$ at time $t=10$ for the numerical solution of the NLSE3, as a function of the time step $\Delta t$, for different composition and symmetric affine splittings. (Right) Relative error $\epsilon_M(t=10)$ as a function of the computational cost. The parameters are those of Figure \ref{fig:NLSE3_soliton_fourier_err_inf}.} 
    \label{fig:NLSE3_fourier_soliton_err_M}
\end{figure}
Next we investigate the complementary aspect of the numerical preservation of conserved quantities, namely the mass \eqref{eq:NLSE3_mass} and Hamiltonian \eqref{eq:NLSE3_hamiltonian}. To this end, in the left panel of Figure \ref{fig:NLSE3_fourier_soliton_err_M} we plot the mass relative error $\epsilon_M$ at $t=10$, as a function of the step size $\Delta t$. In accordance with the mass preservation property of composition schemes,  constructed by composing propagators that separately preserve the mass (see the computation of propagators in  \ref{sec:propagators}), it can be seen that these methods set the lower bound for $\epsilon_M$ across the entire range of explored $\Delta t$, with the second-order Strang scheme reaching the lowest error as a consequence of the lower operations count per step, followed by the fourth and sixth-order composition schemes. It is also evident for the composition schemes the almost linear monotone growth of $\epsilon_M$ with diminishing step size, due to the accumulation of round-off errors as the total number of operations increases. On the other hand, and in spite of not being designed to enforce this preservation property, high-order affine schemes perform equivalently well (with errors near machine precision) if the step size is below a certain threshold (about $10^{-1}$ for the sixth-order scheme and $10^{-2}$ for the fourth-order scheme).

\begin{figure}[t]
    \centering
    \includegraphics[width=0.49\textwidth]{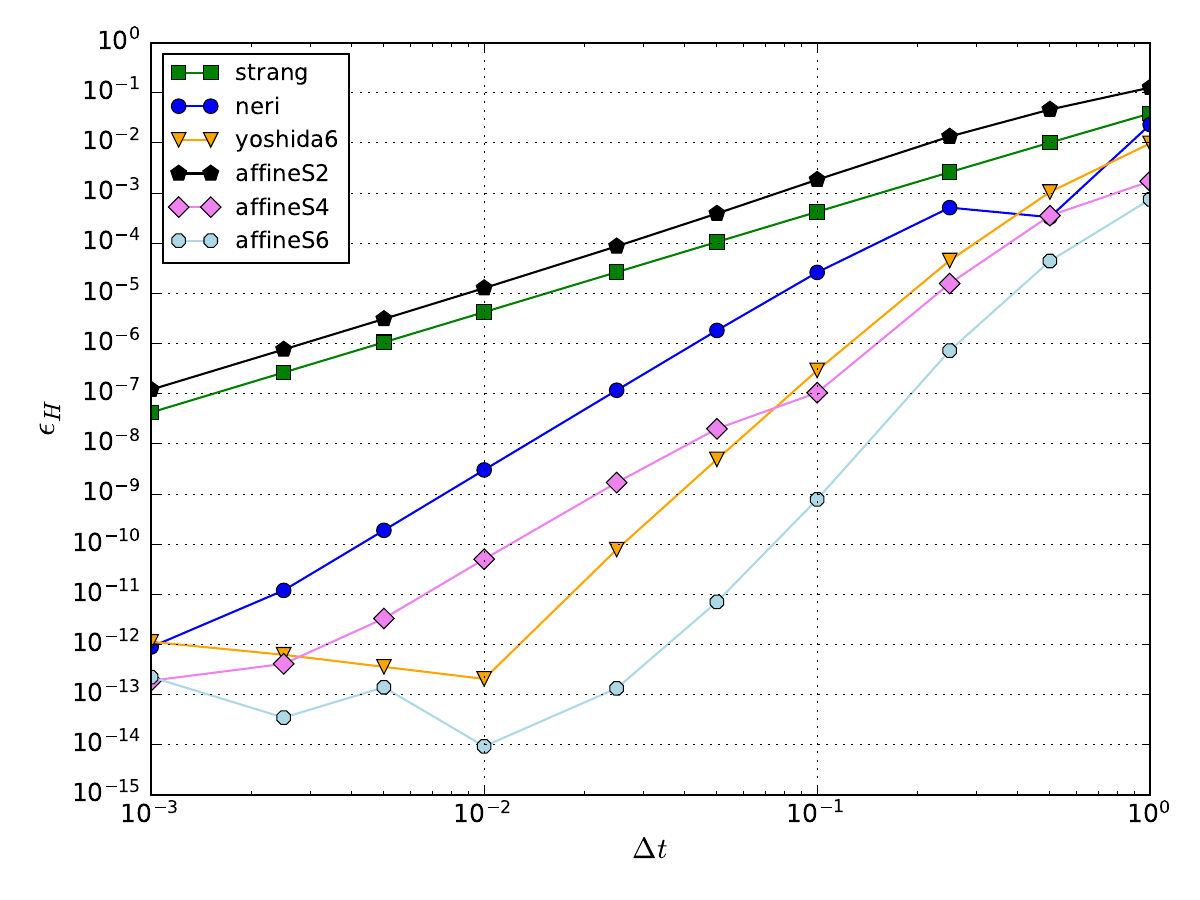}
    \includegraphics[width=0.49\textwidth]{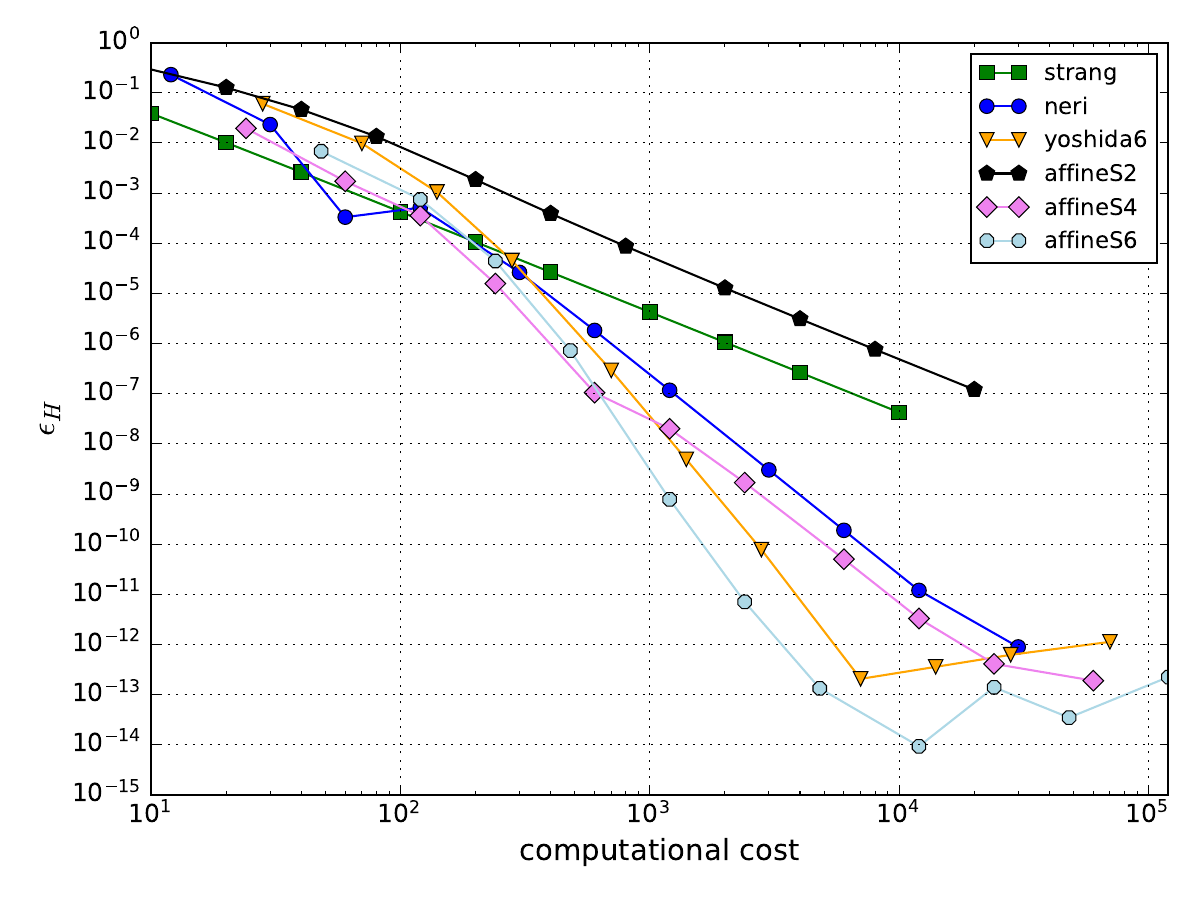}
    \caption{(Left) Relative error for the Hamiltonian $\epsilon_H$ at time $t=10$, as a function of the time step $\Delta t$, for different composition and symmetric affine splittings. (Right) Idem for the computational cost. The parameters are those of Figure \ref{fig:NLSE3_soliton_fourier_err_inf}.} 
    \label{fig:NLSE3_fourier_soliton_err_H}
\end{figure}

Regarding the Hamiltonian conservation, the left panel of Figure \ref{fig:NLSE3_fourier_soliton_err_H} shows the good behavior of affine schemes of order higher than two in terms of accuracy for all the explored step sizes. For a given step size, the difference in accuracy between different types of schemes of the same order can reach two orders of magnitude. Also, the computational cost required for a given accuracy is lower for affine schemes than for composition schemes of equivalent order, when the relative error is below $\epsilon_H\approx 10^{-3}$, as can be seen in the right panel of the same figure. The sixth-order affine scheme shows a particularly good overall performance.

\begin{figure}[t]
    \centering
    \includegraphics[width=0.49\textwidth]{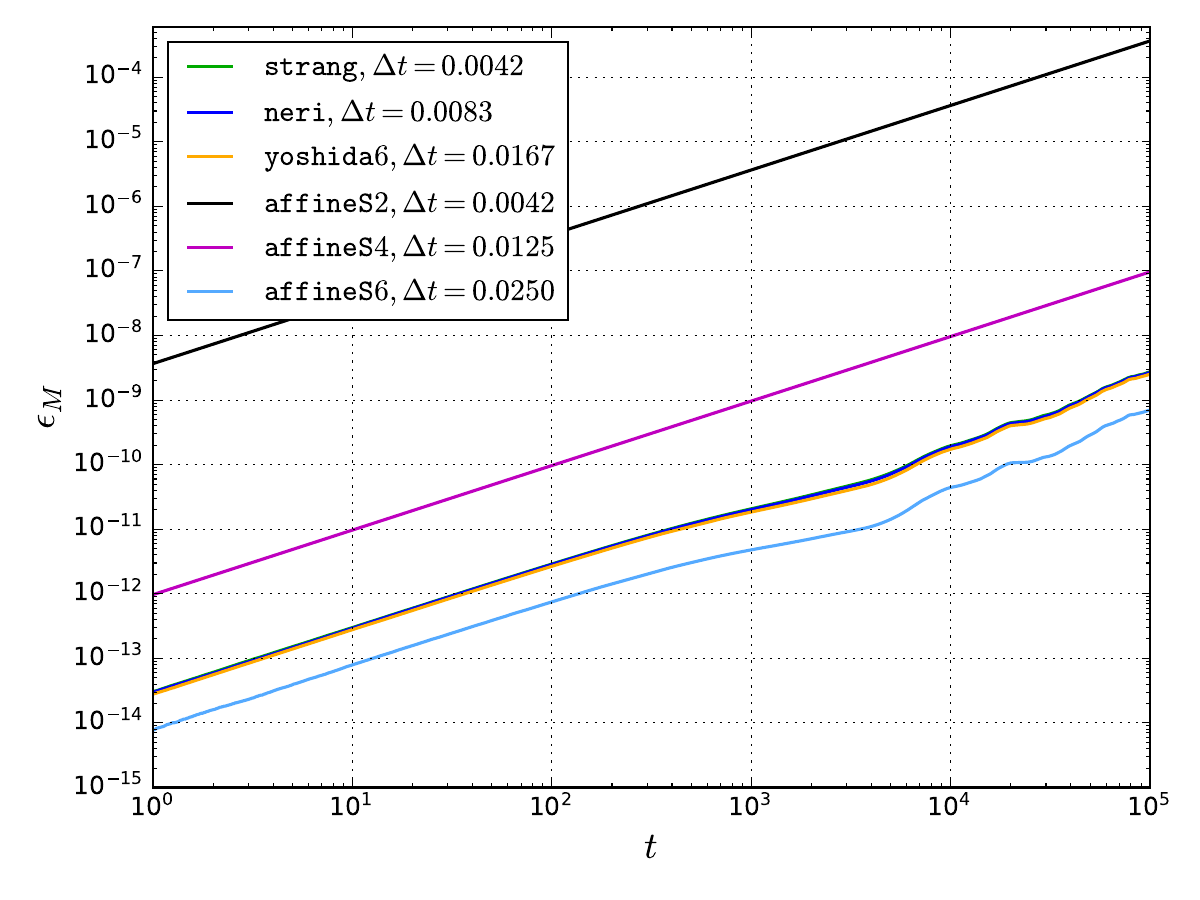}
    \includegraphics[width=0.49\textwidth]{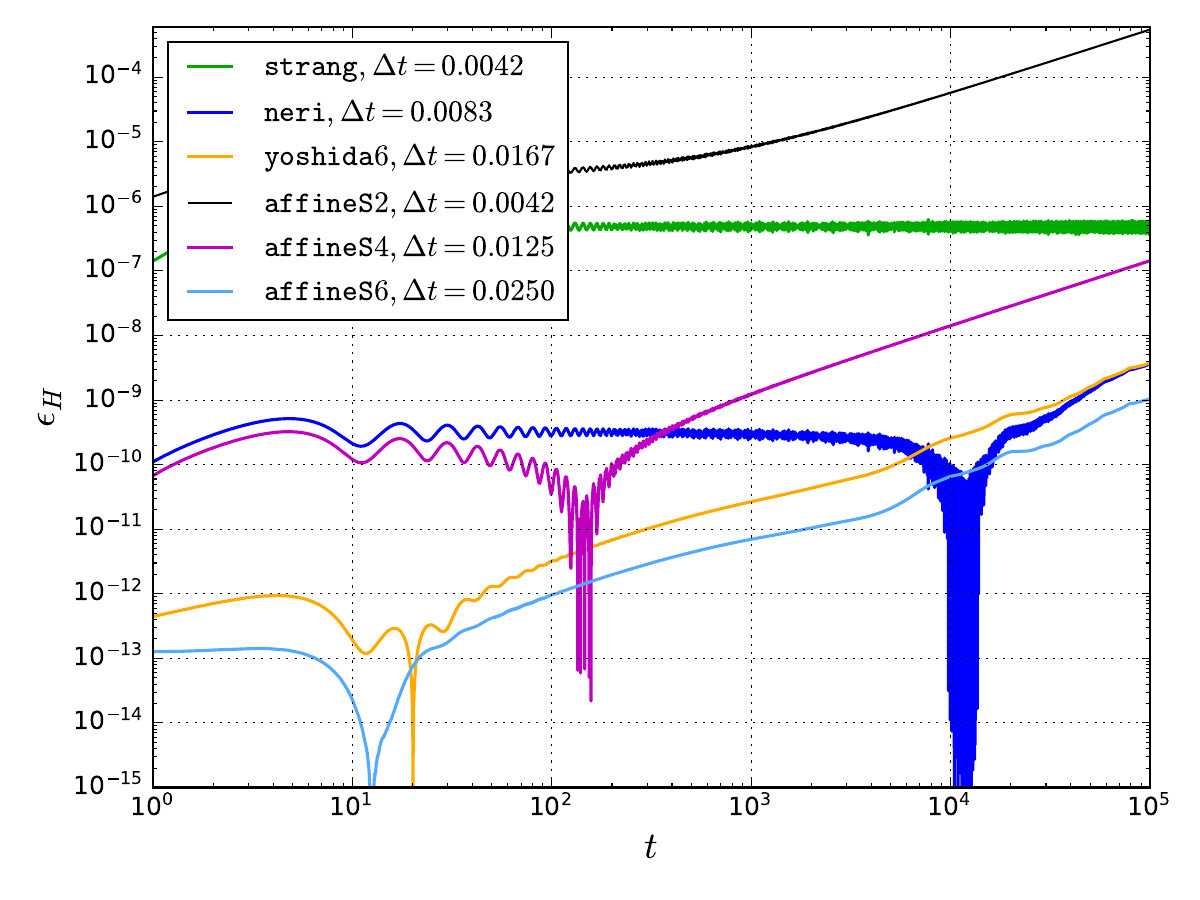}
    \caption{(Left) Relative error $\epsilon_M$ for the mass  as a function of time $t$, for different composition and symmetric affine splittings. (Right) Relative error $\epsilon_H$ for the Hamiltonian. The step size $\Delta t$ indicated for each method is selected in order to give the same computational cost $\mathcal{C}\approx4.8\cdot 10^{7}$, and $t \in [0,10^5]$. The parameters for the reference solution are $\eta=1,c=0.001,x_0=\phi_0=0$.} 
    \label{fig:NLSE3_fourier_soliton_err_H_M_long_times}
\end{figure}

To investigate whether the favorable preservation characteristics of affine methods are solely a result of approximating the exact solution with high accuracy, we present the long-time behavior  (up to $t=10^5$) of relative errors for both mass (left panel) and Hamiltonian (right panel) in Figure \ref{fig:NLSE3_fourier_soliton_err_H_M_long_times}. In this experiment, the interval is $I=[-50\pi, 50\pi]$ and the soliton speed is $c=0.001$, ensuring that the numerical support of the solution remains well within the interval boundaries for all $t$. The choice of the step size  $\Delta t$  for each scheme aims to maintain a consistent computational cost $\mathcal{C}\approx 4.8 \cdot 10^{7}$.
In the left panel, we observe that the mass relative error  $\epsilon_M$ linearly grows with time for all methods due to the accumulation of round-off errors. Composition schemes exhibit similar performance with very low initial errors, consistent with the fact that each stage of composition is nearly mass-preserving (subject to round-off and pseudo-spectral discretization errors). They clearly outperform affine schemes of order 2 and 4. However,  the sixth-order affine scheme yields the lowest mass error in the most efficient manner even though it does not consist of mass-preserving propagators. On the other hand, the right panel shows the long-time behavior of the Hamiltonian relative error $\epsilon_H$.  It is theoretically expected for composition methods that the error in the Hamiltonian remains bounded and decreases with the order of the method  \cite{hairer2006geometric, blanes2017concise}.  However, the figure shows that only the Strang-Marchuk scheme adheres to this theoretical upper bound, and this assertion holds only within the considered time span. All the other schemes exhibit an initially bounded error, then a transition zone and finally a linear growth, which can be attributed to the accumulation of round-off errors. Also in this case the sixth-order affine method emerges as the most efficient and accurate.

In summary, while composition schemes theoretically preserve the mass and the Hamiltonian for the NLSE3, in practical (numerical) terms their exactness is limited by machine precision, so in this regard we can consider the sixth-order affine method as highly competitive and almost-preserving within the standard double-precision limits of typical computers. To leverage the theoretical properties of composition schemes would require to resort to computationally-demanding multi-precision arithmetic calculations without direct hardware support (which, of course, could be justified in certain applications, e.g. high-precision long-time simulations in celestial mechanics).  

\subsubsection{Fractional cubic nonlinear Schrödinger equation (fNLSE3)}
The NLSE3 has been generalized recently by defining the fractional derivative via the Fourier multiplier
\begin{equation}
    \partial_x^{\alpha/2} u = \mathcal{F}^{-1}\big\{ (\mathrm{i} |k|)^{\alpha/2} \mathcal{F} u\big\},
\end{equation}
and introducing the fractional Laplace operator \cite{weitzner2003fractional, guo2012standingfNLS, klein2014fractional, duo2016spectralfNLS, mao2017fractional}
\begin{equation}
    (-\partial_{x}^2)^{\alpha / 2} u(x) := \mathcal{F}^{-1}\left(|k|^{\alpha}\mathcal{F}u\right)=\frac{1}{\sqrt{2\pi}}\int_\mathbb{R} |k|^{\alpha}\hat{u}(k)\mathrm{e}^{\mathrm{i}kx}\mathrm{d}k,
    \label{eq:frac_laplacian}
\end{equation}
where $1 < \alpha \leq 2$ is the {\em Lévy index} (the usual Laplacian corresponds to the case $\alpha=2$), giving the fractional cubic NLSE (fNLSE3)
\begin{equation}
    \mathrm{i} \partial_t u = \frac{1}{2}(-\partial_{x}^2)^{\alpha/2} u \pm |u|^{2} u.
    \label{eq:fNLSE3}
\end{equation}
In the focusing case (minus sign in the nonlinear term), the fNLSE \eqref{eq:fNLSE3} admits standing wave solutions of the form $u(x,t) = \psi(x) \mathrm{e}^{\mathrm{i} \omega t}$ with $\omega \in \mathbb{R}$, where $\psi(x)$ solves the nonlinear elliptic equation \cite{guo2012standingfNLS, frank2013fractional}
\begin{equation}
\frac{1}{2} \left(-\partial_x^2\right)^{\alpha/2} \psi + \omega \psi -|\psi|^{2} \psi  =  0.
\label{eq:ground_state_fNLS}
\end{equation}
For $\alpha<2$, no explicit solutions of \eqref{eq:ground_state_fNLS} are known, so they must be found numerically.
Like the standard NLSE3, the fractional version conserves the mass and the Hamiltonian, which in this case is given by
\begin{equation}
    H(u) = \frac{1}{2}  \int_\mathbb{R} \left( |\partial_x^{\alpha/2} u|^2  \pm |u|^{4} \right) \mathrm{d}x.
        \label{eq:fNLSE3_hamiltonian}
\end{equation}
\begin{figure}[t]
    \centering
    \includegraphics[width=0.49\textwidth]{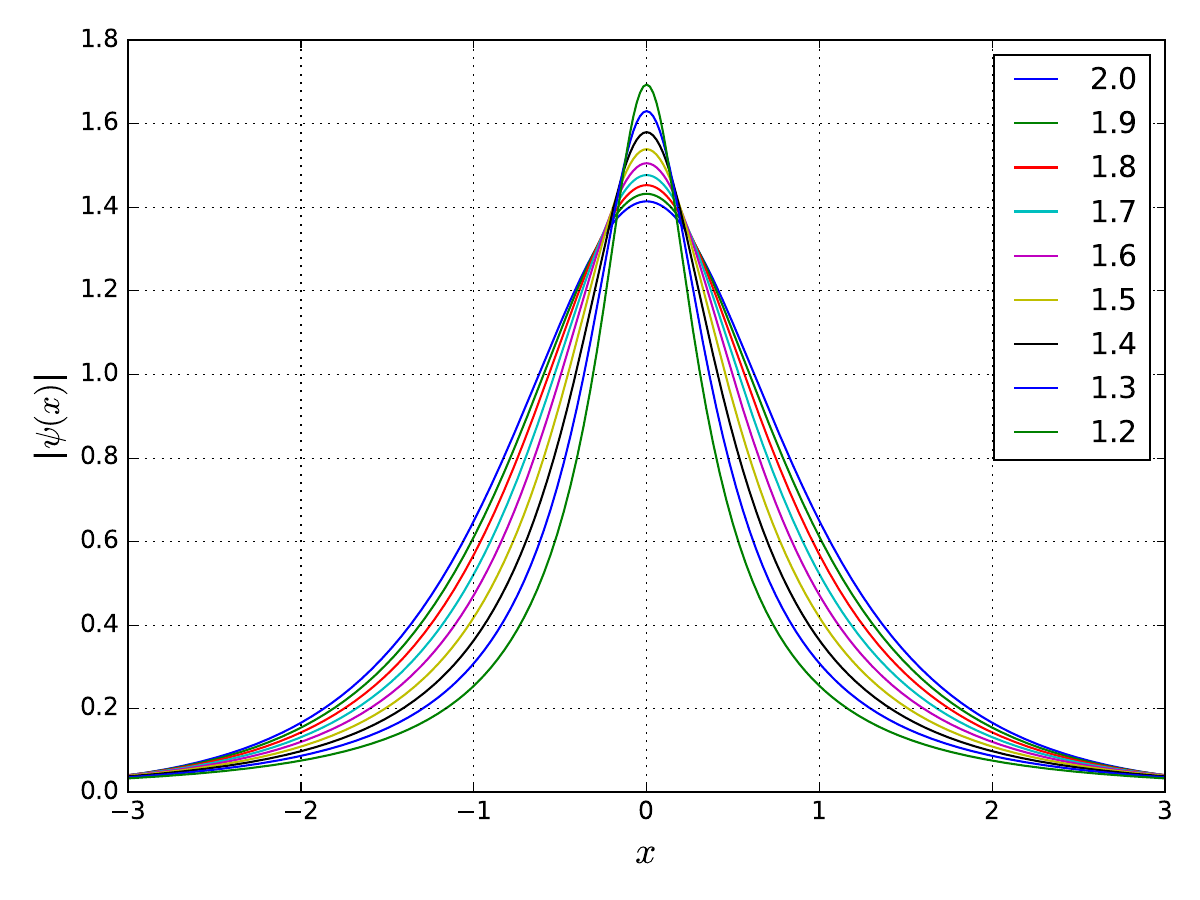}
    \includegraphics[width=0.49\textwidth]{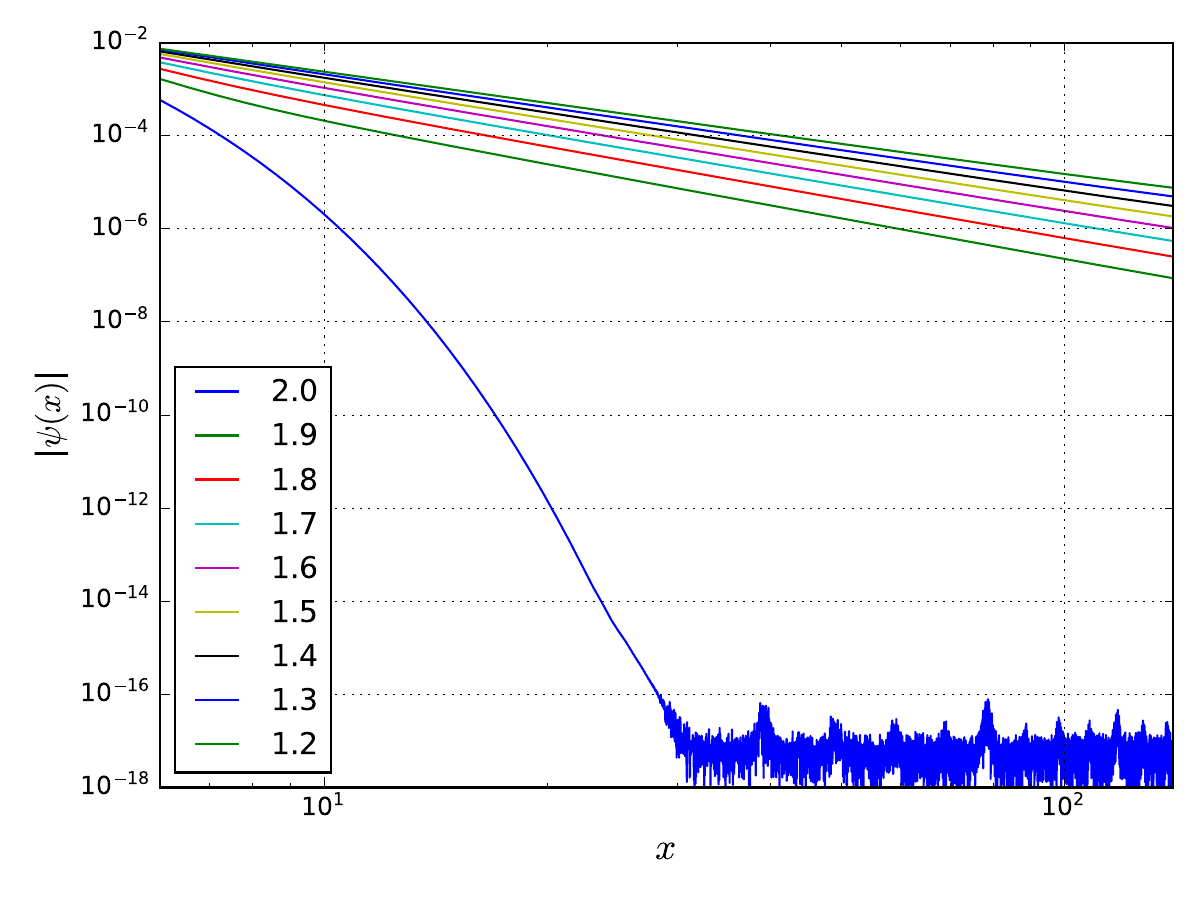}\\
    \includegraphics[width=0.49\textwidth]{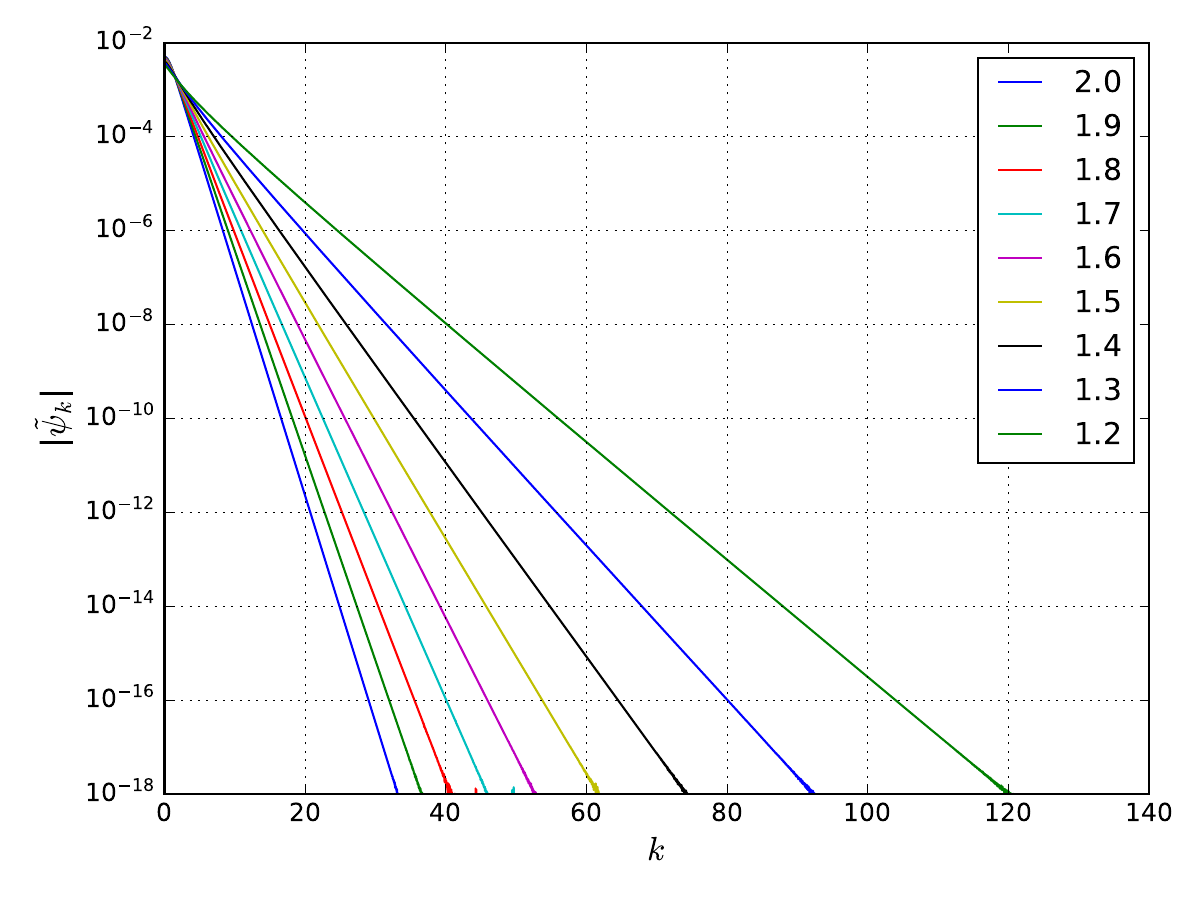}
    \caption{(Left) Absolute value of the numerically calculated standing wave solution of the fNLSE3 for various $\alpha$ and $\omega=1$. (Right) Asymptotic spatial decay of the absolute value $|\psi(x)|$ in log-log scale. When $\alpha=2$ (standard NLSE3) the decay is exponential, while for $\alpha<2$ the decay is algebraic, and is faster for higher values of $\alpha$. (Lower) Modulus of Fourier coefficients of the fNLSE3 standing wave solutions in semi-log scale, which show exponential decay.}
    \label{fig:fNLSE3_ground_states}
\end{figure}

The fractional {\em linear} Schrödinger equation has been investigated by Laskin as a model for fractional quantum mechanics \cite{laskin2000fractional}. Later, Longhi \cite{longhi2015fractional} presented a proposal for its optical implementation, which was experimentally realized by \cite{liu2023experimental}. The fractional nonlinear Schrödinger equation has been studied theoretically in \cite{weitzner2003fractional, guo2012standingfNLS, frank2013fractional} and numerically in \cite{klein2014fractional, duo2016spectralfNLS, mao2017fractional}. 

Although it shares the hamiltonian character of the standard NLSE, no exact analytical solutions are known for the fNLSE. This fact forces us to resort to numerically computed reference solutions to test our methods. To this end, we construct standing wave solutions whose amplitude is obtained by solving the nonlinear equation \eqref{eq:ground_state_fNLS} with high accuracy, using the Newton-Krylov algorithm and the LGMRES solver implemented in the \texttt{scipy} package of the numerical Python ecosystem \cite{virtanen2020scipy}. The states are discretized using a Fourier pseudo-spectral method with $N=2^{15}$ modes on the interval $I = [-300, 300]$. The amplitude of the standing wave is guaranteed to satisfy equation \eqref{eq:ground_state_fNLS} with an absolute tolerance of $10^{-12}$ in the max-norm. In Figure \ref{fig:fNLSE3_ground_states}, we display the amplitude of standing waves of the fNLSE3 (left panel) and their spatial decay (right panel) for $\omega=1$ and values of the Lévy index $\alpha$ ranging from 1.2 to 2.0. The spatial decay clearly follows an algebraic law for $\alpha<2$ and the decay rate slows down as the Lévy index $\alpha$ decreases. The case $\alpha=2$ corresponding to the stationary state of the standard NLSE3 is included for comparison, with exponential spatial decay.
The lower panel of Figure \ref{fig:fNLSE3_ground_states}, on the other hand, exhibits the modulus of the pseudo-spectral Fourier coefficients for the same stationary states in semi-log scale, which in all cases show an exponential decay. The decay rate is lower for higher values of $\alpha$.  These figures highlight that as $\alpha \to 1$, more Fourier modes and a wider interval are required to obtain an accurate pseudo-spectral approximation of the standing wave envelope in the fractional case.

In order to test the splitting methods we take as reference solution the standing wave  with $\omega=1$, i.e. $u_\mathrm{ref}(x,t)=\psi(x)\mathrm{e}^{\mathrm{i}t}$, where $\psi(x)$ is obtained numerically as previously described. To avoid redundancy, we show the results and make some general comments, since the procedure is similar to that employed for the standard Schrödinger equation in the previous subsection.
\begin{figure}[t]
    \centering
    
    \includegraphics[width=0.49\textwidth]{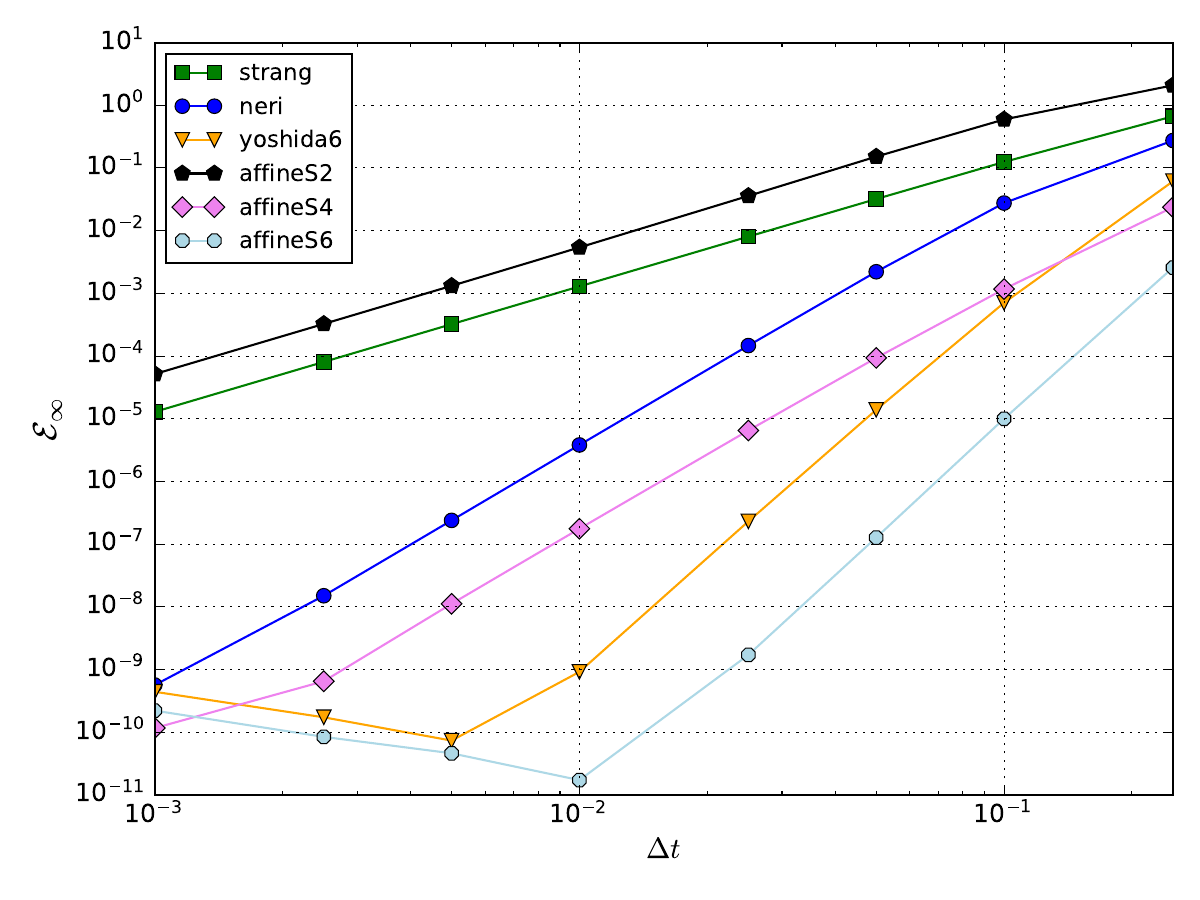}
    \includegraphics[width=0.49\textwidth]{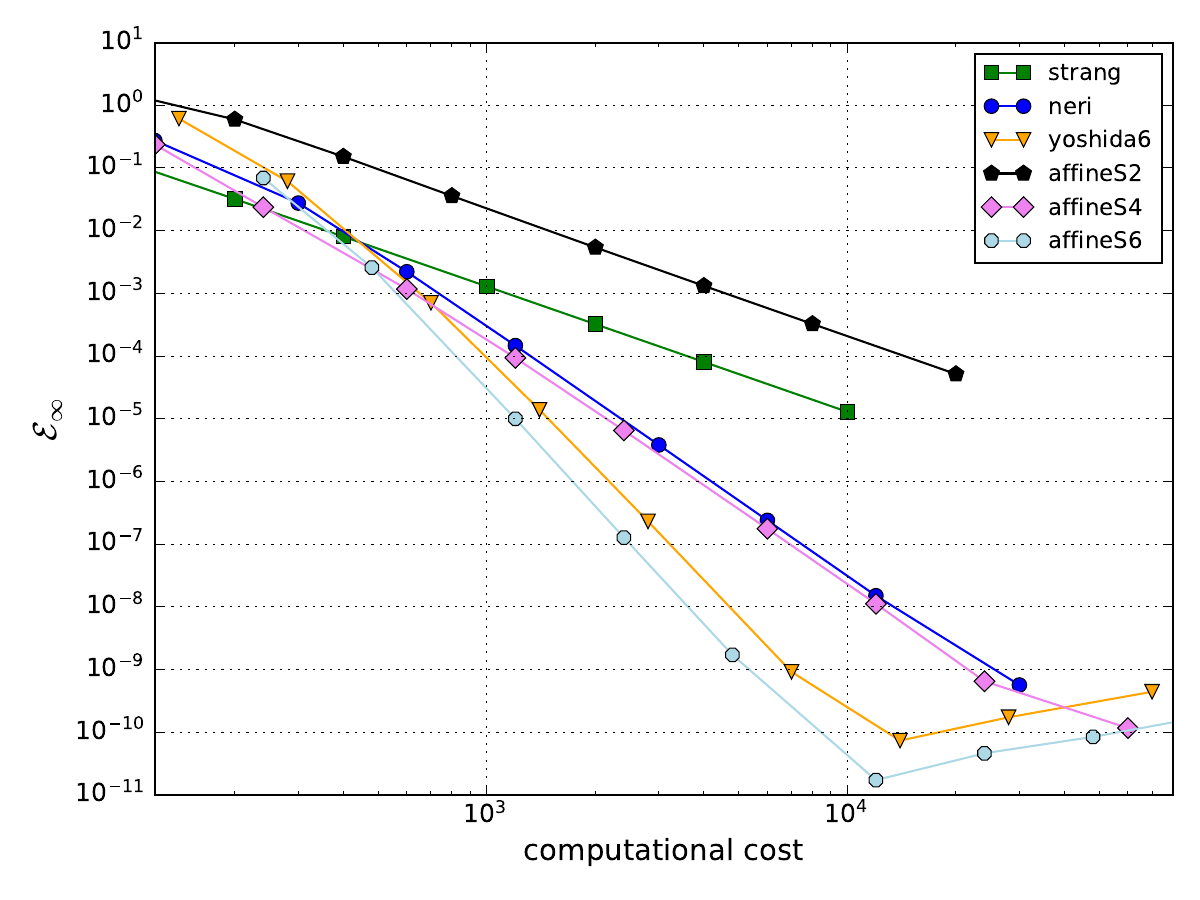} \\
    
    \includegraphics[width=0.49\textwidth]{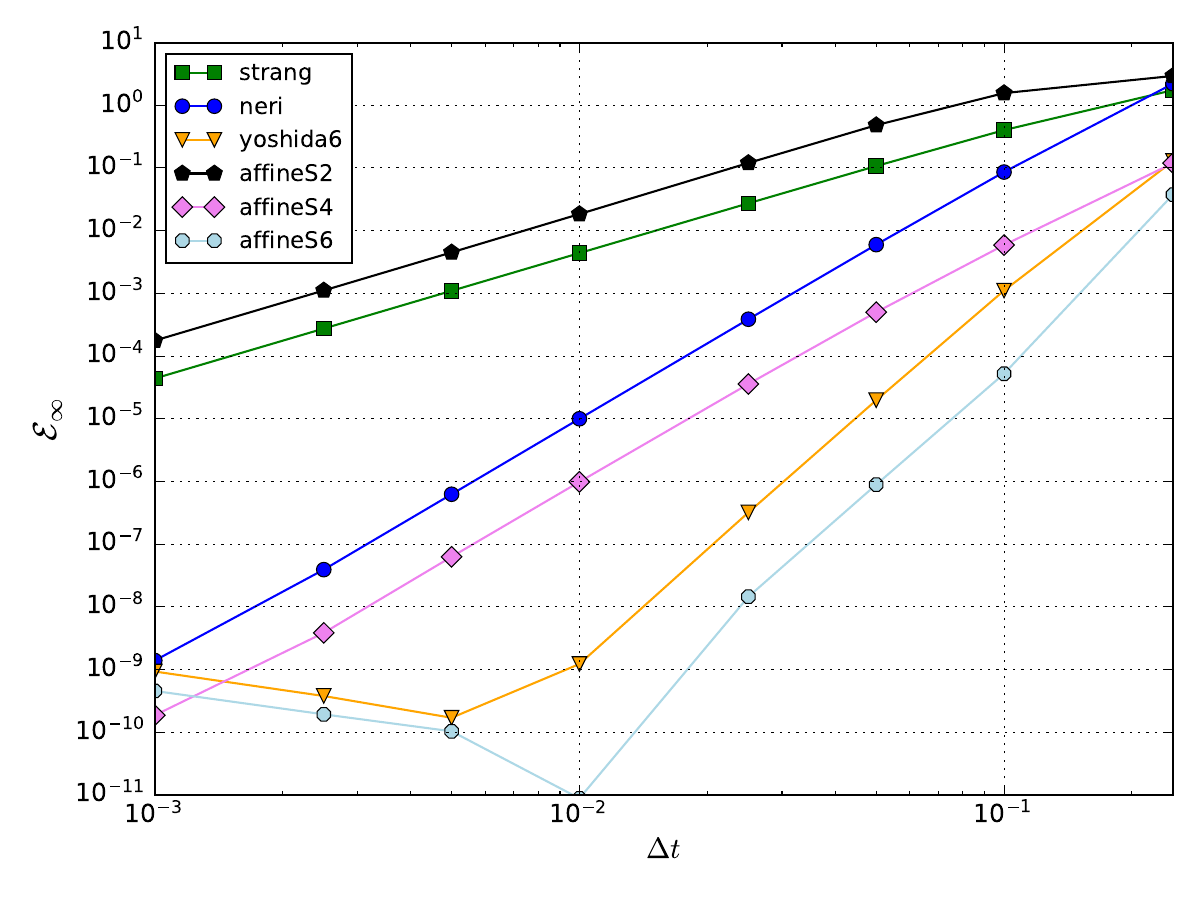}
     \includegraphics[width=0.49\textwidth]{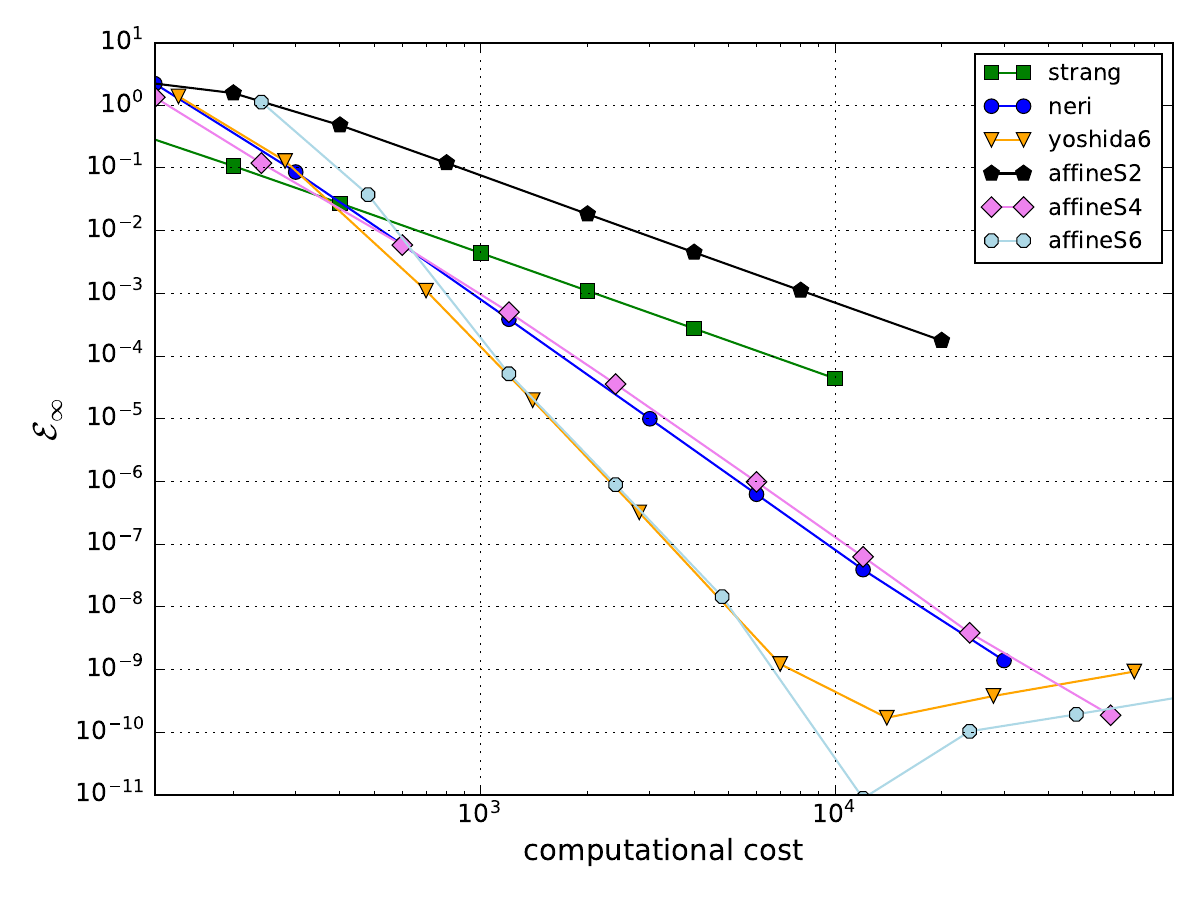}
    \caption{(Left) Absolute error $\mathcal{E}_\infty (t=10)$ for the numerical solution of the focusing fNLSE3, as a function of the time step $\Delta t$, for different composition and symmetric affine splittings. (Right) Efficiency plot. (Upper) $\alpha=1.8$. (Lower) $\alpha=1.3$. Reference solution is $u_\mathrm{ref}(x,t)=\psi(x)\mathrm{e}^{\mathrm{i}t}$, where $\psi(x)$ is the numerical solution of \eqref{eq:ground_state_fNLS}. States and operators are approximated using a Fourier pseudo-spectral method with $N=2^{15}$ modes on the interval $I=[-300,300]$.}
    \label{fig:fNLSE3_fourier_stationary_err_inf}
\end{figure}
In the left panels of Figure \ref{fig:fNLSE3_fourier_stationary_err_inf} we plot the error of the calculated solution at $t=10$ as a function of the time step $\Delta t$ for the investigated numerical schemes, while the right panels show the corresponding efficiency plot. The upper panels correspond to $\alpha=1.8$ and the lower ones to $\alpha=1.3$. We observe that the fourth and sixth-order affine schemes outperform composition schemes of the same order in both metrics for $\alpha=1.8$. When $\alpha=1.3$, fourth and sixth-order affine schemes exhibit similar efficiency when compared with the composition ones, but with lower absolute errors for identical step sizes. Overall, our numerical investigation reveals that affine schemes of high order perform increasingly better than composition schemes as the Lévy index grows from $\alpha=1$ to $\alpha=2$.

\begin{figure}[t]
    \centering
    \includegraphics[width=0.49\textwidth]{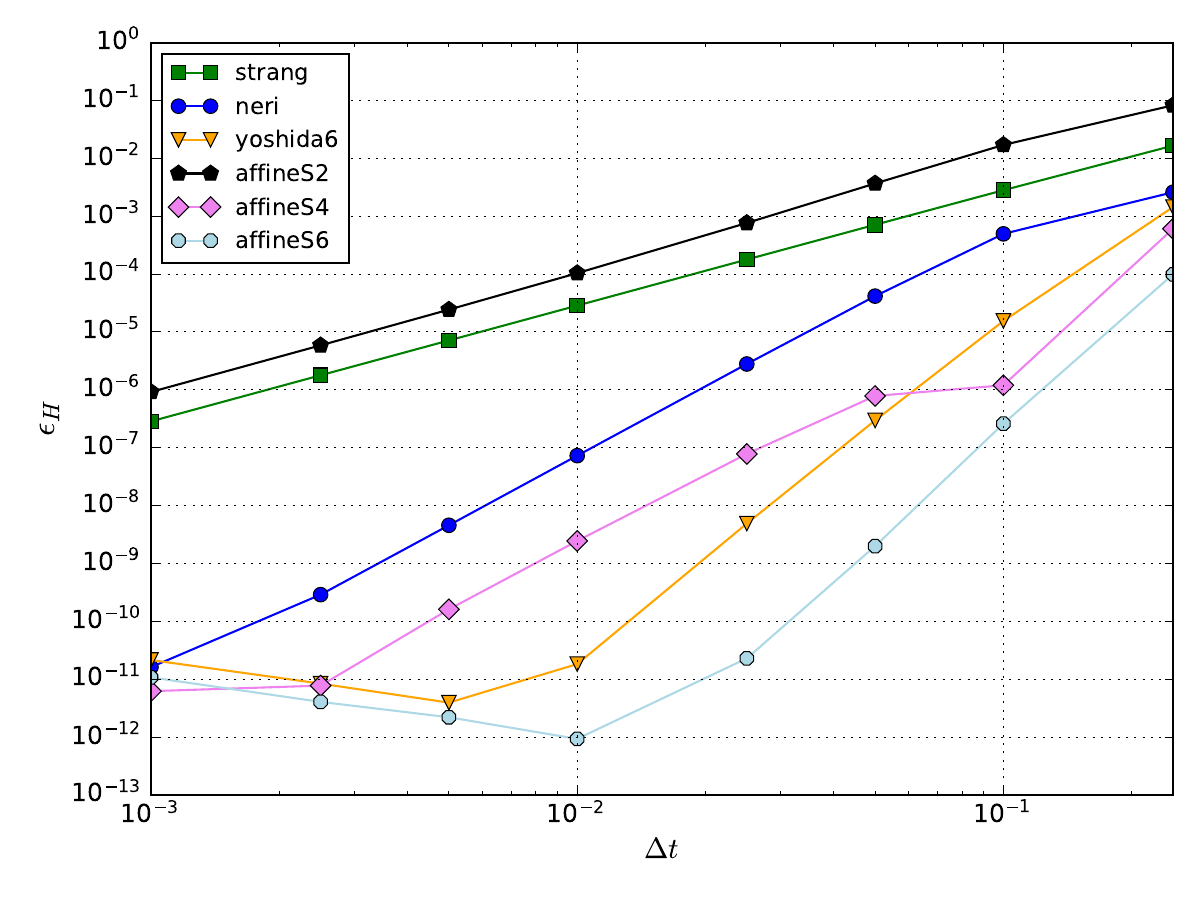}
    \includegraphics[width=0.49\textwidth]{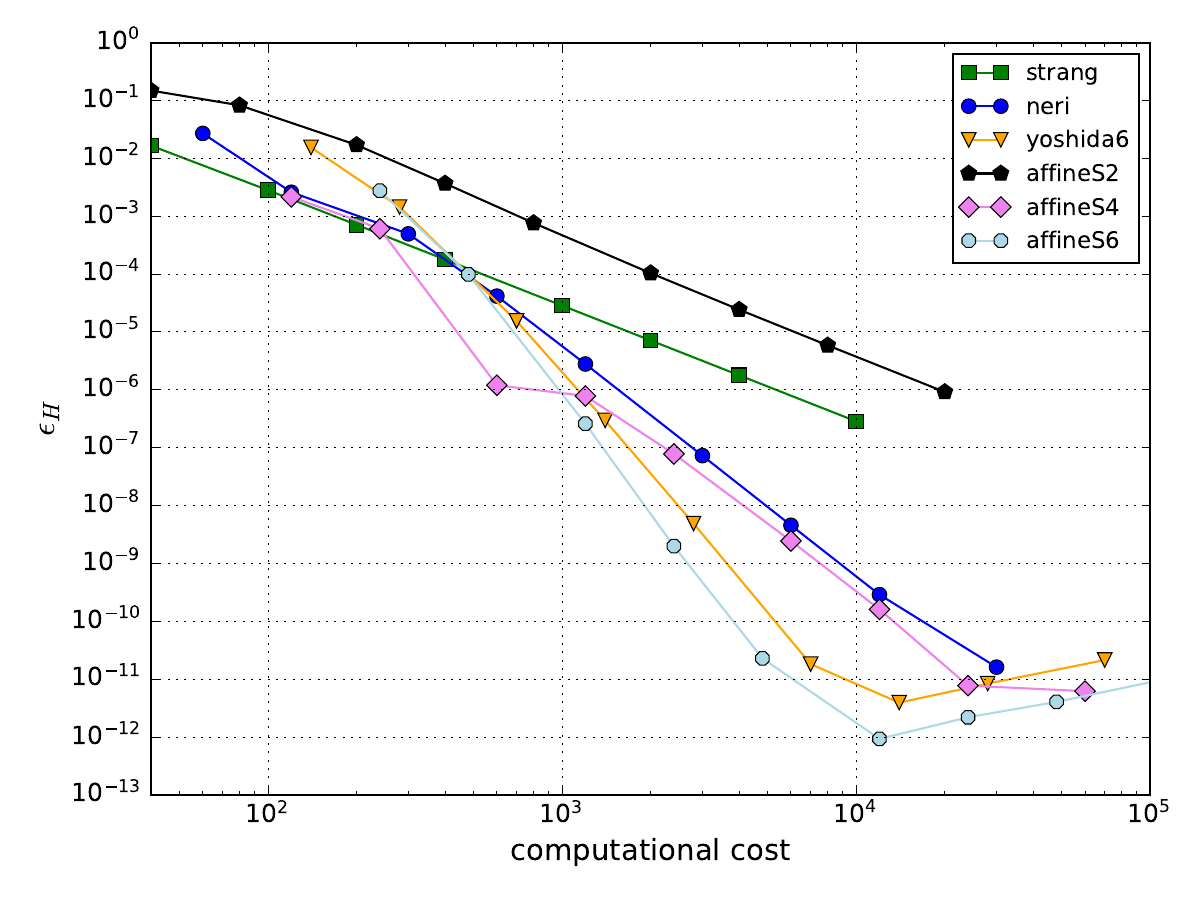} \\
    \includegraphics[width=0.49\textwidth]{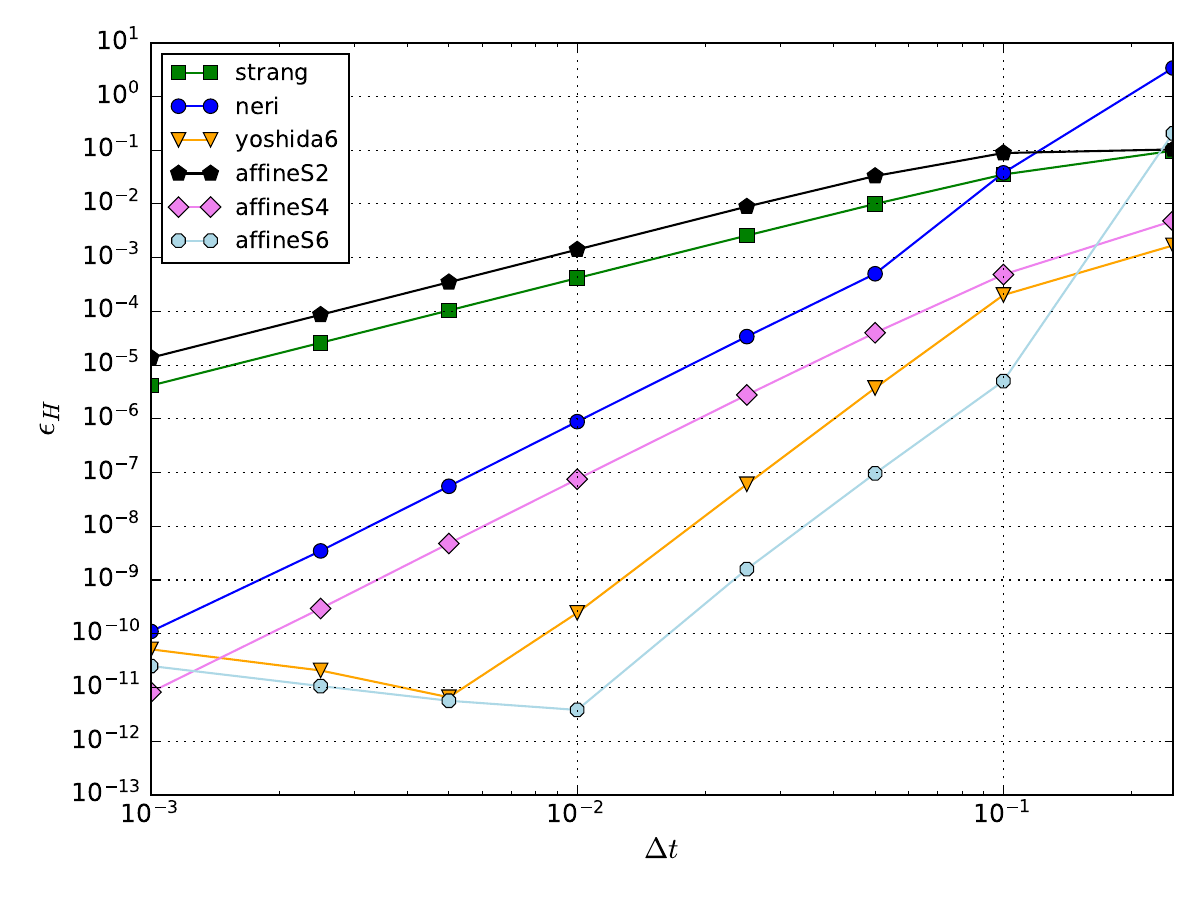}
    \includegraphics[width=0.49\textwidth]{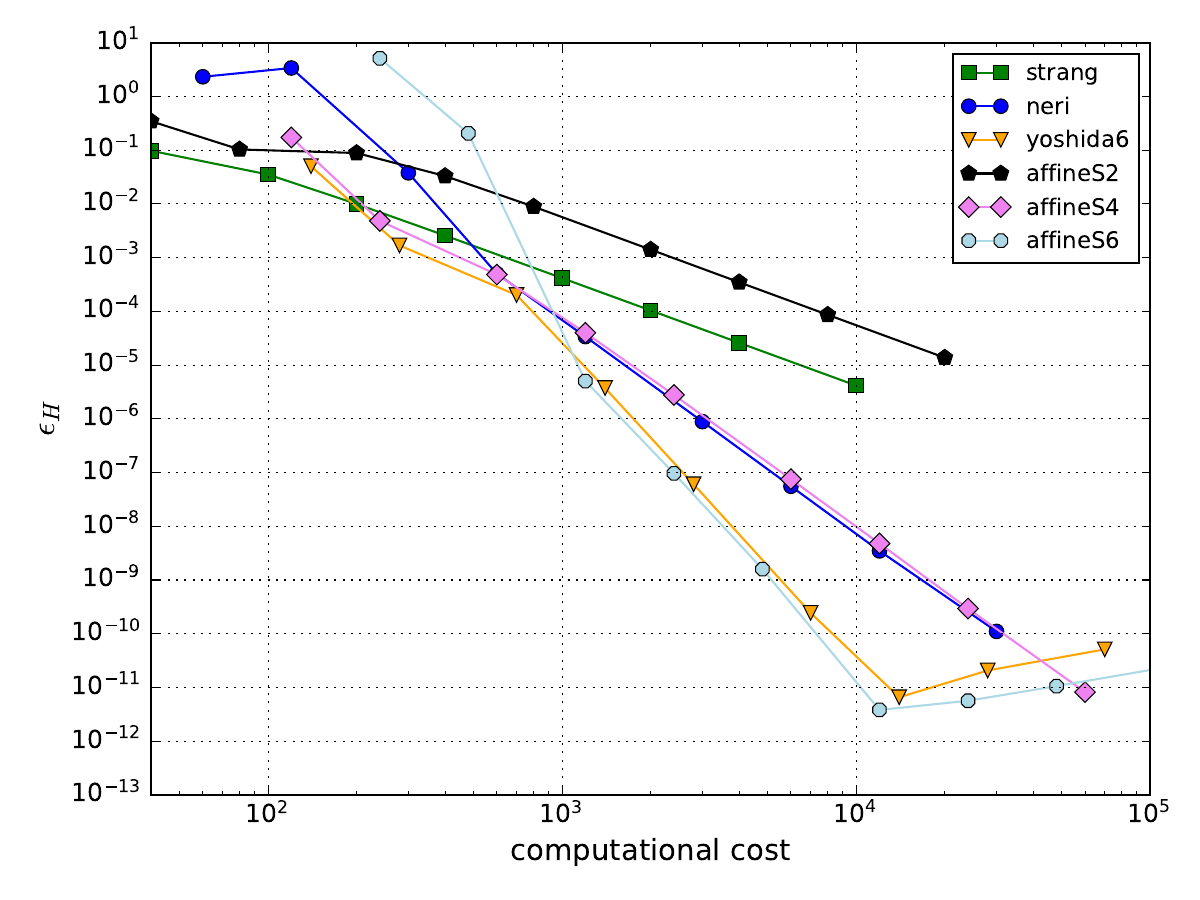}
    \caption{(Left) Hamiltonian relative error $\epsilon_H(t=10)$ for the solution of the focusing fNLSE3, as a function of the time step $\Delta t$, for different composition and symmetric affine splittings. (Right) Efficiency plot. The reference solution is the same of Figure \ref{fig:fNLSE3_fourier_stationary_err_inf}  for $\alpha=1.8$ (upper panels) and $\alpha=1.3$ (lower panels).}
    \label{fig:fNLSE3_fourier_stationary_err_H}
\end{figure}

Moreover, affine schemes preserve the Hamiltonian better than composition ones for all values of $\alpha$ and with lower computational cost, as exemplified in the right panels of Figure \ref{fig:fNLSE3_fourier_stationary_err_H} for $\alpha=1.8$ (upper panel) and $\alpha=1.3$ (lower panel).
With respect to mass conservation, the affine schemes perform equivalently to composition ones only for step sizes lower than a specific threshold, as was the case with the standard NLSE3. This threshold is approximately given by the step size for which the affine scheme minimizes the absolute error $\mathcal{E}_\infty$ of the numerical solution.

\begin{figure}[t]
    \centering
    \includegraphics[width=0.49\textwidth]{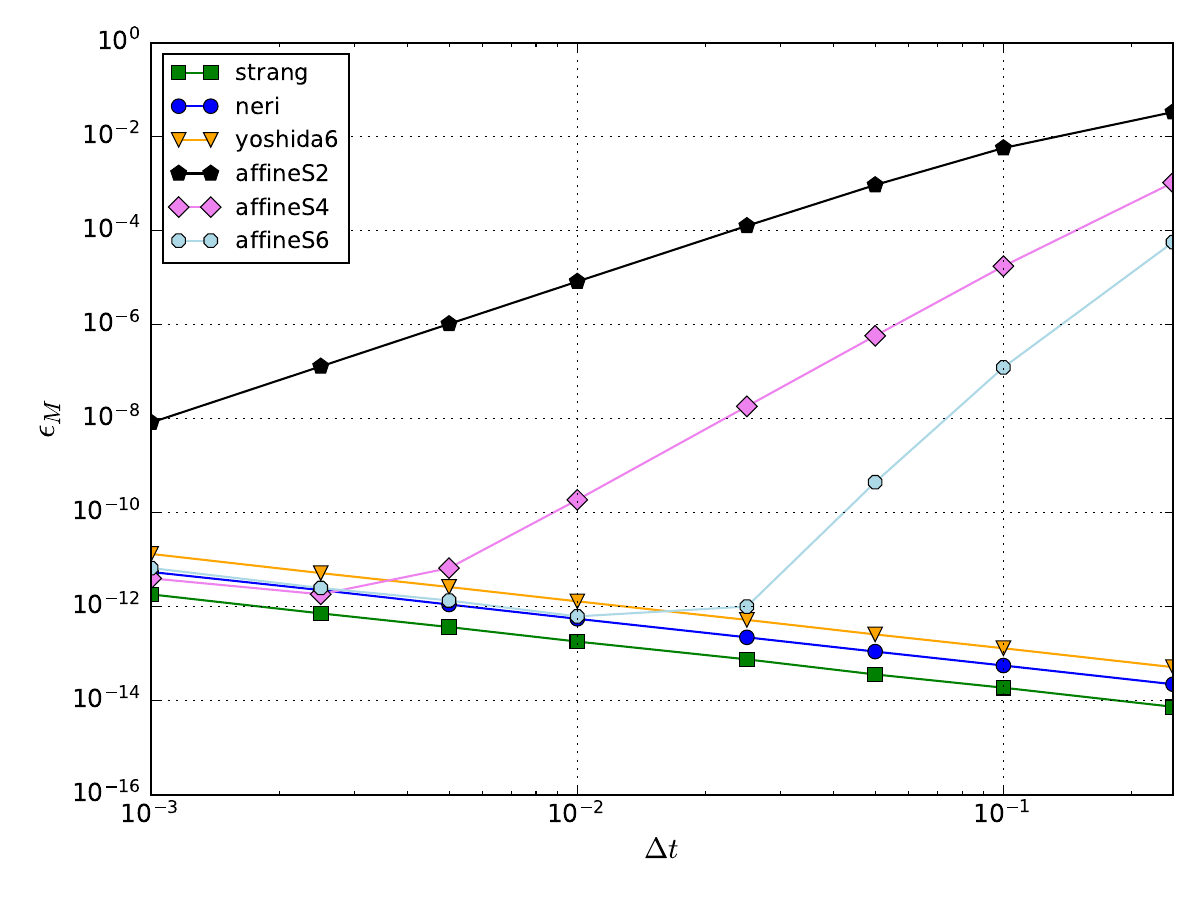}
    \includegraphics[width=0.49\textwidth]{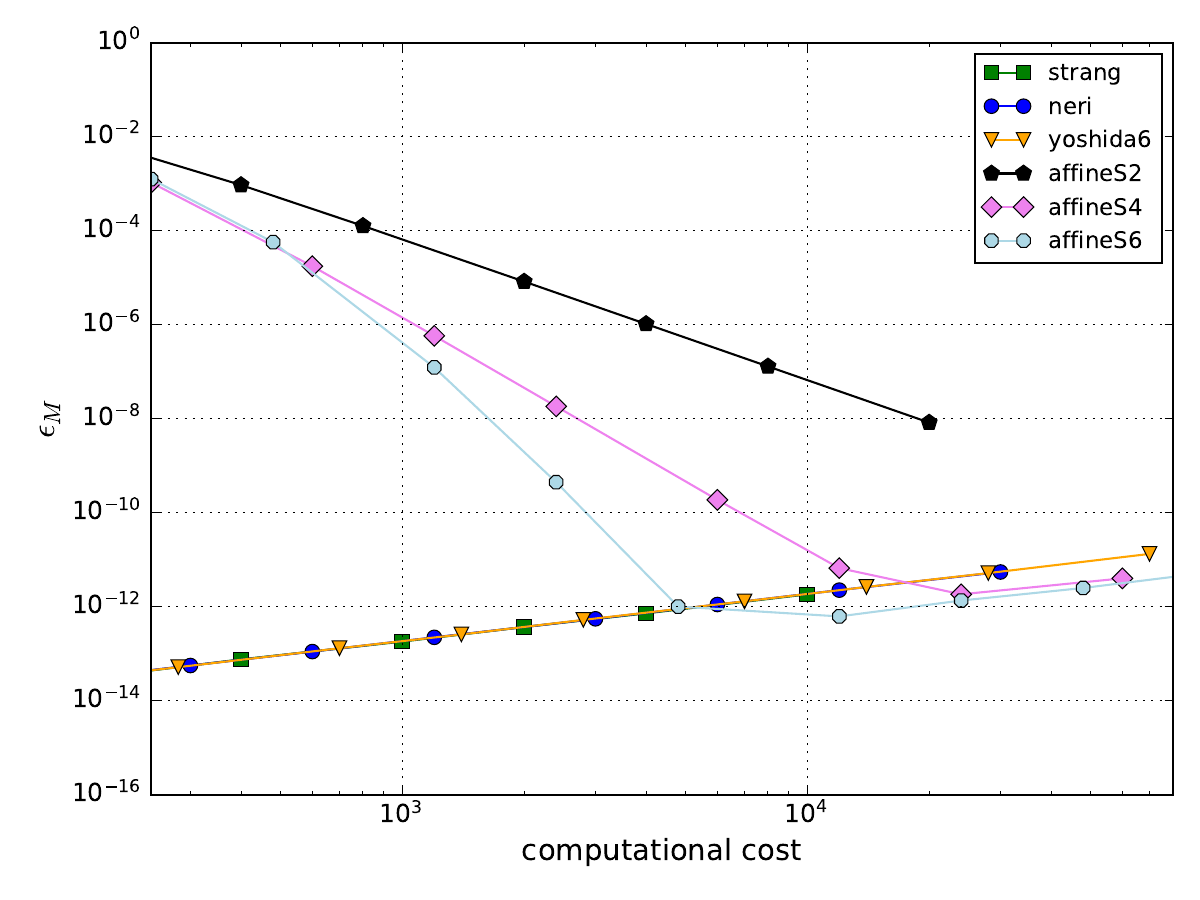} \\
    \includegraphics[width=0.49\textwidth]{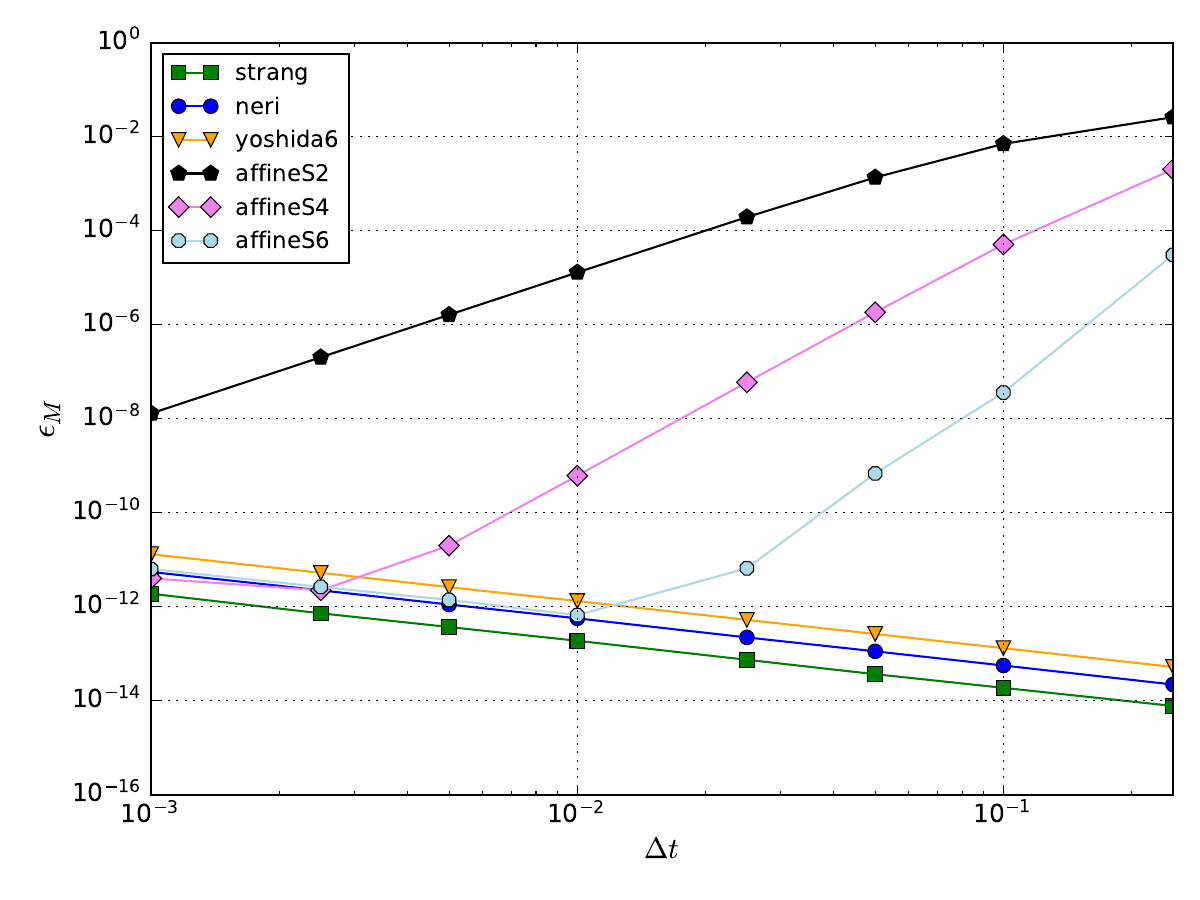}
    \includegraphics[width=0.49\textwidth]{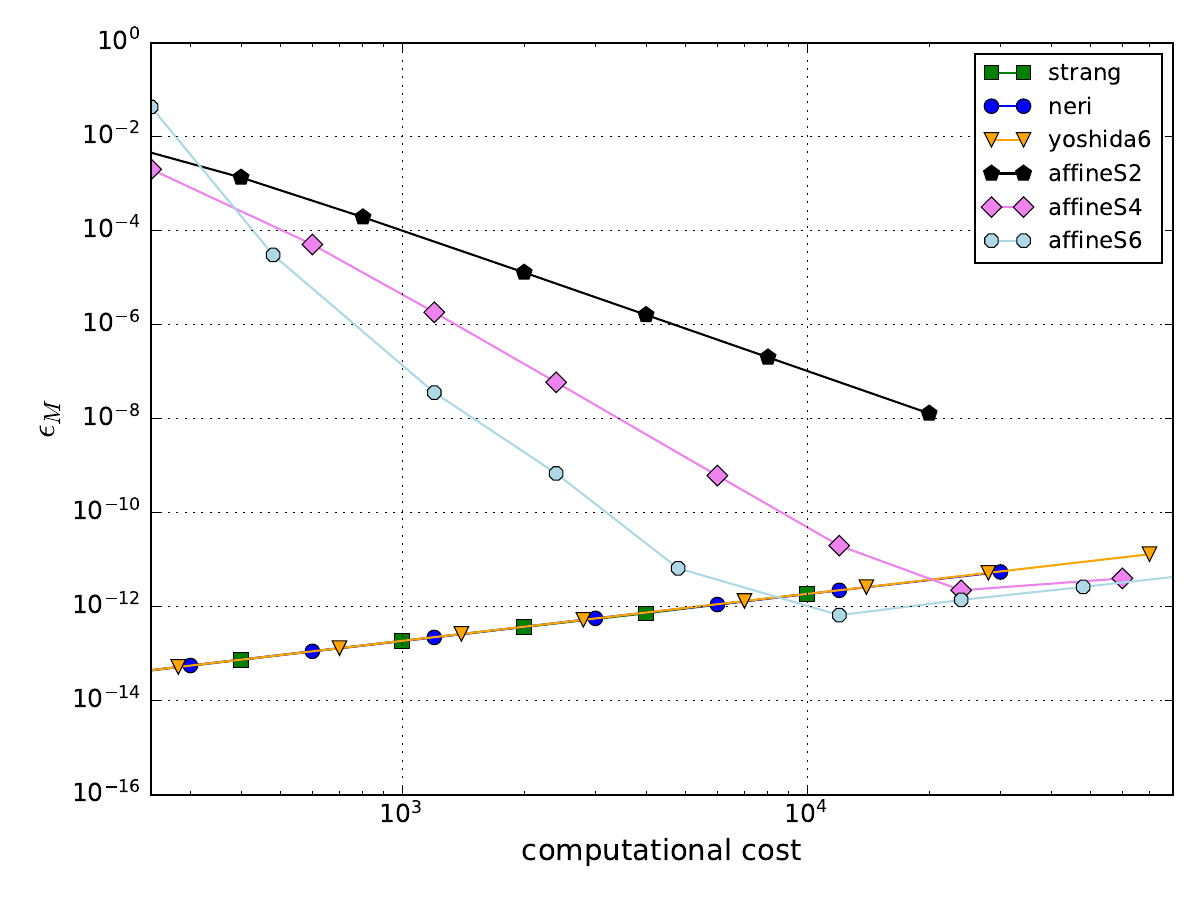}
    \caption{(Left) Mass relative error $\epsilon_M(t=10)$ for the solution of the focusing fNLSE3, as a function of the time step $\Delta t$, for different composition and symmetric affine splittings. (Right) Efficiency plot. The reference solution is the same of Figure \ref{fig:fNLSE3_fourier_stationary_err_inf}  for $\alpha=1.8$ (upper panels) and $\alpha=1.3$ (lower panels).}
    \label{fig:fNLSE3_fourier_stationary_err_M}
\end{figure}

The extensive numerical evidence presented supports the conclusion that high-order affine schemes are efficient and accurate for solving the standard and fractional cubic nonlinear Schrödinger equations, outperforming composition schemes in most cases. In particular, the sixth-order affine scheme shows great promise as a general-purpose solver for achieving high accuracy in all the investigated metrics with modest step sizes and low computational cost.

\subsection{Fisher's reaction-diffusion equation}
In order to test the proposed methods in a non-conservative irreversible model, we resort to the well-known reaction-diffusion equation of Fisher, 
\begin{equation}
    \partial_t u(x,t) = \alpha \partial_x^2 u(x,t) + \beta u(x,t)(1-u(x,t)),\qquad x\in\mathbb{R}, t>0,
\end{equation}
where $\alpha$ is a non-negative diffusion coefficient and $\beta$ a real parameter. This equation was proposed originally as a model for gene propagation, having also applications in neutron diffusion in nuclear reactors, combustion problems and population dynamics \cite{Sari2015}. The diffusive character of the equation poses a challenge for methods with negative time steps.

For the numerical experiment we take the initial state $u_0(x)=\sech^2(10 x)$ discretized using the Fourier pseudo-spectral method with $N=2^{10}$ modes on the interval $I=[-80, 80]$. The parameters for the model are $\alpha=2.0$ and $\beta=0.5$. We solve the resulting ODE by means of the adaptive Dormand--Prince 8(5,3) method with tolerances near machine precision and take this numerical solution as reference for measuring the error of the splitting schemes. 
\begin{figure}[t]
    \centering
    \includegraphics[width=0.49\textwidth]{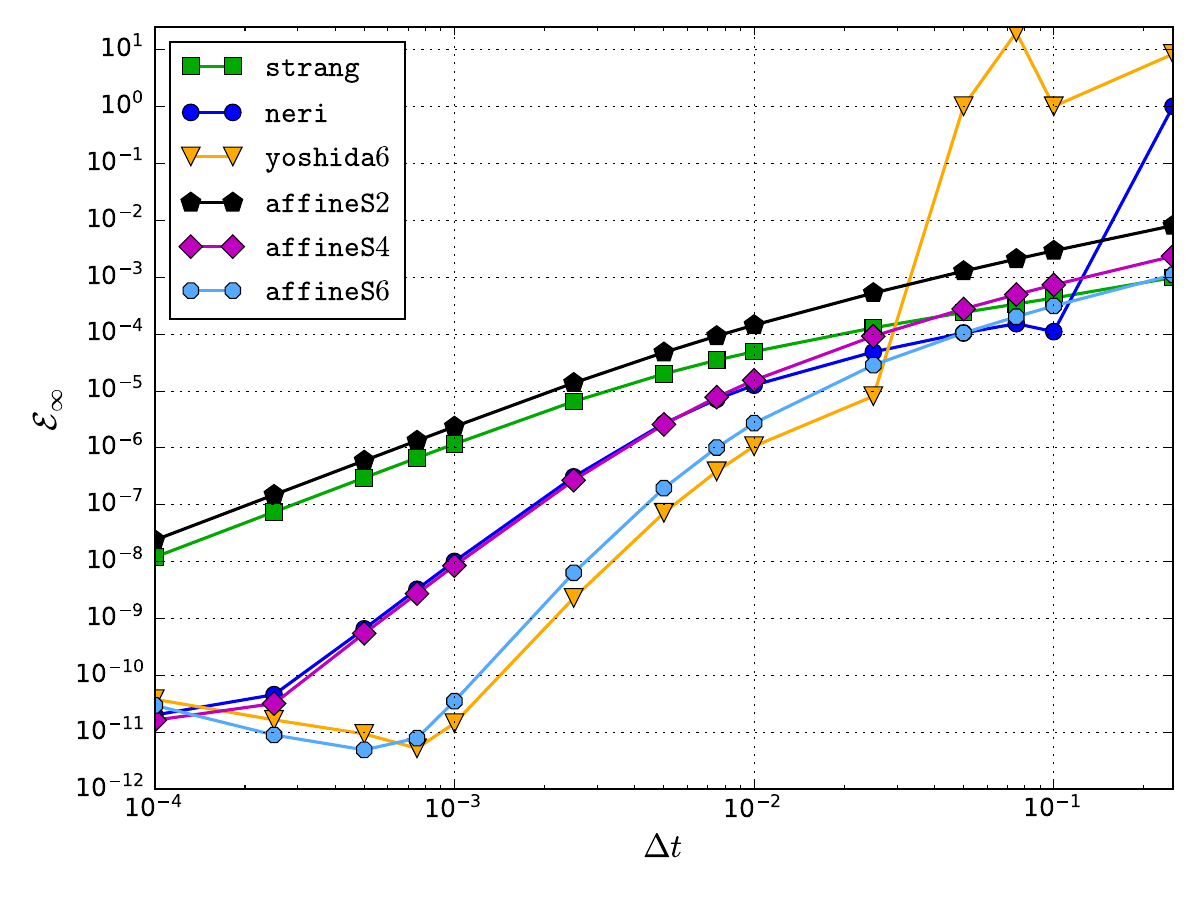}
    \includegraphics[width=0.49\textwidth]{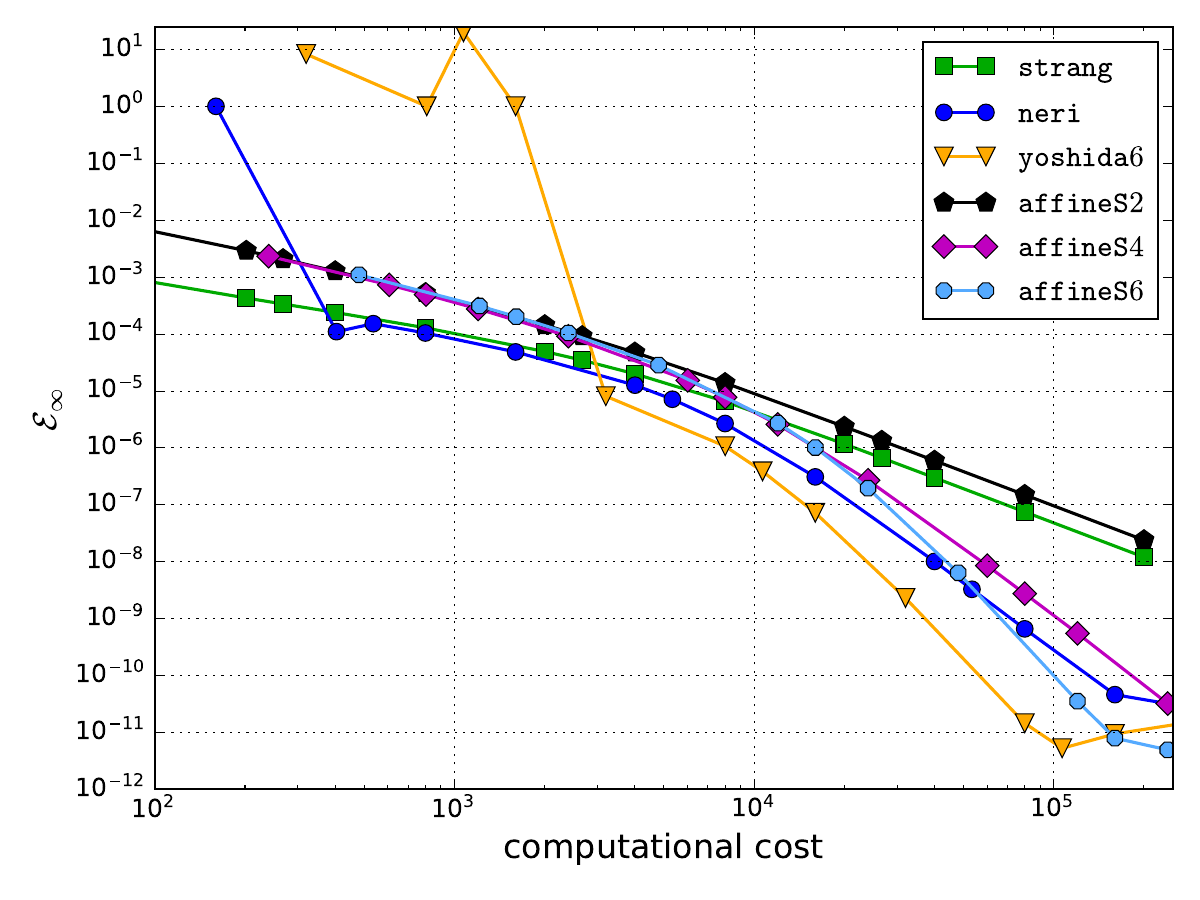}
    \caption{(Left) Absolute error $\mathcal{E}_\infty$  as a function of time $t$, for different composition and symmetric affine splittings. (Right) Computational cost. The parameters for the reference solution are $\alpha=2.0, \beta=0.5$.} 
    \label{fig:fisher_fourier_err}
\end{figure}
In Figure \ref{fig:fisher_fourier_err} we show the error $\mathcal{E}_\infty$ as function of the time step $\Delta t$ (left panel) and of computational cost (right panel). The error for methods of the same order is similar for sufficiently small time steps, but higher-order composition methods show instabilities for large time steps while affine methods exhibit a regular behavior. From the point of view of computational cost, composition methods (particularly the sixth-order Yoshida scheme) are slightly more efficient when they are stable. Additional numerical experiments show that the stability of high-order composition methods depends on the initial state and deteriorates as the diffusion coefficient increases. The sixth-order affine scheme guarantees both stability and efficiency when high accuracy is needed.

\subsection{The complex Ginzburg-Landau equation}
The complex Ginzburg-Landau equation (CGLE) is a canonical model for the description of light propagation in nonlinear dissipative media, phase transitions, chemical oscillations, and pattern formation \cite{aranson2002review, garcia-morales2012CGLE}. It can be seen as a generalized NLSE which includes linear and nonlinear gain/loss and diffusion. Models with cubic and quintic nonlinearities have been studied in the last decades, mainly in the context of nonlinear optics phenomena. Recently, these models have played a role in the theoretical and experimental investigation of \emph{dissipative solitons}, localized structures whose dynamics and interactions exhibit a variety of interesting behavior, including pulsating, exploding and bonding, as a consequence of the complex balance between linear, nonlinear and dissipative effects \cite{akhmediev1996singularities, akhmediev2001pulsating, akhmediev2001solitons, akhmediev2005dissipative, akhmediev2008dissipative, grelu2012dissipative, ferreira2022dissipative}. 

Stressing the analogy with \eqref{eq:NLSE3}, we write the cubic-quintic CGLE (CGLE5) in the form commonly studied by the nonlinear optics community
\begin{equation}
    \mathrm{i}\partial_t u = (\frac{1}{2} - \mathrm{i} \beta) (-\partial_{x}^2) u + \mathrm{i} \delta u + (\gamma + \mathrm{i}\varepsilon) |u|^{2} u +
    (-\nu + \mathrm{i} \mu) |u|^{4} u,
    \label{eq:CGLE5}
\end{equation}
where $\beta > 0$ is the diffusion coefficient that accounts for spectral filtering, $\delta$ is the linear gain/loss factor, and $\varepsilon$ is related to nonlinear gain processes. The quintic (last) term, responsible for high-order nonlinear effects, is essential for the existence of pulsating solutions \cite{akhmediev2001pulsating, akhmediev2001solitons, akhmediev2005dissipative, akhmediev2008dissipative}.

\subsubsection{Cubic complex Ginzburg-Landau equation (CGLE3)}

\begin{figure}[t]
    \includegraphics[width=0.49\textwidth]{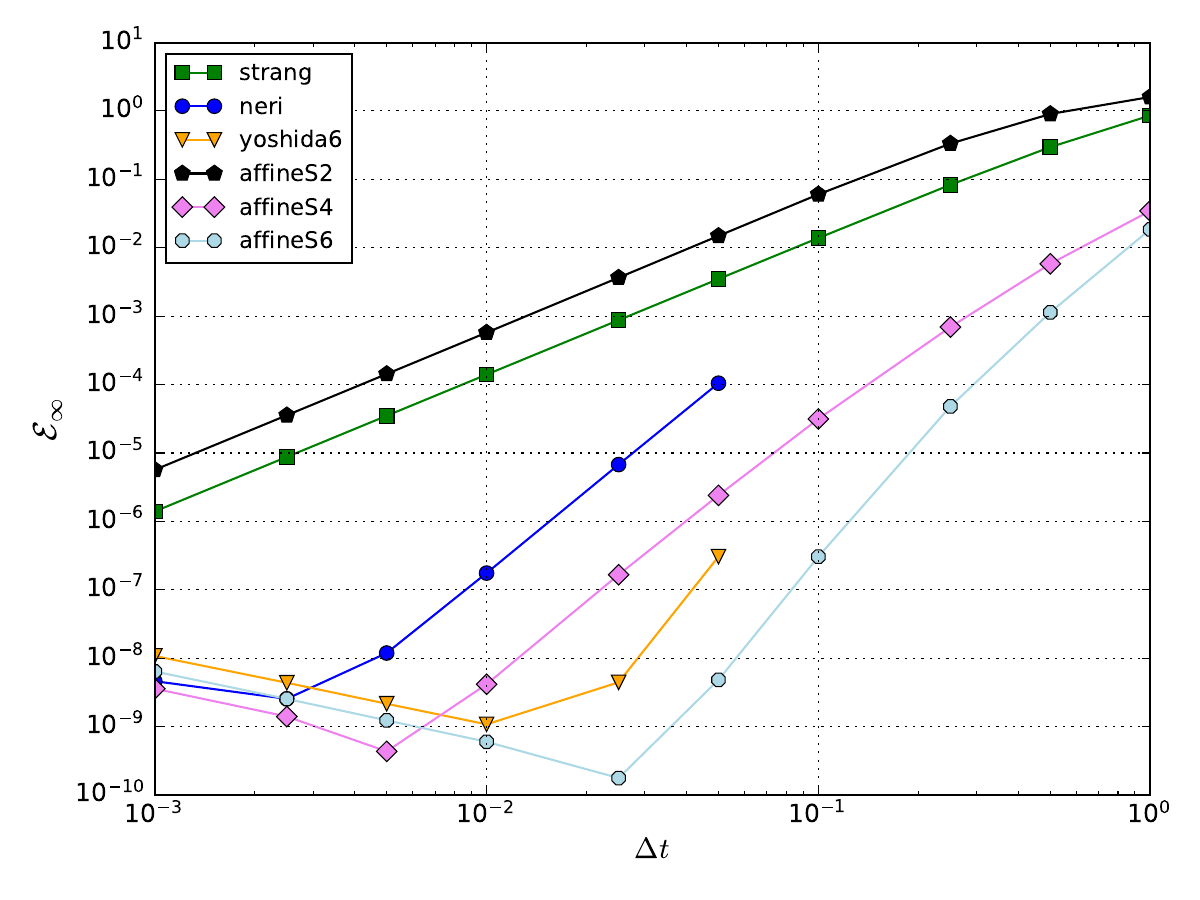}
    \includegraphics[width=0.49\textwidth]{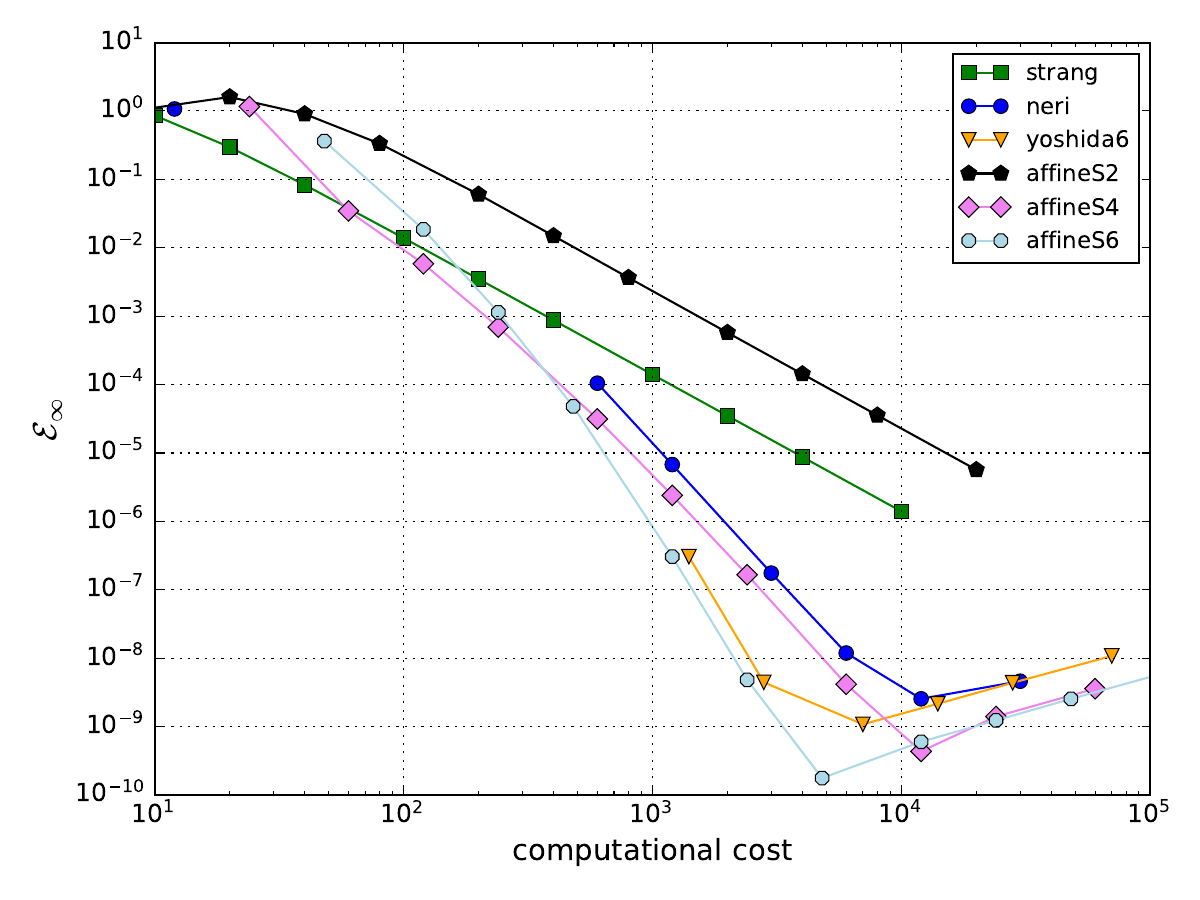}
    \caption{(Left) Absolute error $\mathcal{E}_\infty(t=10)$ for the numerical solution of the CGLE3, as a function of the time step $\Delta t$, for composition and symmetric splittings. (Right) Efficiency plot. The reference solution is the exact soliton \eqref{eq:CGLE3_soliton} of equation \eqref{eq:CGLE5} with $\beta=0.25$, $G=1$, and $\phi_0=0$, i.e. $u_\mathrm{ref}(x,t)=1.072\ \mathrm{sech}(x)\exp\left[\mathrm{i} \left(0.236\, \ln\left(1.072\, \mathrm{sech}(x)\right)-\omega t\right)\right]$. States and operators are discretized using a Hermite pseudo-spectral method with $N=300$ and scaling $s=1$.}
    \label{fig:CGLE3_fourier_soliton_err}
\end{figure}
Although the CGLE is nonintegrable and has no known conserved quantities, some exact soliton solutions have been found \cite{akhmediev1996singularities}. In particular, in the cubic case ($\mu=\nu=0$) if $\delta=0$ there exist stable arbitrary amplitude solitons of the form 
\begin{equation}
    u^\mathrm{CGL}_\mathrm{sol}(x,t) = \varphi(x) \exp(-\mathrm{i} \omega t),
    \label{eq:CGLE3_soliton}
\end{equation}  with
\begin{equation}
    \varphi(x) = \tilde{\varphi}(x) \exp\left[\mathrm{i}(\phi_0 + d \ln \tilde{\varphi}(x))\right], \quad
    \tilde{\varphi}(x) = G F \  \mathrm{sech}(G x),
\end{equation}
where $G>0$, and $d$, $\omega$ and $F$ are given by
\begin{equation}
   \lambda = \sqrt{1 + 4 \beta^2}, \quad \omega =  - \frac{d \lambda^2 G^2}{2\beta}, \quad 
   d = \frac{\lambda-1}{2\beta}, \quad  F = \frac{(2+9\beta^2) \lambda (\lambda - 1)}{2\beta^2(3\lambda - 1).}
\end{equation}
Taking this exact soliton with $\beta=0.25$ and $G=1$ as a reference solution, in Figure \ref{fig:CGLE3_fourier_soliton_err} we show as usual the dependency of the error on the time step (left panel) and on the computational cost (right panel), using in this case a Hermite pseudo-spectral discretization with $N=300$ and $s=1$. Again, the superiority of the higher-order affine splittings is evident. Besides, if the time step is larger than approximately $\Delta t=0.025$, composition methods of order 4 and 6 become unstable and are unable to compute a solution, as expected due to the presence of the diffusive term which renders the equation irreversible. Indeed, our numerical experiments show that the instability of high-order composition methods worsens as $\beta$ grows.

\subsubsection{Fractional cubic-quintic complex Ginzburg-Landau equation (fCGLE5)}

Recently, a fractional version of equation \eqref{eq:CGLE5} has been investigated  \cite{qiu2020solitonfCGLE}. In this variant, as with the fractional NLSE, the Laplace operator is substituted by its fractional generalization \eqref{eq:frac_laplacian}:
\begin{equation}
    \mathrm{i}\partial_t u = (\frac{1}{2} - \mathrm{i} \beta) (-\partial_{x}^2)^{\alpha/2} u + \mathrm{i} \delta u + (\gamma + \mathrm{i}\varepsilon) |u|^{2} u +
    (-\nu + \mathrm{i} \mu) |u|^{4} u.
    \label{eq:fCGLE5}
\end{equation}
Please note that our notation differs slightly from that used in \cite{qiu2020solitonfCGLE} and is more similar to the one introduced in \cite{akhmediev2001pulsating}. 
As far as we know, no exact solutions were found for this equation, so approximate and numerical methods are essential for the investigation of its properties.
\begin{figure}[t]
    \centering
    \includegraphics[width=0.49\textwidth]{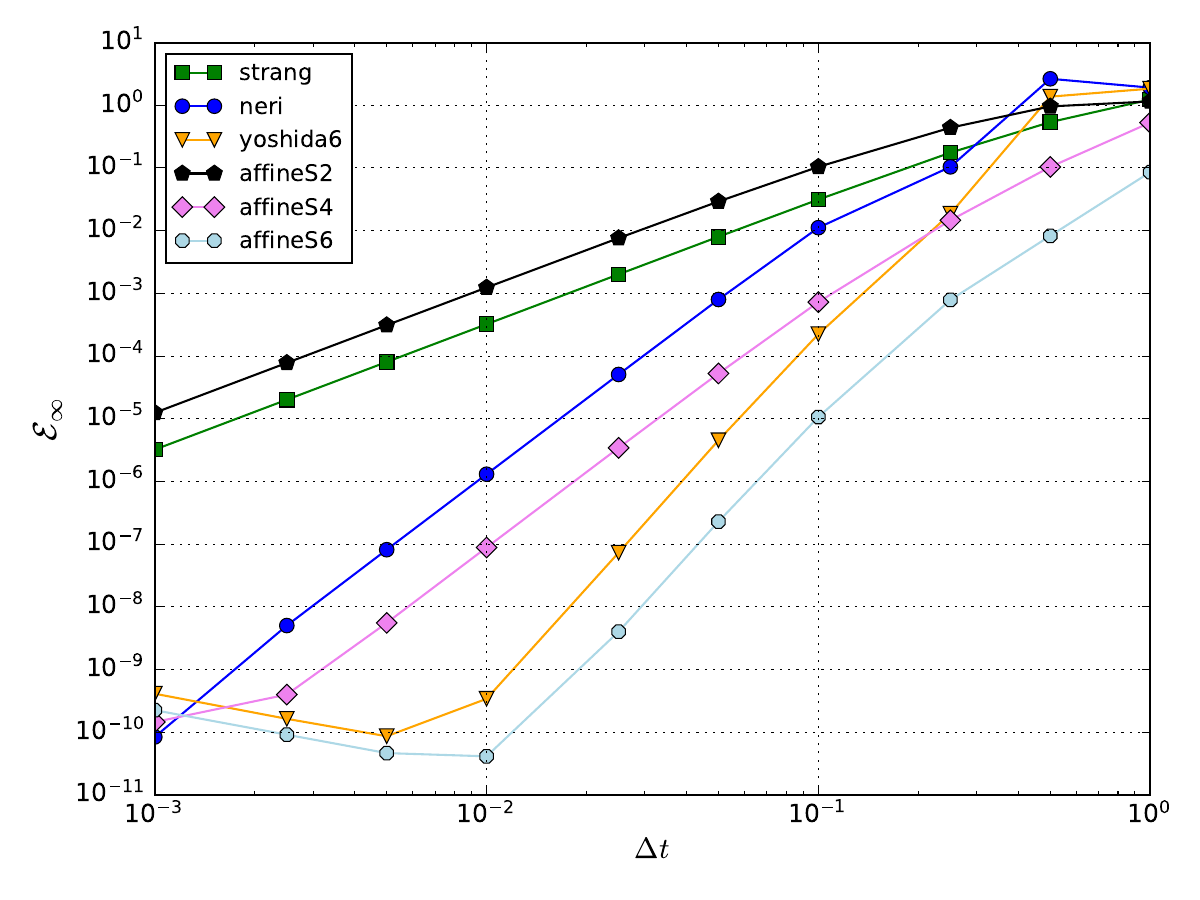}
    \includegraphics[width=0.49\textwidth]{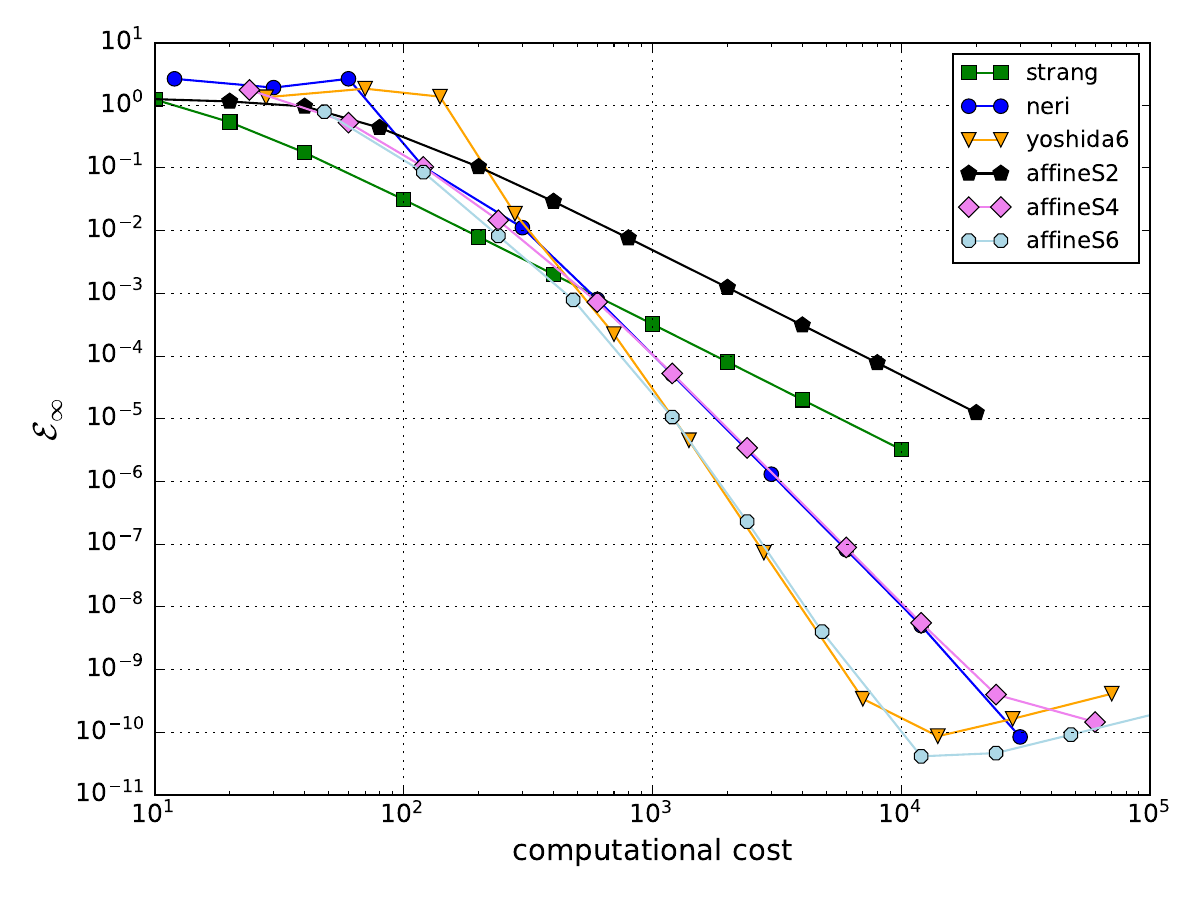} \\
    \includegraphics[width=0.49\textwidth]{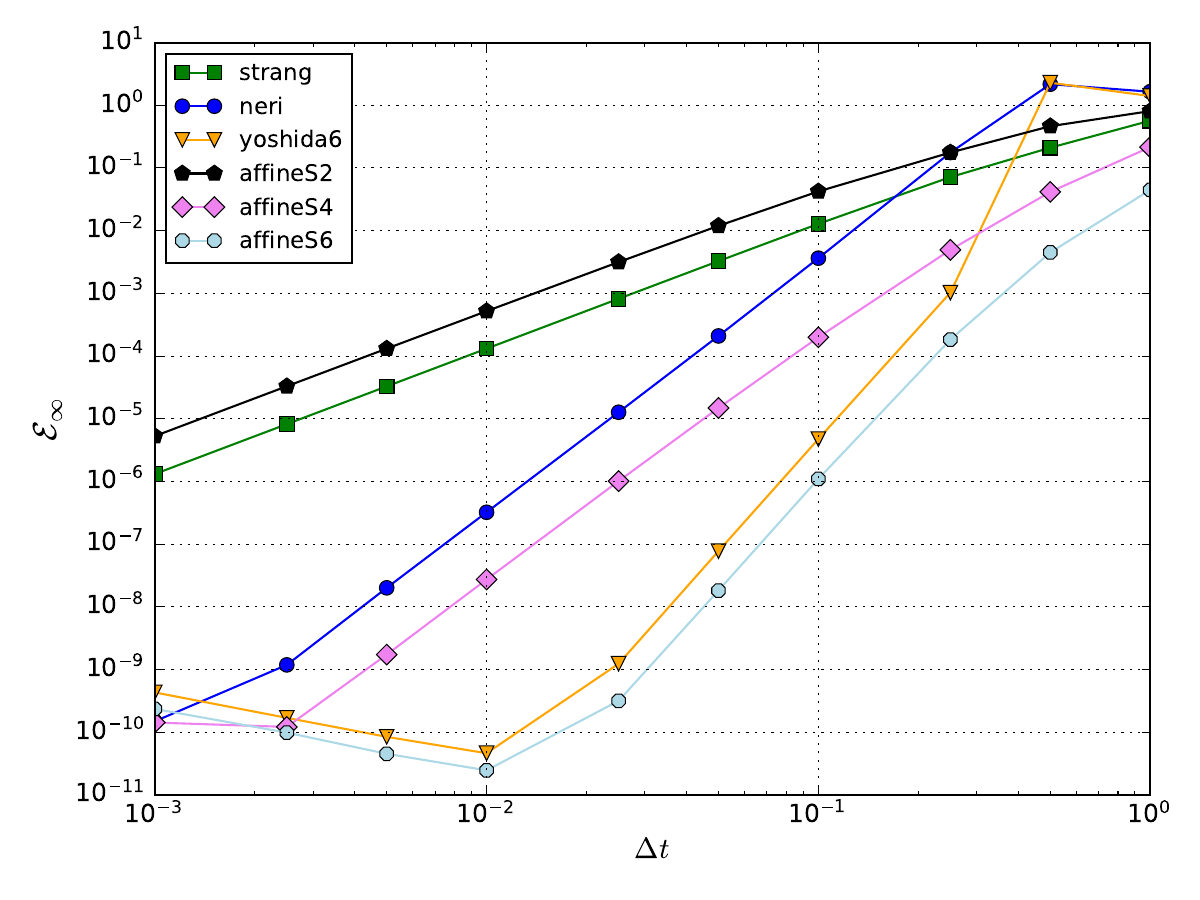}
    \includegraphics[width=0.49\textwidth]{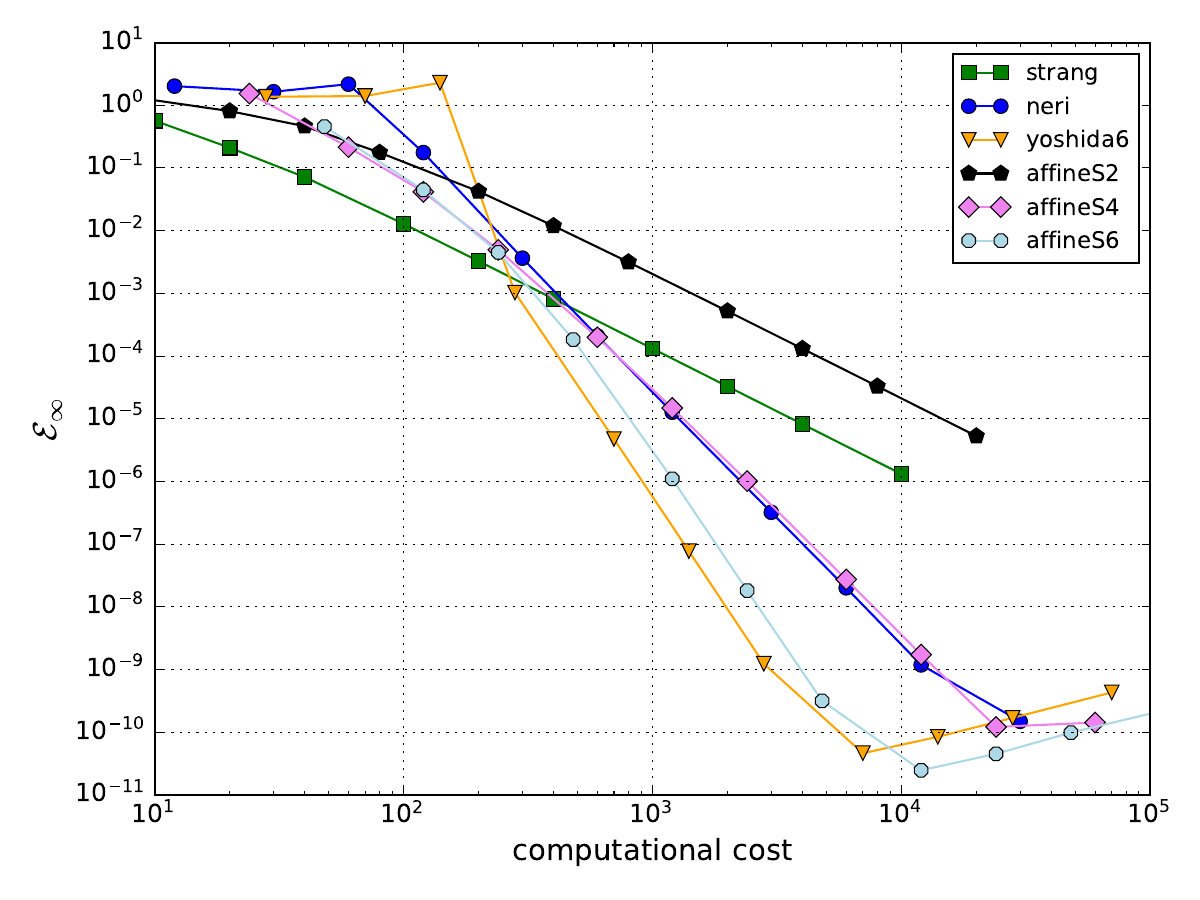}
    \caption{(Left) Absolute error $\mathcal{E}_\infty$ at $t=10$ of the numerical solution of fCGLE5 as a function of step size $\Delta t$.  (Right) Efficiency plot. (Upper) $\alpha=1.8$. (Lower) $\alpha=1.1$. Model parameters are $\beta = 0.1, \gamma = -1.0, \delta = -0.2, \varepsilon = 1.7, \nu = -0.115, \mu = -1.0$.}
    \label{fig:fCGLE5_solution}
\end{figure}
In this subsection we compute the numerical solutions of the fCGLE5  \eqref{eq:fCGLE5} and compare them to a reference solution obtained with an adaptive eighth-order Dormand-Prince 8(5,3) scheme, with tolerances near machine precision \cite{hairer1993solving} (see \ref{sec:propagators} for the implementation details). Given that no explicit analytical expression can be found for the partial propagator $\phi_B$ of the nonlinear cubic-quintic term, this propagator is also approximated by the same numerical integrator.

We take an initial state of gaussian shape $u_0(x) = 1.2 \exp(-x^2/2)$ and model parameters set to $\beta = 0.1, \gamma = -1.0, \delta = -0.2, \varepsilon = 1.7, \nu = -0.115, \mu = -1.0$, for it has been recently reported that in this case the solution converges to a dissipative soliton \cite{qiu2020solitonfCGLE}. 

In Figure \ref{fig:fCGLE5_solution} we plot the dependence of the absolute error $\mathcal{E}_\infty$ at $t=10$ on the time step $\Delta t$ (left panel) and on the computational cost (right panel), in double logarithmic scale, for values of the Lévy index $\alpha=1.8$ (upper panels) and $\alpha=1.1$ (lower panels), respectively. Regardless of the value of the Lévy index, the figure shows that higher order affine integrators achieve lower errors than composition schemes of the same order for identical time steps, with the exception of the symmetric affine scheme of order two. Particularly, for $\alpha=1.1$ the difference can be as high as one order of magnitude. It is worth noting that for time steps that are relatively long, specifically those with $\Delta t \geq 5\cdot 10^{-1}$, the higher order composition methods display a slightly unstable behavior. This instability is evidenced by the non-monotonic relationship between the error and the decreasing step size. Although not depicted, the instability becomes more pronounced as the time step increases, as a consequence of the negative fractional steps required by these methods and the ill-posedness of the backward evolution driven by the diffusive term.

\begin{figure}[h!]
    \centering
    \includegraphics[width=0.45\textwidth]{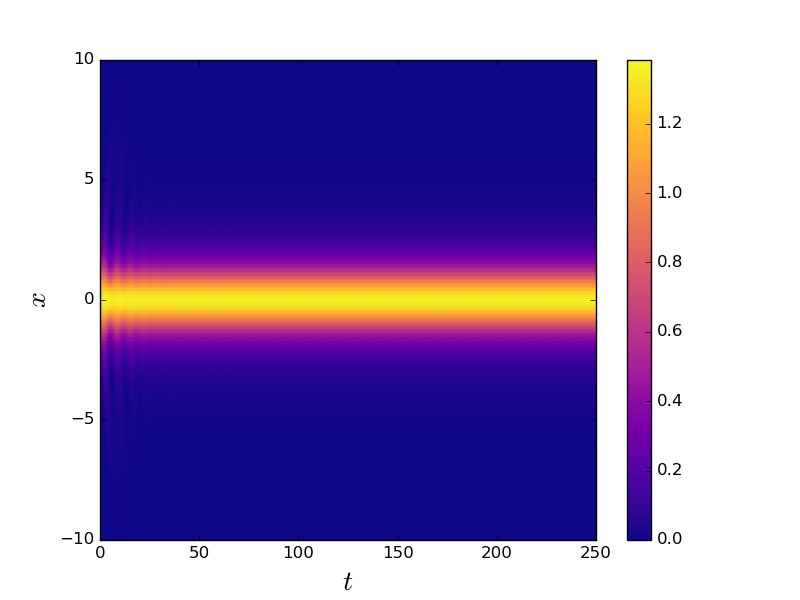} \includegraphics[width=0.54\textwidth]{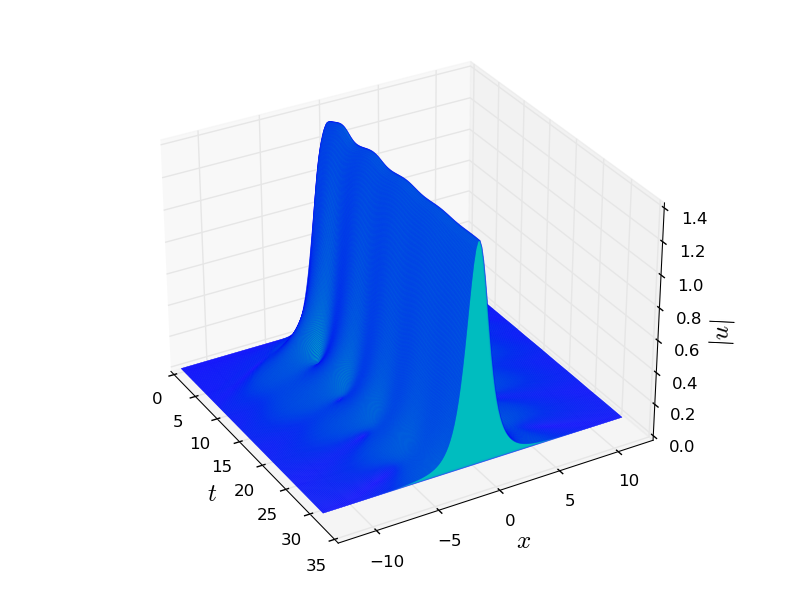}\\
    \includegraphics[width=0.45\textwidth]{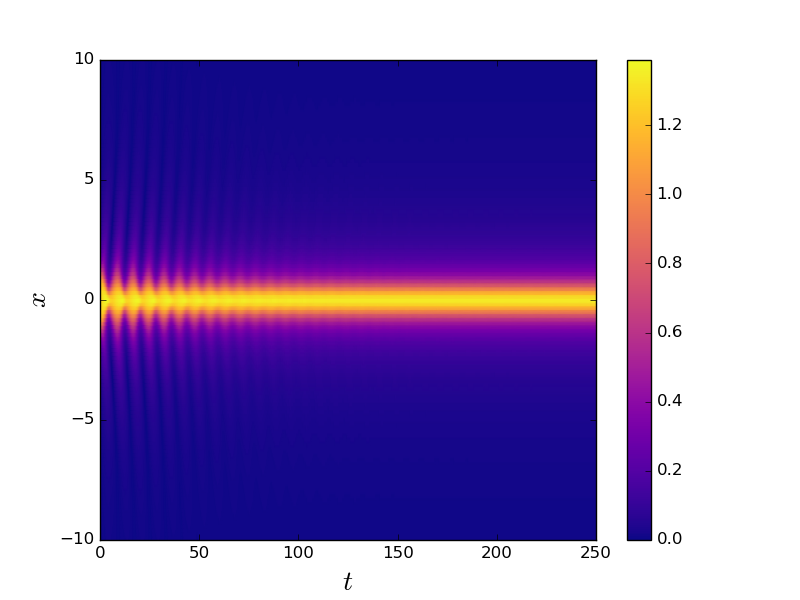} \includegraphics[width=0.54\textwidth]{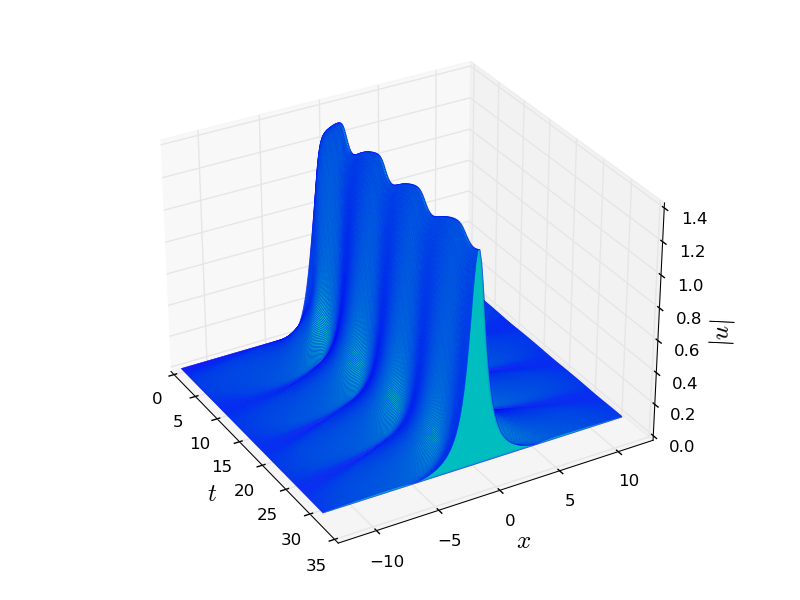}
    \caption{(Left) Evolution of the absolute value of an initially gaussian state up to $t=250$, which converges to a dissipative soliton. (Right) Same evolution up to $t=30$, showing details of the initial dynamics. The upper panels corresponds to Lévy index $\alpha=1.8$, while the lower ones to $\alpha=1.1$. The remaining parameters are $\beta = 0.1, \gamma = -1.0, \delta = -0.2, \varepsilon = 1.7, \nu = -0.115, \mu = -1.0$. The solutions are computed by means of a sixth-order affine splitting with step size $\Delta t=0.025$ using a Hermite pseudo-spectral discretization with $N=300$.}
    \label{fig:fCGLE5_abs_solution_alpha_1.1}
\end{figure}

Complementary, the efficiency plots displayed in the right panels of the the same figure demonstrate that, in general, high-order affine schemes perform similarly to composition schemes of equivalent order, again with the exception of the order-two affine method. For $\alpha=1.1$, the sixth-order composition scheme has better accuracy than the corresponding affine method for the same computational cost. However, this difference diminishes as the Lévy index increases, becoming negligible as $\alpha \to 2$. It is important to recall that, as was mentioned in the introductory remarks, the computational cost of composition schemes is calculated for the most favorable case, and therefore may be underestimated. In light of this fact, affine methods results are highly competitive. On the other hand, the affine schemes have not been subjected to any special optimizations, despite their compatibility with time-adaptivity and parallelization. Taken together, these findings suggest that there is significant potential for further optimization and improvement of affine methods beyond what has been demonstrated in the present study, with the possibility of achieving even better performance. 

The actual shape of the solutions is depicted in Figure \ref{fig:fCGLE5_abs_solution_alpha_1.1}, as calculated by a sixth-order affine Hermite pseudo-spectral method. The absolute value of the solution for $\alpha=1.8$ up to time $t=250$ can be seen in the upper left panel, where the convergence to a stable solution is evident, while the upper right panel shows in more detail the solution dynamics up to $t=30$, revealing a pulsating behavior prior to stabilization. The same information is conveyed in the lower panels for a lower Lévy index $\alpha=1.1$. It can be observed that the solution stabilizes faster in the former case, and converges to a spatially wider soliton after a shorter period of oscillations. Contrarily, the oscillations are wilder for $\alpha=1.1$. This is indeed the general trend as we decrease the Lévy index from $\alpha=2.0$ to $\alpha=1.0$. The qualitative behavior shown complements results given in \cite{qiu2020solitonfCGLE}.

\section{Conclusions}
\label{sec:conclusions}
We have demonstrated that the high-order affine splitting integrators first introduced in \cite{deleo2016affine} are competitive in terms of accuracy and efficiency when compared to numerical methods based on composition splittings, performing better in most situations with respect to various metrics. In particular, we have shown that the sixth-order affine splitting scheme offers remarkably good overall efficiency across a wide range of cases when high accuracy is required, also for models with conserved quantities. Thus, high-order affine splitting integrators can serve as an efficient alternative of composition schemes in both Hamiltonian and irreversible problems, particularly in the last case where composition schemes of order higher than two can become unstable due to negative time steps. 

Moreover, the combination of affine time-splitting schemes with pseudo-spectral space discretizations, based on Fourier and Hermite function expansions, yields a robust and efficient method for solving evolutionary nonlinear partial differential equations in simple geometries. This approach is effective also for non-local fractional models. Indeed, a wide variety of interesting cases have been successfully addressed in the present paper. Through extensive numerical simulations, we have provided evidence of the effectiveness of these methods in solving one-dimensional nonlinear Schrödinger, reaction-diffusion, and complex Ginzburg-Landau equations, including cases with fractional Laplacian operators which are of contemporary theoretical and experimental interest in the fields of nonlinear optics, fractional quantum mechanics and anomalous diffusion. As a byproduct of our simulations, we have found some interesting qualitative relationships between the Lévy index and properties of dissipative solitons of the fCGLE5, which as far as we know were not reported elsewhere. While this paper addressed exclusively one-dimensional problems, the proposed methods admit generalizations to higher dimensions by using tensor product pseudo-spectral bases in problems posed on sufficiently simple geometries. In more complex geometries, affine splittings can also be combined with other space discretizations, like finite-elements, finite-differences and radial basis functions.

Finally, while not explored in this paper, it is worth noting that the possibility of time adaptivity in affine splitting schemes opens a promising avenue for future improvements. Particularly, the pair of fourth and sixth-order affine integrators offers a potentially efficient way to estimate the local error and to adjust the time step in order to achieve even better accuracy with equivalent or lower computational cost. It is true that adaptivity is also available for some pairs of composition schemes, as reported in \cite{auzinger2023splittings, auzinger2017adaptive, auzinger2019adaptive} where schemes up to order $q=6$ are constructed, but the special structure of affine schemes makes them particularly attractive for the straightforward construction of adaptive schemes of any order. Affine splitting schemes are also parallelizable, which could significantly reduce the computational time required for large-scale simulations. The combination of time adaptivity and parallelization presents exciting opportunities for advancing the field of splitting methods for partial differential equation solvers, and we are currently pursuing these options.

%The Appendices part is started with the command \appendix;
%appendix sections are then done as normal sections

\appendix

\section{Pseudo-spectral methods}
\label{sec:pseudo-spectral}
In what follows we give the general approach for the discretization of functions and operators using the pseudo-spectral method \cite{fornberg1996pseudospectral, trefethen2000spectral, boyd2000chebyshev, canuto2007spectral, hesthaven2007spectral, kopriva2009implementing, shen2011spectral}. 

If $\mathcal{H}$ is a Hilbert space and $u:\Omega \to \mathbb{C} \in \mathcal{H}$, we can expand $u$ in an orthonormal basis $\{\varphi_j\}_{j=0}^\infty$ 
\begin{equation}
    u(x) = \sum_{j=0}^{\infty} \widehat{u}_j \varphi_j(x),
    \label{eq:spectral_series}
\end{equation}
where the basis functions satisfy an orthogonality condition with respect to a weight function $w(x)$
\begin{equation}
    \big( \varphi_j, \varphi_k \big) = \int_\Omega w(x) \varphi_j(x) \varphi^\ast_k(x) \mathrm{d}x = \delta_{jk} \lVert \varphi_j\rVert^2, 
\end{equation}
and the spectral coefficients $\widehat{u}_j$ are given by the inner products
\begin{equation}
    \widehat{u}_j = \big(u, \frac{\varphi_j}{\lVert \varphi_j\rVert^2}\big) = \int_\Omega w(x) u(x) \frac{\varphi^\ast_j(x)}{\lVert \varphi_j\rVert^2} \mathrm{d}x.
    \label{eq:spectral_coefficients}
\end{equation}
In the pseudo-spectral method, $u$ is approximated by truncating the expansion \eqref{eq:spectral_series} to $N$ terms 
\begin{equation}
    u(x) \approx \sum_{j=0}^{N-1} \tilde{u}_j \varphi_j(x),
\end{equation}
and approximating the spectral coefficients by means of a suitable numerical quadrature formula, obtaining the {\em pseudo-spectral coefficients} 
\begin{equation}
    \tilde{u}_j := \sum_{n=0}^{N-1} w_n u(x_n) \frac{\varphi^\ast_j(x_n)}{\lVert \varphi_j\rVert^2}.
    \label{eq:forward_discrete_transform}
\end{equation}
The values $w_n,x_n$ are, respectively, the {\em weights} and {\em nodes} of the quadrature formula, defined to  {\em exactly} compute the inner product \eqref{eq:spectral_coefficients} when $u$ is a linear combination of $\{\varphi_j\}_{j=0}^{N-1}$, which also implies the collocation condition
\begin{equation}
    u_n := u(x_n) = \sum_{j=0}^{N-1} \tilde{u}_j \varphi_j(x_n).
    \label{eq:backward_discrete_transform}
\end{equation}
Formulas \eqref{eq:forward_discrete_transform} and \eqref{eq:backward_discrete_transform} are called {\em forward} and {\em backward discrete transforms}, respectively, and constitute the fundamental operations of pseudo-spectral methods. In what follows we denote by $\mathbf{u}$ the vector whose components are $u_n$, and by $\tilde{\mathbf{u}}$ the vector with components $\tilde{u}_j$.

\subsubsection*{Fourier pseudo-spectral method}
For problems posed on a bounded interval $I=[a,b]\subset\mathbb{R}$ with periodic boundary conditions, it is natural to chose a basis of complex exponentials $\varphi_j(x):=\exp\big(\mathrm{i} k_j x \big)$ with $k_j:=2\pi j/\lvert I \rvert, j \in \mathbb{Z}$; and to use a  trapezoidal quadrature formula with $x_n=a + n (b-a)/N$ and $w_n=(b-a)/N, n=0,\dots,N-1$, for approximating the inner products \eqref{eq:spectral_coefficients}. In this case, the pseudo-spectral coefficients are given by the Discrete Fourier Transform
\begin{equation}
    \tilde{u}_j = \frac{1}{N} \sum_{n=0}^{N-1} u_n \mathrm{e}^{-\mathrm{i}\frac{2\pi}{N}jn},\quad j=-\frac{N}{2}, \dots, \frac{N}{2}-1,
\end{equation}
and can be computed by the Fast Fourier Transform (FFT) with complexity $\mathcal{O}(N\log_2 N)$. Conversely, the inverse FFT allows to obtain the values $u_n$ from the pseudo-spectral coefficients $\tilde{u}_j$
\begin{equation}
    u_n = \sum_{j=-N/2}^{N/2-1} \tilde{u}_j \mathrm{e}^{\mathrm{i}\frac{2\pi}{N}jn}.
\end{equation}

\begin{rmk}
    If the problem is posed on $\mathbb{R}$ and $u$ decays rapidly as $\lvert x \rvert \to \infty$, the Fourier pseudo-spectral method can also  be applied by truncating this unbounded domain to a sufficiently large interval $I$, in such a way that the value of the function is negligible (i.e. near machine precision) at the interval boundaries and outside it. An adequate selection of the interval length $\lvert I \rvert$ and number of modes $N$ is crucial to correctly resolve the function in both space and frequency. In the examples considered in this paper, these parameters have been suggested by numerical experimentation.
\end{rmk}

\subsubsection*{Hermite pseudo-spectral method}

For problems posed on unbounded domains $\Omega = \mathbb{R}$, Hermite functions are frequently chosen for the expansion \eqref{eq:spectral_series} \cite{bao2005splitting, boyd2000chebyshev, mao2017fractional, tang1993hermite, ma2005scaling,  thalhammer2009splitting, gauckler2011convergence,  xia2021scalingmoving, chou2023adaptive}. They are given by 
\begin{equation}
    \varphi_j(x) = c_j \mathrm{e}^{-\frac{x^2}{2}} H_j(x), 
\end{equation}
where $H_j(x)$ is the $j$-th order Hermite polynomial and $c_j=\big(2^j j! \sqrt{\pi}\big)^{-1/2}$ are normalization constants such that $\int_\mathbb{R} \varphi_j(x) \varphi_k(x)\mathrm{d}x=\delta_{jk}$ and $\lVert \varphi_j \rVert = 1$. 
These functions can be evaluated stably by the three-term recursion\footnote{In order to overcome the pitfalls related to double-precision arithmetic when evaluating Hermite functions of high-order for $\lvert x \rvert \geq 38.6$, we implemented the strategy given in Appendix A of \cite{bunck2009algorithm}, where also an efficient algorithm is presented.} 
\begin{align}
    \varphi_{0}(x) & =\pi^{-1/4}e^{-\frac{x^{2}}{2}},\quad\varphi_{1}(x)=\pi^{-1/4}\sqrt{2}e^{-\frac{x^{2}}{2}}x,\\
    \varphi_{j}(x) & =\sqrt{\frac{2}{j}}\ x\ \varphi_{j-1}(x)-\sqrt{\frac{j-1}{j}}\ \varphi_{j-2}(x),\quad j\geq2.\label{eq:hermite_function_3_term_recurrence}
\end{align}
An important and useful property of Hermite functions is that they are eigenfunctions of the Fourier Transform
\begin{equation}
\mathcal{F}\big\{\varphi_{j}\big\}=\widehat{\varphi}_{j}=(-\mathrm{i})^{j}\ \varphi_{j},
\label{eq:hermite_eigenfunctions_FT}
\end{equation}
and also of the quantum harmonic oscillator.

With this basis, the quadrature formula \eqref{eq:forward_discrete_transform} requires the nodes $x_n$ to be the zeros of the $N$-th order Hermite polynomial, while the weights are given by
\begin{equation}
    w_n=\frac{1}{N\varphi_{N-1}^{2}(x_{n})},\qquad 0\leq n\leq N-1.\label{eq:modified_weights}
\end{equation}

In practical problems it is usually useful to translate and scale the spatial coordinate by the affine transformation $x' = sx + c$, defining the modified basis functions $\varphi_j^{s,c}(x')=\varphi_j(\frac{x'-c}{s})$. Then the pseudo-spectral coefficients $\tilde{u}_j$ are given by \eqref{eq:forward_discrete_transform} but the function $u$ has to be evaluated on a grid of {\em modified nodes} $x'_n = s x_n + c$, i.e. $u_n=u(s x_n + c)$. $s$ is the {\em scaling factor} and $c$ is the {\em center} of the modified basis. It is known that an adequate selection of the scaling factor $s$ allows to improve the rate of convergence of the series expansion in some situations \cite{boyd2000chebyshev, tang1993hermite, ma2005scaling}. Also, scaling and translating is required for implementing adaptive schemes \cite{xia2021scalingmoving, chou2023adaptive}. 

Under these conditions, the {\em forward discrete Hermite transform} \eqref{eq:forward_discrete_transform} is a matrix-vector product
\begin{equation}
    \tilde{\mathbf{u}} = \mathbf{H}^{+}\  \mathbf{u},
\end{equation}
with the matrix $\mathbf{H}^{+}$ given by
\begin{equation}
    \mathbf{H}^{+} = \mathbf{\varphi} \cdot \mathrm{diag}(w_0, w_1, \dots, w_{N-1}).
\end{equation}
$\mathbf{\varphi}$ is the matrix whose rows are the values of Hermite functions evaluated on the quadrature nodes $\varphi_{jn}=\varphi_j(x_n)$, $\tilde{\mathbf{u}}$ is the vector of pseudo-spectral coefficients $\tilde{u}_j$, and $\mathbf{u}$ is the vector of approximations $u_n$ of $u$ on the grid $x'_n$. Conversely, {\em the backward discrete Hermite transform} \eqref{eq:backward_discrete_transform} is given by
\begin{equation}
    \mathbf{u} = \mathbf{H}^{-}\  \tilde{\mathbf{u}},
\end{equation}
with $\mathbf{H}^{-}=\mathbf{\varphi}^\mathsf{T}$. Both transforms can be implemented with cost $\mathcal{O}(N^2)$.

\section{Computation of propagators}
\label{sec:propagators}
\subsection*{Linear propagator}
The linear propagator $\phi_A$ solves the linear subproblem \eqref{eq:sub_A}, which after a Fourier Transform reads
\begin{equation}
    \mathrm{i}\partial_t \widehat{u}(k,t) = \widehat{Au}(k,t) = \mathcal{A}(k) \widehat{u}(k,t).
\end{equation}
Formally, we have the solution
\begin{equation}
    \widehat{u}(k, t+\Delta t) = \exp \big(-\mathrm{i} \Delta t \mathcal{A}(k) \big)\, \widehat{u}(k,t)=\phi_A(\Delta t) \widehat{u}(k,t),
\end{equation} 
where the symbol is $\mathcal{A}(k)=\frac{1}{2}\lvert k \rvert^\alpha$ for the fNLSE, $\mathcal{A}(k)=-\mathrm{i}(\alpha k^2 - \beta)$ for Fisher's equation, and $\mathcal{A}(k)=(1/2 -\mathrm{i}\beta)\lvert k \rvert^\alpha + \mathrm{i} \delta$ for the fCGLE (notice that the interpretation of parameters $\alpha, \beta$ and $\delta$ depends on each specific model). Particularly, $\mathcal{A}(k)$ is real for the fNLSE3 so we have $\int_\mathbb{R} \lvert \widehat{u}(k, t+\Delta t) \rvert^2 \mathrm{d}k = \int_\mathbb{R} \lvert \widehat{u}(k, t) \rvert^2 \mathrm{d}k$, and by the unitarity of the FT the mass is conserved under the linear evolution (this fact is essential for the mass-preservation property of composition schemes). The numerical computation of the linear propagator $\phi_A$ is basis-dependent and in what follows we describe it using Fourier and Hermite methods, respectively. 

\subsubsection*{Fourier method}
Using the Fourier pseudo-spectral discretization, the time evolution of pseudo-spectral coefficients is approximated by
\begin{equation}
    \tilde{u}_j(t+\Delta t) = \exp \big(-\mathrm{i} \mathcal{A}(k_j) \Delta t\big)\, \tilde{u}_j(t),\quad j=-\frac{N}{2}, \dots, \frac{N}{2}-1.
\end{equation}
In vector form
\begin{equation}
    \tilde{\mathbf{u}}(t+\Delta t) = \mathbf{exp}\big(-\mathrm{i}\Delta t\mathcal{A}(\mathbf{k})\big) \odot \tilde{\mathbf{u}}(t),
    \label{eq:fourier_linear_prop}
\end{equation}
where $\odot$ denotes element-wise multiplication of vectors, $\mathbf{exp}$ is the vectorized exponential function and $\mathcal{A}(\mathbf{k})$ is the vector whose components are the values of the symbol $\mathcal{A}$ at wave-numbers $k_j$. The approximation $\mathbf{\Phi}_A$ of the linear propagator $\phi_A$  by the Fourier pseudo-spectral method is thus given by
\begin{equation}
    \mathbf{\Phi}_A(\Delta t) \mathbf{u}(t) := \mathbf{u}(t+\Delta t) = \mathfrak{F}^{-1}_N\Big( \mathbf{exp}\big(-\mathrm{i}\Delta t\mathcal{A}(\mathbf{k})\big) \odot \mathfrak{F}_N \mathbf{u}(t) \Big),
\end{equation}
where $\mathfrak{F}_N$ denotes the Discrete Fourier Transform for $N$ points.

\subsubsection*{Hermite method}
For the Hermite method, we approximate $u$ by its projection on the linear subspace spanned by the first $N$ Hermite functions $\{\varphi_j\}_{j=0}^{N-1}$, equation \eqref{eq:spectral_series}. The action of the linear operator $A$ on this subspace is completely determined by the matrix of inner products
\begin{equation}
    \widehat{A}_{mn} = \big(A \varphi_m, \varphi_n\big) = \big(\widehat{A\varphi}_n, \widehat{\varphi}_m\big) = (-\mathrm{i})^{n-m} \int_\mathbb{R} \mathcal{A}(k)\varphi_m(k)  \varphi_n (k) \mathrm{d}k,
    \qquad 0\leq m,n \leq N-1,
    \label{eq:A_matrix_elements}
\end{equation}
where we used the definition of $A$ as a Fourier multiplier \eqref{eq:A_symbol}, the unitarity of the FT \eqref{eq:unitaryFT} and the eigenfunction property \eqref{eq:hermite_eigenfunctions_FT} of Hermite functions. Given that is generally not possible to exactly calculate the integral \eqref{eq:A_matrix_elements}, we resort to a Gauss-Hermite quadrature with $N'$ points
\begin{equation}
    \widehat{A}_{mn} \approx \tilde{A}_{mn} :=(-\mathrm{i})^{n-m} \sum_{j=0}^{N'-1}  w_j \mathcal{A}(k_j) \varphi_m(k_j)  \varphi_n(k_j),  \qquad 0\leq m,n \leq N-1.
    \label{eq:A_matrix_quadrature}
\end{equation}
The quadrature is  {\em exact} if $\exp(k^{2})\mathcal{A}(k)\varphi_{p}(k)\varphi_{q}(k)=\mathcal{A}(k)H_{p}(k)H_{q}(k)$ is a polynomial of order $2N'-1$ or less, namely if $\mathcal{A}(k)$ is a polynomial of order not higher than $2(N'-N)+1$. This condition is fulfilled in all the non-fractional models considered in this paper by selecting $N'=N+1$, for the symbol is a complex polynomial of degree not higher than two. For models with fractional Laplacians, more quadrature points are needed in order to give an accurate Hermite pseudo-spectral representation of the linear operator. However, this is not particularly troublesome since the matrix must be calculated only once and this can be done with any required precision without an essential penalty on performance. 
The pseudo-spectral coefficients of $Au$ in the Hermite basis are then
\begin{equation}
    \big(\tilde{Au}\big)_m = \sum_{n=0}^{N-1} \tilde{A}_{mn} \tilde{u}_n,
\end{equation}
and the Hermite pseudo-spectral discretization of $A$ thus transforms the infinite-dimensional partial linear subproblem \eqref{eq:sub_A} into a finite-dimensional system of ODE
\begin{equation}
    \mathrm{i} \partial_t \tilde{u}_m = \sum_{n=0}^{N-1} \tilde{A}_{mn} \tilde{u}_n, \qquad m=0, \dots, N-1,
    \label{eq:linear_evolution}
\end{equation}
that in matrix form reads
\begin{equation}
    \mathrm{i} \partial_t \tilde{\mathbf{u}} = \tilde{\mathbf{A}} \tilde{\mathbf{u}}.
\end{equation}
Its solution is
\begin{equation}
    \tilde{\mathbf{u}}(t+\Delta t) = \exp(-\mathrm{i}\Delta t \tilde{\mathbf{A}}) \tilde{\mathbf{u}}(t).
\end{equation}
 The matrix exponential is computed via the \texttt{expm} function in \texttt{scipy} \cite{virtanen2020scipy}. To obtain the discrete approximation $\mathbf{\Phi}_A$ of the propagator $\phi_A$ in grid space we apply the forward and backward discrete Hermite transforms 
\begin{equation}
    \mathbf{\Phi}_A(\Delta t) \mathbf{u}(t) := \mathbf{u}(t+\Delta t) = \mathbf{H}^{-} \exp(-\mathrm{i}\Delta t \tilde{\mathbf{A}})  \mathbf{H}^{+} \mathbf{u}(t),
\end{equation}
and thus the action of the discrete propagator $\mathbf{\Phi}_A(\Delta t)$ on the discretization of $u$ obtained by the Hermite pseudo-spectral method can be computed as a matrix-vector product with cost $\mathcal{O}(N^2)$.  
\subsection*{Nonlinear propagator}
For all the models considered, the nonlinear propagator $\phi_B$ is diagonal in space. In what follows, we give analytical solutions for the nonlinear propagators of all models investigated in this paper, which can be numerically computed by means of an element-wise function evaluation on the grid points $x_j$ with cost $\mathcal{O}(N)$.

The nonlinear sub-problem for the cubic nonlinearity of the fNLS3 and CGLE3 models is
\begin{equation}
    \mathrm{i}\partial_t u(x,t) = \big(\gamma +\mathrm{i} \epsilon\big) \lvert u(x,t) \rvert ^2 u(x,t), 
\end{equation}
and has the solution
\begin{equation}
    u(x,t) = u(x,0) \exp\Big((\epsilon -\mathrm{i}\gamma)\int_0^t \lvert u(x,t)\rvert^2 \mathrm{d}t\Big),
\end{equation}
which can be obtained exactly by taking into account that
\begin{equation}
    \frac{\partial}{\partial t} \lvert u \rvert^2 =  2 \Re{\frac{\partial u}{\partial t} u^\ast} = 2 \Re{-\mathrm{i} \big(\gamma +\mathrm{i} \epsilon\big)  \lvert u \rvert ^4}  = 2\epsilon\lvert u \rvert^4.
\end{equation}
In particular, for the fNLSE3, $\epsilon=0$ and $\gamma=\pm1$, and thus
\begin{equation}
    u(x,t) = u(x,0) \exp\big(\pm\mathrm{i}\lvert u(x,0)\rvert^2 t\big) =: \phi_B(t) \big(u(x,0)\big),
\end{equation}
so the nonlinear evolution conserves the mass $M=\int_\mathbb{R}\lvert u \rvert^2 \mathrm{d}x$.
For the CGLE3, in general $\epsilon \neq 0$ and 
\begin{equation}
\lvert u(x,t) \rvert^2 = \frac{\lvert u(x,0) \rvert^2}{1-2\epsilon\lvert u(x,0) \rvert^2 t}, 
\label{eq:modulus_evol_cubic_nonlin}
\end{equation}
so the exact solution for the nonlinear subproblem is
\begin{equation}
    u(x,t) = u(x,0) \exp\Big(-\frac{1}{2}\big(1-\mathrm{i}\frac{\gamma}{\epsilon})\ln \big(1-2\epsilon \lvert u(x,0) \rvert^2 t\big)\Big) = \phi_B(t) \big(u(x,0)\big).
\end{equation}
From \eqref{eq:modulus_evol_cubic_nonlin} it is evident that in this case necessarily
\begin{equation}
    t < \frac{1}{2\epsilon \lvert u(x,0) \rvert^2},
\end{equation}
otherwise the solution blows up. This condition imposes an upper bound on the step size for the nonlinear propagator, that depends on $\epsilon$ and on the maximum modulus of the solution.

For a quintic nonlinearity, it can be found an exact solution in an analogous way (although this case was not explored in the present paper).

For the cubic-quintic nonlinearity of the fCGLE5 model, it is not possible to find an \emph{explicit} exact solution, so we resort to a numerical propagator computed using the adaptive Dormand-Prince 8(5,3) scheme implemented in \texttt{scipy} \cite{virtanen2020scipy}. 

For Fisher's equation, the nonlinear subproblem is $\mathrm{i}\partial_t u(x,t) = - \beta u^2(x,t)$, so the nonlinear propagator $\phi_B$ is given by the analytical expression
\begin{equation}
    \phi_B(\Delta t) u(x, t) := u(x,t+\Delta t) = \frac{u(x,t)}{1+u(x,t)\beta\Delta t}.
\end{equation}

\bibliographystyle{elsarticle-num} 
\bibliography{main.bib}

\begin{thebibliography}{10}
\expandafter\ifx\csname url\endcsname\relax
  \def\url#1{\texttt{#1}}\fi
\expandafter\ifx\csname urlprefix\endcsname\relax\def\urlprefix{URL }\fi
\expandafter\ifx\csname href\endcsname\relax
  \def\href#1#2{#2} \def\path#1{#1}\fi

\bibitem{sulem1999nls}
C.~Sulem, P.-L. Sulem, The Nonlinear {Schrödinger} Equation – Self-focusing
  and Wave Collapse, Springer, 1999.

\bibitem{scott2004encyclopedia}
A.~Scott, Encyclopedia of Nonlinear Science, Routledge, 2005.
\newblock \href {https://doi.org/10.4324/9780203647417}
  {\path{doi:10.4324/9780203647417}}.

\bibitem{ablowitz2011nonlinearwaves}
M.~J. Ablowitz, Nonlinear Dispersive Waves: Asymptotic Analysis and Solitons,
  Cambridge Texts in Applied Mathematics, Cambridge University Press, 2011.
\newblock \href {https://doi.org/10.1017/CBO9780511998324}
  {\path{doi:10.1017/CBO9780511998324}}.

\bibitem{agrawal2013nonlinearfiber}
G.~Agrawal, Nonlinear Fiber Optics, 5th Edition, Academic Press, 2013.

\bibitem{Bao2003GPE}
W.~Bao, D.~Jaksch, P.~A. Markowich, Numerical solution of the
  {Gross-Pitaevskii} equation for {Bose-Einstein} condensation, Journal of
  Computational Physics 187~(1) (2003) 318 – 342.
\newblock \href {https://doi.org/10.1016/S0021-9991(03)00102-5}
  {\path{doi:10.1016/S0021-9991(03)00102-5}}.

\bibitem{bao2005splitting}
W.~Bao, J.~Shen, A fourth-order time-splitting {Laguerre-Hermite}
  pseudospectral method for {Bose-Einstein} condensates, SIAM J. Sci. Comp.
  26~(6) (2005) 2010–2028.

\bibitem{aranson2002review}
I.~S. Aranson, L.~Kramer, The world of the complex {Ginzburg-Landau} equation,
  Rev. Mod. Phys. 74 (2002) 99--143.
\newblock \href {https://doi.org/10.1103/RevModPhys.74.99}
  {\path{doi:10.1103/RevModPhys.74.99}}.

\bibitem{garcia-morales2012CGLE}
V.~García-Morales, K.~Krischer, The complex {Ginzburg–Landau} equation: an
  introduction, Contemporary Physics 53~(2) (2012) 79--95.
\newblock \href {https://doi.org/10.1080/00107514.2011.642554}
  {\path{doi:10.1080/00107514.2011.642554}}.

\bibitem{mclahan2002splitting}
R.~I. McLachlan, G.~Reinout, W.~Quispel, Splitting methods, Acta Numerica 11
  (2002) 341–434.

\bibitem{hairer2006geometric}
E.~Hairer, C.~Lubich, G.~Wanner, Geometric Numerical Integration:
  Structure-Preserving Algorithms for Ordinary Differential Equations, Vol.~31,
  Springer Science \& Business Media, 2006.

\bibitem{blanes2008splitting}
S.~Blanes, F.~Casas, A.~Murua, Splitting and composition methods in the
  numerical integration of differential equations, SeMA Journal: Bolet{\'\i}n
  de la Sociedad Espa{\~n}ola de Matem{\'a}tica Aplicada 45 (2008) 89--146.

\bibitem{holden2010splitting}
H.~Holden, K.~H. Karlsen, K.-A. Lie, Splitting methods for partial differential
  equations with rough solutions: Analysis and MATLAB programs, Vol.~11,
  European Mathematical Society, 2010.

\bibitem{borgna2015splitting}
J.~P. Borgna, M.~De~Leo, D.~Rial, C.~S\'anchez de~la Vega, General splitting
  methods for abstract semilinear evolution equations, Commun. Math. Sci.
  13~(1) (2015) 83–101.

\bibitem{glowinski2017splitting}
R.~Glowinski, S.~J. Osher, W.~Yin, Splitting methods in communication, imaging,
  science, and engineering, Springer, 2017.

\bibitem{blanes2017concise}
S.~Blanes, F.~Casas, A concise introduction to geometric numerical integration,
  CRC Press, 2017.

\bibitem{fornberg1996pseudospectral}
B.~Fornberg, A practical guide to pseudospectral methods, Cambridge University
  Press, 1996.

\bibitem{trefethen2000spectral}
L.~N. Trefethen, Spectral Methods in {MATLAB}, Society for Industrial and
  Applied Mathematics, 2000.
\newblock \href {https://doi.org/10.1137/1.9780898719598}
  {\path{doi:10.1137/1.9780898719598}}.

\bibitem{boyd2000chebyshev}
J.~P. Boyd, Chebyshev and {Fourier} Spectral Methods, 2nd Edition, Dover, 2000.

\bibitem{canuto2007spectral}
C.~Canuto, M.~Y. Hussaini, A.~Quarteroni, T.~A. Zang, Spectral methods:
  fundamentals in single domains, Springer Science \& Business Media, 2007.

\bibitem{hesthaven2007spectral}
J.~S. Hesthaven, S.~Gottlieb, D.~Gottlieb, Spectral Methods for Time-Dependent
  Problems, Cambridge University Press., 2007.

\bibitem{kopriva2009implementing}
D.~A. Kopriva, Implementing spectral methods for partial differential
  equations: Algorithms for scientists and engineers, Springer Science \&
  Business Media, 2009.

\bibitem{shen2011spectral}
J.~Shen, T.~Tang, L.-L. Wang, Spectral Methods – Algorithms, Analysis and
  Applications, Springer, 2011.

\bibitem{taha1984analytical}
T.~R. Taha, M.~I. Ablowitz, Analytical and numerical aspects of certain
  nonlinear evolution equations. {II.} {Numerical}, nonlinear {Schr{\"o}dinger}
  equation, Journal of computational physics 55~(2) (1984) 203--230.

\bibitem{weideman1986split}
J.~Weideman, B.~M. Herbst, Split-step methods for the solution of the nonlinear
  {Schr{\"o}dinger} equation, SIAM Journal on Numerical Analysis 23~(3) (1986)
  485--507.

\bibitem{pathria1990pseudo}
D.~Pathria, J.~L. Morris, Pseudo-spectral solution of nonlinear
  {Schr{\"o}dinger} equations, Journal of Computational Physics 87~(1) (1990)
  108--125.

\bibitem{muslu2005higher}
G.~M. Muslu, H.~Erbay, Higher-order split-step {Fourier} schemes for the
  generalized nonlinear {Schr{\"o}dinger} equation, Mathematics and Computers
  in Simulation 67~(6) (2005) 581--595.

\bibitem{antoine2013computational}
X.~Antoine, W.~Bao, C.~Besse, Computational methods for the dynamics of the
  nonlinear {Schr{\"o}dinger/Gross--Pitaevskii} equations, Computer Physics
  Communications 184~(12) (2013) 2621--2633.

\bibitem{akhmediev1996singularities}
N.~N. Akhmediev, V.~V. Afanasjev, J.~M. Soto-Crespo, Singularities and special
  soliton solutions of the cubic-quintic complex {Ginzburg-Landau} equation,
  Phys. Rev. E 53 (1996) 1190--1201.
\newblock \href {https://doi.org/10.1103/PhysRevE.53.1190}
  {\path{doi:10.1103/PhysRevE.53.1190}}.

\bibitem{akhmediev2001pulsating}
N.~Akhmediev, J.~M. Soto-Crespo, G.~Town, Pulsating solitons, chaotic solitons,
  period doubling, and pulse coexistence in mode-locked lasers: Complex
  {Ginzburg-Landau} equation approach, Phys. Rev. E 63 (2001) 056602.
\newblock \href {https://doi.org/10.1103/PhysRevE.63.056602}
  {\path{doi:10.1103/PhysRevE.63.056602}}.

\bibitem{akhmediev2001solitons}
N.~Akhmediev, A.~Ankiewicz, Solitons of the complex {Ginzburg--Landau}
  equation, in: S.~Trillo, W.~Torruellas (Eds.), Spatial Solitons, Springer
  Berlin Heidelberg, 2001, pp. 311--341.

\bibitem{akhmediev2005dissipative}
N.~Akhmediev, A.~Ankiewicz, Dissipative Solitons, Lecture Notes in Physics,
  Springer Berlin Heidelberg, 2005.

\bibitem{akhmediev2008dissipative}
N.~Akhmediev, A.~Ankiewicz, Dissipative solitons: from optics to biology and
  medicine, Vol. 751, Springer Science \& Business Media, 2008.

\bibitem{ferreira2022dissipative}
M.~F. Ferreira, Dissipative Optical Solitons, Vol. 238, Springer Nature, 2022.

\bibitem{grelu2012dissipative}
P.~Grelu, N.~Akhmediev, Dissipative solitons for mode-locked lasers, Nature
  photonics 6~(2) (2012) 84--92.

\bibitem{laskin2000fractional}
N.~Laskin, Fractional quantum mechanics, Phys. Rev. E 62 (2000) 3135--3145.
\newblock \href {https://doi.org/10.1103/PhysRevE.62.3135}
  {\path{doi:10.1103/PhysRevE.62.3135}}.

\bibitem{weitzner2003fractional}
H.~Weitzner, G.~Zaslavsky, Some applications of fractional equations,
  Communications in Nonlinear Science and Numerical Simulation 8~(3) (2003)
  273--281.
\newblock \href {https://doi.org/10.1016/S1007-5704(03)00049-2}
  {\path{doi:10.1016/S1007-5704(03)00049-2}}.

\bibitem{guo2012standingfNLS}
B.~Guo, D.~Huang, Existence and stability of standing waves for nonlinear
  fractional {Schrödinger} equations, Journal of Mathematical Physics 53~(8)
  (2012) 083702.
\newblock \href {https://doi.org/10.1063/1.4746806}
  {\path{doi:10.1063/1.4746806}}.

\bibitem{frank2013fractional}
R.~L. Frank, E.~Lenzmann, {Uniqueness of non-linear ground states for
  fractional {Laplacians} in ${\mathbb{R}}$}, Acta Mathematica 210~(2) (2013)
  261 -- 318.
\newblock \href {https://doi.org/10.1007/s11511-013-0095-9}
  {\path{doi:10.1007/s11511-013-0095-9}}.

\bibitem{klein2014fractional}
C.~Klein, C.~Sparber, P.~Markowich, Numerical study of fractional nonlinear
  {Schrödinger} equations, Proc. R. Soc. A 470 (2014) 0364.

\bibitem{longhi2015fractional}
S.~Longhi, Fractional {Schrödinger} equation in optics, Opt. Lett. 40~(6)
  (2015) 1117--1120.
\newblock \href {https://doi.org/10.1364/OL.40.001117}
  {\path{doi:10.1364/OL.40.001117}}.

\bibitem{duo2016spectralfNLS}
S.~Duo, Y.~Zhang, Mass-conservative {Fourier} spectral methods for solving the
  fractional nonlinear {Schrödinger} equation, Computers \& Mathematics with
  Applications 71~(11) (2016) 2257--2271.
\newblock \href {https://doi.org/10.1016/j.camwa.2015.12.042}
  {\path{doi:10.1016/j.camwa.2015.12.042}}.

\bibitem{mao2017fractional}
Z.~Mao, J.~Shen, Hermite spectral methods for fractional {PDEs} in unbounded
  domains, SIAM J. Sci. Comput. 39~(5) (2017) A1928–A1950.

\bibitem{qiu2020solitonfCGLE}
Y.~Qiu, B.~A. Malomed, D.~Mihalache, X.~Zhu, L.~Zhang, Y.~He, Soliton dynamics
  in a fractional complex {Ginzburg-Landau} model, Chaos, Solitons \& Fractals
  131 (2020) 109471.
\newblock \href {https://doi.org/10.1016/j.chaos.2019.109471}
  {\path{doi:10.1016/j.chaos.2019.109471}}.

\bibitem{liu2023experimental}
S.~Liu, Y.~Zhang, B.~A. Malomed, E.~Karimi, Experimental realisations of the
  fractional {Schr{\"o}dinger} equation in the temporal domain, Nature
  Communications 14~(1) (2023) 222.

\bibitem{goldman1996nonreversible}
D.~Goldman, T.~J. Kaper, N-th order operator splitting schemes and
  nonreversible systems, SIAM journal on numerical analysis 33~(1) (1996)
  349--367.

\bibitem{blanes2005necessity}
S.~Blanes, F.~Casas, On the necessity of negative coefficients for operator
  splitting schemes of order higher than two, Applied Numerical Mathematics
  54~(1) (2005) 23--37.

\bibitem{deleo2016affine}
M.~De~Leo, D.~Rial, C.~S{\'a}nchez de~la Vega, High-order time-splitting
  methods for irreversible equations, IMA Journal of Numerical Analysis 36~(4)
  (2016) 1842--1866.

\bibitem{ruth1983splitting}
R.~D. Ruth, A canonical integration technique, IEEE Transactions on Nuclear
  Science 30~(4) (1983) 2669--2671.
\newblock \href {https://doi.org/10.1109/TNS.1983.4332919}
  {\path{doi:10.1109/TNS.1983.4332919}}.

\bibitem{neri1987lie}
F.~Neri, Lie algebras and canonical integration, Department of Physics report,
  University of Maryland (1987).

\bibitem{yoshida1990symplectic}
H.~Yoshida, Construction of higher order symplectic integrators, Physics
  Letters A 150~(5) (1990) 262--268.
\newblock \href {https://doi.org/10.1016/0375-9601(90)90092-3}
  {\path{doi:10.1016/0375-9601(90)90092-3}}.

\bibitem{auzinger2023splittings}
W.~Auzinger, O.~Koch, Coefficients of various splitting methods,
  \url{http://www.asc.tuwien.ac.at/ ~winfried/splitting/} (Accessed April 7,
  2023).

\bibitem{stillfjord2018adaptive}
T.~Stillfjord, Adaptive high-order splitting schemes for large-scale
  differential {Riccati} equations, Numerical Algorithms 78~(4) (2018)
  1129--1151.

\bibitem{harris2020numpy}
C.~R. Harris, K.~J. Millman, S.~J. Van Der~Walt, R.~Gommers, P.~Virtanen,
  D.~Cournapeau, E.~Wieser, J.~Taylor, S.~Berg, N.~J. Smith, et~al., Array
  programming with {NumPy}, Nature 585~(7825) (2020) 357--362.

\bibitem{virtanen2020scipy}
P.~Virtanen, R.~Gommers, T.~E. Oliphant, M.~Haberland, T.~Reddy, D.~Cournapeau,
  E.~Burovski, P.~Peterson, W.~Weckesser, J.~Bright, et~al., Scipy 1.0:
  fundamental algorithms for scientific computing in {Python}, Nature methods
  17~(3) (2020) 261--272.

\bibitem{raviola2023pseudosplit}
L.~Raviola, Code repository for the paper,
  \url{https://github.com/raviola/pseudosplit-paper/} (2023).

\bibitem{tang1993hermite}
T.~Tang, The {Hermite} spectral method for {Gaussian-type} functions, SIAM
  Journal on Scientific Computing 14~(3) (1993) 594--606.
\newblock \href {https://doi.org/10.1137/0914038} {\path{doi:10.1137/0914038}}.

\bibitem{ma2005scaling}
H.~Ma, W.~Sun, T.~Tang, Hermite spectral methods with a time-dependent scaling
  for parabolic equations in unbounded domains, SIAM Journal on Numerical
  Analysis 43~(1) (2005) 58--75.
\newblock \href {https://doi.org/10.1137/S0036142903421278}
  {\path{doi:10.1137/S0036142903421278}}.

\bibitem{tang2018fractional}
T.~Tang, H.~Yuan, T.~Zhou, Hermite spectral collocation methods for fractional
  {PDEs} in unbounded domains, Communications in Computational Physics 24~(4)
  (2018) 1143--1168.
\newblock \href {https://doi.org/10.4208/cicp.2018.hh80.12}
  {\path{doi:10.4208/cicp.2018.hh80.12}}.

\bibitem{Sari2015}
M.~Sari, Fisher's equation, in: B.~Engquist (Ed.), Encyclopedia of Applied and
  Computational Mathematics, Springer Berlin Heidelberg, Berlin, Heidelberg,
  2015, pp. 550--553.
\newblock \href {https://doi.org/10.1007/978-3-540-70529-1_340}
  {\path{doi:10.1007/978-3-540-70529-1_340}}.

\bibitem{hairer1993solving}
E.~Hairer, S.~P. N{\"o}rsett, G.~Wanner, Solving ordinary differential
  equations. 1, Nonstiff problems, Springer-Verlag, 1993.

\bibitem{auzinger2017adaptive}
W.~Auzinger, H.~Hofst{\"a}tter, D.~Ketcheson, O.~Koch, Practical splitting
  methods for the adaptive integration of nonlinear evolution equations. {Part
  I}: {Construction} of optimized schemes and pairs of schemes, BIT Numerical
  Mathematics 57 (2017) 55--74.

\bibitem{auzinger2019adaptive}
W.~Auzinger, I.~Březinová, H.~Hofstätter, O.~Koch, M.~Quell, Practical
  splitting methods for the adaptive integration of nonlinear evolution
  equations. {Part II}: {Comparisons} of local error estimation and
  step-selection strategies for nonlinear {Schr{\"o}dinger} and wave equations,
  Computer Physics Communications 234 (2019) 55--71.
\newblock \href {https://doi.org/10.1016/j.cpc.2018.08.003}
  {\path{doi:10.1016/j.cpc.2018.08.003}}.

\bibitem{thalhammer2009splitting}
M.~Thalhammer, M.~Caliari, C.~Neuhauser, High-order time-splitting {Hermite}
  and {Fourier} spectral methods, Journal of Computational Physics 228 (2009)
  822–832.

\bibitem{gauckler2011convergence}
L.~Gauckler, Convergence of a split-step {Hermite} method for the
  {Gross–Pitaevskii} equation, IMA Journal of Numerical Analysis 31 (2011)
  396--415.

\bibitem{xia2021scalingmoving}
M.~Xia, S.~Shao, T.~Chou, Efficient scaling and moving techniques for spectral
  methods in unbounded domains, SIAM Journal on Scientific Computing 43~(5)
  (2021) A3244--A3268.
\newblock \href {https://doi.org/10.1137/20M1347711}
  {\path{doi:10.1137/20M1347711}}.

\bibitem{chou2023adaptive}
T.~Chou, S.~Shao, M.~Xia, Adaptive {Hermite} spectral methods in unbounded
  domains, Applied Numerical Mathematics 183 (2023) 201--220.
\newblock \href {https://doi.org/10.1016/j.apnum.2022.09.003}
  {\path{doi:10.1016/j.apnum.2022.09.003}}.

\bibitem{bunck2009algorithm}
B.~Bunck, A fast algorithm for evaluation of normalized {Hermite} functions,
  BIT Numer. Math. 49 (2009) 281–295.

\end{thebibliography}


%\begin{thebibliography}{00}



%\end{thebibliography}
\end{document}